\documentclass{amsart}
\usepackage{comment}
\usepackage[foot]{amsaddr}
\usepackage[utf8]{inputenc}
\usepackage[dvipsnames, table]{xcolor}
\usepackage[top=1in, bottom=1.25in, left=1.25in, right=1.25in]{geometry}
\usepackage[hidelinks]{hyperref}
\newcommand{\ms}[1]{{ \color{magenta}{[MS - #1]}}}
\newcommand{\fb}[1]{{ \color{teal}{[FB - #1]}}}
\newcommand{\ti}[1]{{ \color{blue}{[TI - #1]}}}

\input{packages.sty}
\usepackage[hidelinks]{hyperref} 
\usepackage{algpseudocode}
\usepackage{algorithm}
\usepackage{comment}
\usepackage[dvipsnames, table]{xcolor}
\usepackage{float}
\usepackage{tikz}
\usepackage{multirow}
\usepackage{amsthm}
\newtheorem{remark}{Remark}
\newtheorem{theorem}{Theorem}

\begin{document}
\title[Control-Reg-ROMs]{%Control-Reg-ROM
%An Evolve-Filter-Relax Regularization for Feedback Control of the Navier-Stokes Equations
%A 
New Feedback Control %Strategy 
and Adaptive Evolve-Filter-Relax Regularization for the Navier-Stokes Equations {in the Convection-Dominated Regime}
}
\author[M. Strazzullo, F. Ballarin, T. Iliescu, C. Canuto]{Maria Strazzullo$^{1,*}$, Francesco Ballarin$^{2,*}$, Traian Iliescu$^3$, Claudio Canuto$^{1,*}$}
\date{\today}
\address{$^1$ Politecnico di Torino, Department of Mathematical Sciences ``Giuseppe Luigi Lagrange'', Corso Duca degli Abruzzi, 24, 10129, Torino, Italy}
\address{$^2$Department of Mathematics and Physics, Università Cattolica del Sacro Cuore, via Garzetta 48, I-25133 Brescia, Italy}

\address{$^3$Department of Mathematics, Virginia Tech, Blacksburg, VA 24061, USA}
\address{$^*$INdAM research group GNCS member}

\maketitle

\begin{abstract}
We propose, analyze, and investigate numerically a novel feedback control strategy for high Reynolds number flows.
For both the continuous and the discrete (finite element) settings, we prove that the new strategy yields accurate results for high Reynolds numbers that were not covered by current results.
We also show that the new feedback control 
yields more accurate results than the 
{current} control {approaches} in marginally-resolved numerical simulations of a two-dimensional flow past a circular cylinder at Reynolds numbers $Re=1000$.
We note, however, that for realistic control parameters, the stabilizing effect of the new feedback control strategy is not sufficient in the convection-dominated regime.
Our second contribution is the development of 
an adaptive evolve-filter-relax (aEFR) regularization that stabilizes marginally-resolved simulations in the convection-dominated regime and increases the accuracy of the new feedback control in realistic parameter settings.\
For the finite element setting, we prove that the novel feedback control equipped with the new aEFR method yields accurate results for high Reynolds numbers.
Furthermore, our numerical investigation shows that the new strategy yields accurate results for reduced order models that dramatically decrease the size of the feedback control problem.
\end{abstract}

\section{Introduction}
    \label{sec:intro}

Flow control is central in numerous applications~\cite{gunzburger2002perspectives}. 
{The scientific interest in optimal and suboptimal control strategies for {the} Navier-Stokes equations is amply documented in the literature, see, for example, \cite{BanschBenner2015,COCV_2003__9__197_0,Barbu2004,Benner20201653,benner2019optimal,Deckelnick2004297,10.1007/978-3-642-59709-1_11,gunzburger1991analysis,GunzburgerOptimalControl2000,Heinkenschloss2008,doi:10.1137/S036301290241246X,Hinze2000273,raymond2005local,Raymond2006,RAYMOND2007627}.}
The main goal of flow control is to steer the fluid toward a %beneficial 
{desired} configuration.
For laminar flows, flow control strategies (e.g., feedback control) have been successful at both the theoretical and {the} computational levels.
Turbulent flows, however, still pose significant challenges to classical flow control approaches.
One of these challenges is the 
{convection-dominated regime, in which under-resolved or marginally-resolved full order models (FOMs), i.e., computational models obtained with classical numerical discretizations (e.g., the finite element method), yield inaccurate results, usually in the form of numerical oscillations.
Furthermore, although the flow control strategies generally have a stabilizing effect on the underlying numerical simulations, this is often not sufficient to stabilize the simulations for realistic control parameters. 
Another significant challenge for flow control of turbulent flows is the high computational cost of FOMs.  \\
Indeed, resolving all the spatial scales in a turbulent flow, down to the Kolmogorov scale, can require billions of degree of freedom, %.
{which makes the repeated FOM simulation in {an optimal control setting} impractical.}
{Approaches based on optimal control lead to coupled systems to be solved based on adjoint strategies relying on initial and final conditions in time. The optimality system should be solved in space-time domains or by solving matrix equations. Solving optimality systems can be unfeasible for complex problems in fluids, such as convection-dominated or turbulent Navier-Stokes equations. This contribution aims to provide a less expensive technique chosen a priori that stabilizes the solution and steers it toward a desired profile. To this end, we employ a sub-optimal feedback control law that does not increase the computational time needed for the simulations. We refer the reader to \cite{novo2,novo1,GUNZBURGER2000803,hou,Lee20212533} for sub-optimal a priori feedback laws for the Navier-Stokes equations for control and data assimilation. 
}
%high computational cost of full order models (FOMs), i.e., computational models obtained with classical numerical discretizations (e.g., the finite element method).

{
In this paper, we propose strategies for enabling an efficient and accurate \reviewerAA{distributed (i.e., acting on the whole physical domain)} feedback control of convection-dominated flows modeled by the Navier-Stokes equations. 
{Although the solution is not optimal and \reviewerAA{the distributed nature is not suited for realistic applications, the provided feedback strategy is simple to code and converges exponentially fast toward the desired configuration, representing a first step to efficiently stabilize convection-dominated flows}.}
Indeed, for both the continuous and the finite element settings, we prove that the new feedback control yields accurate results for high Reynolds numbers that were not covered by current results.
This significant improvement over the current feedback control strategies is enabled by introducing a new control forcing term that avoids the strict constraints imposed on the Reynolds number at a theoretical level in current approaches~\cite{GUNZBURGER2000803}. 
Furthermore, we show that the new feedback control yields more accurate results than the current control approaches in marginally-resolved numerical simulations.

Despite its improved theoretical and numerical properties, the new feedback control yields inaccurate results in the convection-dominated regime when realistic control parameters are used, i.e., when the control action is weaker. %smaller}.
We note that, although the control can be often seen as a stabilizer since it steers the approximation to the steady state solution~\cite{COCV_2003__9__197_0,Barbu2004, Raymond2006, RAYMOND2007627}, the convection-dominated regime poses additional stability challenges~\cite{Collis2004237}.
To address the numerical instability of the new feedback control, we propose an adaptive regularization, which leverages spatial filtering to increase the stability of the new feedback control.
The novel {\it adaptive evolve-filter-relax (aEFR)} strategy \reviewerA{for control problems} consists of three simple steps:
%{\red{Change EFR to aEFR.}}
(i) Evolve the current velocity to the next step using the standard FOM discretization.
(ii) Filter the intermediate FOM approximation obtained in step (i).
(iii) Relax the filtered intermediate approximation obtained in step (ii). 
The aEFR approximation at the next time step is precisely the relaxed approximation obtained in step (iii).
\reviewerA{In contrast with standard regularized methods, in aEFR,
steps (i)--(iii) are repeated until a good agreement with the desired state is reached. Namely, the regularization is performed when the simulation is not similar to the desired profile, i.e., when spurious oscillations are more likely to occur. Numerical results show that aEFR velocity converges exponentially in time toward the desired velocity profile}.
The aEFR algorithm has two appealing properties:
(a) %{\it Simplicity:} 
aEFR is one of the %(if not the) 
simplest stabilized %strategy.
{strategies}.
Indeed, one can start with a FOM code, add a simple filtering function (described in Section \ref{sec:aEFRsec}), and obtain the aEFR code in a matter of minutes.
(b) %{\it Modularity:} 
The three steps of the aEFR algorithm {(i.e., running the FOM algorithm for one step to get the intermediate solution, filtering this intermediate solution, and then relaxing it)} are highly modular.
For the finite element setting, we prove and show numerically that the novel feedback control equipped with the new aEFR method yields accurate results for high Reynolds numbers.

The new feedback control and aEFR strategies address the numerical instability of classical flow control strategies for convection-dominated flows.
To address the high computational cost of current flow control approaches, we propose reduced order models (ROMs), which 
}
%Reduced order models (ROMs) 
represent an appealing alternative to the expensive FOMs since the ROM dimensions can be orders of magnitude lower than the FOM dimensions.
% \ms{
%The scientific interest in optimal and suboptimal control strategies for Navier-Stokes equations is amply documented in the literature, see for example \cite{BanschBenner2015,COCV_2003__9__197_0,Barbu2004,Benner20201653,benner2019optimal,Deckelnick2004297,10.1007/978-3-642-59709-1_11,gunzburger1991analysis,gunzburger2000analysis,GunzburgerOptimalControl2000,Heinkenschloss2008,doi:10.1137/S036301290241246X,Hinze2000273,raymond2005local,Raymond2006,RAYMOND2007627}
% . 
%Despite the indisputable usefulness of controlled systems in many fields, their application is still limited due to the computational resources needed to solve optimization tasks. %For these motivations, ROMs for control problems represent a broad research topic.
% }
% \ti{Thus, ROMs represent a natural fit for control problems.}
The following is a far from complete list of papers in which ROMs have been used in control of linear and nonlinear systems ~\cite{afs19, afs20, Alla20173091,StrazzulloZuazua,Dede2012, dolgov2022data,falcone2023approximation,karcher_grepl_2014,karcher2018certified,kunisch2008proper,lee2021efficient,negri2015reduced,negri2013reduced,quarteroni2007reduced,Strazzullo2, StrazzulloRB,  Strazzullo3}, and, in particular, control of the Navier-Stokes equations~\cite{10.1007/978-3-319-23413-7_120,alla2020, Benner20201653,Hinze2000273,Pichi20221361,ZakiaMaria}.
Galerkin ROMs (G-ROMs), which leverage data to build the basis in the Galerkin framework, have been %relatively 
{successful} in the numerical simulation of laminar flows.
However, for turbulent flows, under-resolved G-ROMs yield inaccurate results, usually in the form of spurious numerical oscillations.
%To alleviate/eliminate these numerical oscillations, stabilized ROMs~\cite{parish2023residual, farhat} are generally used.
{Just as in the FOM case, to alleviate these numerical oscillations, stabilized ROMs~\cite{girfoglio2021pod,Girfoglio20210,GIRFOGLIO2022105536,GIRFOGLIO2023112127,grimberg2020stability,gunzburger2020leray,kaneko2020towards,parish2023residual,sabetghadam2012alpha,stabile2019reduced,xie2018numerical} are generally used. %\ms{Maybe cite also Michele and Annalisa here?}. \ti{YES!!!! 
 %Absolutely!  And any ROM+stabilization papers we know of, especially for the NSE.}
In this paper, we leverage the aEFR strategy to stabilize the ROM simulations.
Furthermore, we use the resulting aEFR-ROM within the new feedback control setting to ensure both the accuracy and the efficiency of flow control for convection-dominated flows.
We note that, despite their appeal (and the use of regularized FOMs for control~\cite{mallea2021optimal}), to the best of our knowledge, regularized ROMs (such as aEFR-ROM) have never been used for flow control (see, {however,}~\cite{collis2002analysis,Leykekhman20122012,zoccolan2,zoccolan1} for different FOM and ROM stabilization approaches for optimal control).
}
The rest of the paper is organized as follows:
%in 
{In} Section~\ref{sec:efrFOMuncontrolled}, we present the EFR strategy at both the FOM and the ROM levels.
In Section~\ref{sec:problem}, we first present the current feedback control approach and then propose a novel feedback control strategy that targets higher Reynolds number flows.
For both feedback control approaches, we consider both the continuous and the discrete cases.
At the FOM level, %for both strategies, we prove results for their rates of convergence and we emphasize the role played by the Reynolds number in the a priori error bounds.
{we prove that the new feedback control is accurate at higher Reynolds numbers for which the theoretical results for the standard control do not hold.}
%Finally, 
{In Section~\ref{sec:efr-control},} we propose %an adaptive EFR 
the new aEFR strategy {that further increases the numerical stability of the novel feedback control and enables it to accurately approximate convection-dominated flows for realistic, weak control action}.
In Section~\ref{exp1:fafb}, we first numerically investigate the improvement of the new feedback control strategy over the current approach at the FOM level.
Then, we extend the new feedback control strategy to the EFR setting both at the FOM and ROM level. %investigate 
%display the positive EFR effect on the control at both the FOM and the ROM levels.
Furthermore, we show that, using the novel feedback control framework in the numerical simulation of a flow past a cylinder at Reynolds number $1000$, the EFR and aEFR strategies yield more accurate results than the standard noEFR approach. 
In Section~\ref{sec:conclusions}, we draw conclusions and outline possible research directions.
Furthermore, additional numerical results are presented in Appendices \ref{rem:uncontrolled}, \ref{rem:pred}, and \ref{rem:gammapar}.

\section{The Evolve-Filter-Relax (EFR) Strategy for FOM and ROM %the %Uncontrolled System
%{ Full and Reduced Order Model} 
}
\label{sec:efrFOMuncontrolled}
%This section focuses on EFR 
{In this section, we present the EFR strategy}
. For the sake of brevity, we limit the discussion only to the algorithmic viewpoint {and refer the reader to, e.g.,~\cite{Strazzullo20223148,wells2017evolve} for more details}. %The reader is referred to 
%In Appendix \ref{app:EFRapp}, {we provide} %for some 
%numerical evidence in favor of the use of EFR techniques in the full and reduced settings. %We remark that the algorithms follow the notation already used in \cite{Strazzullo20223148}, where:
{In what follows, we use the same acronyms as those used in~\cite{Strazzullo20223148}:}
\begin{itemize}
    \item[$\circ$] noEFR %stays for 
    denotes the {FOM in which} %case when 
    the EFR regularization is not used. %at the FOM level.} 
    %the solution of the system with no EFR regularization at the FOM level.
    \item[$\circ$] EFR %stays for 
    denotes the {FOM in which} %case when} the solution of the system with 
    the EFR regularization is used. %at the FOM level.
    \item[$\circ$] EFR-noEFR %stays for a FOM EFR regularization but no EFR regularization is used at the ROM level,
    denotes the {ROM in which} %case when 
    the EFR regularization is used at the FOM level but not at the ROM level. 
    \item[$\circ$] EFR-EFR %stays for simulations where
    denotes the {ROM in which} %case when 
    the EFR regularization is used both at the FOM and the ROM levels.
\end{itemize}

%Model order reduction aims at building a low-dimensional model to solve simulations more rapidly with respect to the FOM. In section \ref{sec:POD_t}, we briefly describe the Proper Orthogonal Decomposition-Galerkin (POD-Galerkin) approach that is the backbone of the reduction techniques we use in the numerical results. In section \ref{sec:differential-filter} we describe the EFR-EFR algorithm. This approach has been extensively discussed in \cite{Strazzullo20223148}. 
\subsection{{Navier-Stokes Equations}}
In this section, we %introduce 
{present} the incompressible Navier-Stokes equations {(NSE)}, {which are the mathematical model used in theoretical and numerical investigations}. Let $\Omega \subset \mathbb R^{2}$ be the spatial domain. %we are interested in 
{We seek to approximate} the velocity $u(x, t) \doteq u\in \mathbb U \doteq L^2((0,T); H^1_{u_D}(\Omega))$ and pressure $p(x, t) \doteq  p \in \mathbb Q \doteq L^2((0,T); L^2(\Omega))$ \cite{quarteroni2009numerical}:
 % \ti{Should we use $L_{0}^{2}$ to guarantee uniqueness of pressure?} \ms{Do Neumann boundary conditions allow us to take directly $L^2$?} \fb{Yes, they do}
 % \ti{Could you please add a reference?} \ms{I'll give a look to Temam book tomorrow in the library}.
\begin{equation}
\label{eq:NSE}
\begin{cases}
 \displaystyle u_t  - \nu \Delta u + (u \cdot \nabla) u + \nabla p = 0 & \text{in }  \Omega \times (0, T), \\
\nabla \cdot u = 0  & \text{in }  \Omega \times (0, T), \\
u = u_D & \text{on } \Gamma_D \times (0, T), \\
\displaystyle -p n + \nu \frac{{\partial} u}{\partial n} = 0  & \text{on }  \Gamma_N \times (0, T), \\
u(x, 0) = u_0 & \text{in }  \Omega,
\end{cases}
\end{equation}
where $u_0$ is a given initial condition in $\Omega$, $u_D$ is a given boundary condition on $\Gamma_D$, and $\Gamma_D$ and $\Gamma_N$ are the portions of the domain featuring Dirichlet and ``free flow" boundary conditions, respectively, with $\overline{\Gamma_D} \cup \overline{\Gamma_N} = \partial \Omega$, $\Gamma_D \cap \Gamma_N = \emptyset$. %\ms{why $\overline{\Gamma_D}$ and $\overline{\Gamma_N}$? They are already boundary portions.}
%\ti{I don't understand the question, Maria.} \ms{I meant: why do we need the closure symbol over $\Gamma_D$ and $\Gamma_N$ (before they were not overlined). }
%\ti{Oh, I see.  I think the idea is that they can be open intervals, e.g., $(0,1)$ and $(1,2)$.  Their intersection is the empty set, but their union is not $[0,2]$.  You need the union of their closures for that.}
%The Dirichlet condition is given and denoted by $u_D$. 
Moreover, $\nu$ is the kinematic viscosity and $n$ represents the outer normal vector to $\Gamma_N$. The space $H^1_{u_D}(\Omega)$ contains functions in $H^1(\Omega)$ that %verify 
{satisfy} the Dirichlet boundary condition.\\ Denoting by $\overline{U}$ and ${L}$ the characteristic velocity and length scales of the problem at hand, we %can 
define the Reynolds number as 
$Re \doteq {\overline{U}{L}}/{\nu}$.
%The Reynolds number %indicates {determines} the flow regime. Indeed, large 
{Large} Reynolds {numbers determine}  
%denotes 
a convection-dominated regime, where inertial forces dominate the viscous {forces}. %ones.

\subsection{{The EFR %Algorithm
Strategy}}
\label{sec:efr-fom}
%It is well-known that %in 
%c
Convection-dominated regimes %may 
{generally} lead to spurious numerical oscillations in under-resolved or marginally-resolved discretizations, % \ms{cite stuff}.
{i.e., when the spatial resolution is larger than the Kolomogorov lengthscale, which decreases with the Reynolds number~\cite{Fri95,Pop00}.}
The %evolve-filter-relax (EFR) 
{EFR} algorithm is a numerical stabilization strategy that %algorithm 
can alleviate this issue. 
Let us consider $\Delta t$ as time step. Thus, $t_n = n\Delta t$ for $n = 0, \dots, N_T$, and $T = N_T\Delta t$. 
Let us denote the FE velocity $u^h(t) \in \mathbb U^{N_h^{u}}$ %. With similar notations, we consider 
{and the FE} pressure $p^h \in \mathbb Q^{N_h^p}$, %. Namely, 
where ${N_h^{u}}$ and ${N_h^p}$ denote the FE dimension of the two spaces. We call $u^n$ and $p^n$ the FE variables $u^h(t)$ and $p^h(t)$ evaluated at $t^n$. 
%Relying on 
{Using} the implicit Euler discretization in time,
with %$\chi \in (0,1)$ 
the relaxation parameter $\chi \in (0,1)$, the EFR at the time $t^{n+1}$ reads:
%\ti{I think we should use $\widetilde{u}$ for step I and $\overline{u}$ for filtering in step II, which is the traditional notation. Overline is reserved for filtering :)}
\begin{eqnarray}
         &	\text{(I)}& \text{\emph{ Evolve}:} \quad 
\begin{cases}
        	 \displaystyle \frac{{\widetilde{u}}^{n + 1} - u^n}{\Delta t} + (\widetilde{u}^{n+1} \cdot \nabla)    \widetilde{u}^{n+1} - \nu \Delta \widetilde{u}^{n+1} + \nabla p^{n+1} = 0 & \text{in } \Omega , \vspace{1mm}\\
\nabla \cdot \widetilde{u}^{n+1} = 0 & \text{in } \Omega , \vspace{1mm}\\
\widetilde{u}^{n+1} = u_D^{n+1} & \text{on } \Gamma_D , \vspace{1mm}\\
\displaystyle -p^{n+1} n + \nu \frac{\partial {\widetilde{u}^{n+1}}}{\partial n} = 0  & \text{on } \Gamma_N . \\
\end{cases}
            \label{eqn:ef-rom-1}\nonumber \\[0.3cm]
            &	\text{(II)} &\text{\emph{ Filter:}} \quad
\begin{cases} 
        	 -\delta^2 \, \Delta \overline{u}^{n+1} +  \overline{u}^{n+1} = \widetilde{u}^{n+1}& \text{in } \Omega , \vspace{1mm}\\
\overline{u}^{n+1} = u^{n+1}_D &\text{on } \Gamma_D , \vspace{1mm}\\
\displaystyle \frac{\partial \overline{u}^{n+1}}{\partial  n} = 0  & \text{on } \Gamma_N .
\end{cases}
	\label{eqn:ef-rom-2} \nonumber \\[0.3cm]
            &	\text{(III)} &\text{\emph{ Relax:}} \qquad 
        	   u^{n+1}
            = (1 - \chi) \, \widetilde{u}^{n+1}
            + \chi \, \overline{u}^{n+1} \, .
            \label{eqn:ef-rom-3}\nonumber
\end{eqnarray}

{In step (I), the velocity approximation $u^n$ at the current time is evolved to $\widetilde{u}^{n+1}$, which is the intermediate velocity approximation at the new time, $t^{n+1}$.
In step (II), a \emph{differential filter} (DF) 
with filtering radius $\delta$ is used to filter the intermediate velocity approximation, $\widetilde{u}^{n+1}$, obtained in step (I).
The DF alleviates oscillations by leveraging an elliptic operator that eliminates the high frequencies from $\widetilde{u}^{n+1}$.
}
% Step (III) implies that $u$ is a convex combination of the \emph{evolved velocity} $\overline{u}$, obtained at step (I), and the \emph{filtered velocity} $\widetilde{u}$, computed at step (II).
% The filter step (II) involves a \emph{differential filter} (DF) 
% with \emph{filtering radius} $\delta$.
% The DF alleviates oscillations by means of an elliptic operator which dumps the high frequencies from $u^n$. %As a 
In %the last 
step (III), the filtered %solution 
{velocity} is relaxed in order to diminish the diffusion action \cite{ervin2012numerical, fischer2001filter,mullen1999filtering} and increase the accuracy of the simulation: the reader may refer to \cite{bertagna2016deconvolution} and %to 
\cite{ervin2012numerical} for %some experimental and %the 
{numerical and} theoretical results, respectively. {We stress that other filters might be employed {in the EFR algorithm} (see, for example, \cite{girfoglio2021pod} and the references therein). In our approach, we followed the %approach of 
{strategy used in} \cite{Strazzullo20223148}.}

\subsection{The %POD-Galerkin 
{EFR-noEFR} approach}
\label{sec:POD_t}
{In this section, we present the EFR-noEFR approach, i.e., we use the EFR stabilization at %a 
{the} FOM level (as described in Section~\ref{sec:efr-fom}), but at the ROM level we use the standard ROM (without EFR stabilization).}
POD-based ROM algorithms build reduced spaces based on data compression %over 
{of} the \emph{snapshots}, i.e., FOM simulations collected %in 
{at} specific time instances. %\\
%Let us assume to be provided with 
{Consider} two sets of basis functions, for velocity and pressure, respectively, $\{\phi_j\}_{j=1}^{r_u}$ and $\{\psi_j\}_{j=1}^{r{_p}}$, %They 
{which} span the reduced spaces $\mathbb U^{r_u}$ and $\mathbb Q^{r_p}$. %and, consequently, we 
{We} can expand $u_{r_u} \in \mathbb U^{r_u}$ and $p_{r_p} \in \mathbb Q^{r_p}$ as
$ 
	u_{r_u} \doteq {u}^{r_u}(x,t)
	= \sum_{j=1}^r a_j^{u}(t) \phi_j(x)$ {and} $
  p_{r_p} \doteq {p}^{r_p}(x,t)
	= \sum_{j=1}^r a_j^{p}(t) \psi_j(x),
	\label{eqn:g-rom-1}
$
where we denote the reduced coefficients as $\{a_{j}^{u}(t)\}_{j=1}^{r_u}$ and $\{a_{j}^{p}(t)\}_{j=1}^{r_p}$ \cite{noack2011reduced}. %For now, we are assuming the reduced spaces to have both dimension $r$, this is not the case in our setting.
The bases are built %over the manipulation of 
{from} the snapshots  
$\{u_i\}_{i=1}^{N_{\text{max}}} \subseteq \{u^k\}_{k=1}^{N_T}$ and $\{p_i\}_{i=1}^{N_{\text{max}}} \subseteq \{p^k\}_{k=1}^{N_T}$, where ${N_{\text{max}}}$ is the number of snapshots we consider. Here, for the sake of simplicity, we assume that $N_{\text{max}}$ is the same for velocity and pressure, but, in principle, the velocity and pressure snapshots number may not coincide.
%\ti{Why do we assume that the velocity and pressure ROM spaces to have the same dimension?  And why do we assume the same number of snapshots?  Do we really need that?}\\
In our setting, the snapshots are {obtained from a } regularized {model}, since the EFR approach is applied at the FOM level as described in Section \ref{sec:efrFOMuncontrolled}.
%\ti{Mention model consistency?} \ms{Do you mean the fact that in the EFR-noEFR case there is no consistency between the FOM and the ROM level?}
%We introduce 
{Next, we outline} the {proper orthogonal decomposition (POD)} algorithm \cite{ballarin2015supremizer, burkardt2006pod, hesthaven2015certified} {that} we use to build the reduced spaces. 
First, %of all, 
we stress that we %applied 
{apply} \emph{supremizer stabilization} to enrich $\mathbb U^r$ and guarantee the well-posedness of the reduced system \cite{ballarin2015supremizer,rozza2007stability}. %Namely, supremizer 
{The supremizer} stabilization avoids spurious reduced pressure modes. 
We %also 
note that, when 
%When 
dealing with convection dominated setting, other {stabilization strategies} %techniques 
are needed at the reduced level; %, 
we postpone this topic to Section \ref{sec:efrefr}. %and appendix \ref{app:EFRapp}. 
We define the supremizer operator
$S: \mathbb Q^{N_h^p} \rightarrow{{\mathbb U}}^{N_h^{u}}$ as%:
$(S(p), \tau)_{\mathbb {U}} = \left ( p, \nabla \cdot \tau \right )$ for all $\tau \in \mathbb{U}^{N_h^{u}}$. %\\

The enriched reduced velocity space is
%$$
$
{{\mathbb U}}^{r_{us}} \doteq \text{POD}(\{u_i\}_{i=1}^{N_{\text{max}}}, r_u)\oplus \text{POD}(\{S(p_i)\}_{i=1}^{N_{\text{max}}}, r_s).
\label{eqn:definition-rom-velocity-space}
$
For the pressure, we employ the standard POD procedure 
$
\mathbb Q^{r_p} \doteq \text{POD}(\{p_i\}_{i=1}^{N_{p}}, r_p).
$
%\ti{So we are not using the same number of velocity and pressure modes.}
From these processes, we retain $r_u$, $r_s$, and $r_p$ modes for velocity, supremizer, and pressure, {respectively,} where the enriched velocity space has $r_{us} = r_u + r_s$ modes.
%Namely, 
{Thus,} the bases are of the form $\{\phi_j\}_{j=1}^{r_{us}}$ 
and $\{\psi_j\}_{j=1}^{r_p}$, for the velocity and pressure, respectively. After building the bases, for each time instance $t^{n+1}$, %one projects 
{we project} the system %in 
{onto} this low-dimensional space, solving the weak formulation of the NSE: 
\vspace*{-0.01cm}
%NS equations:
\begin{equation}
\begin{cases}
	\displaystyle \int_{\Omega}\frac{u_{r_{us}}^{n+1} - u_{r_{us}}^{n}} {\Delta t} \cdot \phi_i \; dx
	+ \nu a(u_{r_{us}}^{n+1},  \phi_i)
	+ c(u_{r_{us}}^{n+1}; u_{r_{us}}^{n+1} , \phi_i)  - b(\phi_i, p_{r_p}^{n+1})
	= 0, \\[0.5cm]
b( u_{r_{us}}^{n+1},  \psi_j ) = 0,
\end{cases}
\label{eqn:g-rom-2-n+1}
\end{equation}
for all $ i = 1, \ldots, r_{us},$ and $j = 1, \ldots, r_p$, where, for all $u, v \in \mathbb U$ and $p \in \mathbb Q$
%\ti{Are we using the skew-symmetrized } \ms{not at the continuous level}
%\begin{align*}
    $$ 
    a(u,v) = \int_{\Omega} \nabla u : \nabla v \; dx, \quad b(v,p) = \int_{\Omega} p \nabla \cdot v \; dx, \quad  c(w;u,v) = \int_{\Omega} (w \cdot \nabla)u \cdot v\; dx. 
    $$ %& \forall u,v \in \mathbb U 
    % & \forall p \in \mathbb Q \text{ and }\forall u, v \in \mathbb U,\\
     % & b(v,p) = \int_{\Omega} p \nabla \cdot v \; dx & \forall p \in \mathbb Q \text{ and }\forall v \in \mathbb U,\\
     % & c(w;u,v) = \int_{\Omega} (w \cdot \nabla)u \cdot v\; dx & \forall w, u, v \in \mathbb U.
%\end{align*}
\begin{comment}

The algebraic formulation of \eqref{eqn:g-rom-2-n+1} is
\begin{equation}
\begin{cases}
    \displaystyle \frac{1}{\Delta t}\mathtt{M}(\mathtt u^{n+1} - \mathtt u^{n}) + 
    \nu \mathtt K \mathtt u^{n+1} + \mathtt C(\mathtt u^{n+1})\mathtt u^{n+1} - \mathtt {B^T} \mathtt p^{n+1} = 0,\\
    \mathtt B \mathtt u^{n+1} = 0,
    \end{cases}
\end{equation}
where $\mathtt u^{n+1} \in \mathbb R^{r_{us}}$ are the reduced unknowns of the problem, as defined in \eqref{eqn:g-rom-1} and
\begin{equation}
\label{eq:MK}
    \mathtt M_{ij} \doteq
    \displaystyle \int_{\Omega} \phi_{i}\phi_j \; dx, 
    \quad
    \mathtt K_{ij} \doteq a(\phi_{i},  \phi_j), \quad
    \mathtt C(\mathtt u^{n+1})_{ij} \doteq
    c (\mathtt u^{n+1}; \phi_{i} , \phi_j)
    \quad \text{and} \quad 
    \mathtt B_{ij} \doteq b(\phi_{i},  \psi_j).
\end{equation}
\ms{check the formulation, I think that the nonlinear matrix depends on the high-fidelity vector}

\end{comment}
Following the %notation 
{terminology} in \cite{Strazzullo20223148}, we call %this approach, 
{the resulting ROM,}
which is based on standard Galerkin projection onto POD spaces, as EFR-noEFR ROMs, to highlight that EFR is applied %to the snapshots 
{at the FOM level (to generate the snapshots),} 
but no EFR stabilization is performed at the reduced level.

\subsection{The EFR-EFR approach}
	\label{sec:efrefr}
The EFR-EFR strategy employs the EFR algorithm %in the generation of the snapshots and at the ROM stage, too. 
{at both the FOM and the ROM levels.}
While the reduced spaces are built in the same way as in Section \ref{sec:POD_t}, the reduced system %to be solved 
is different.
Indeed, the EFR strategy is applied also %at 
in the reduced setting:
%\ti{I think we should use $\widetilde{u}$ for step I and $\overline{u}$ for filtering in step II, which is the traditional notation.}

\begin{eqnarray*}
         &	\text{(I)}_r& \text{\emph{}} \quad 
\begin{cases}
	\displaystyle \int_{\Omega}\frac{\widetilde{u}_{r_{us}}^{n+1} - u_{r_{us}}^{n}} {\Delta t} \cdot \phi_i \; dx
	+ \nu a({\widetilde{u}}_{r_{us}}^{n+1},  \phi_i)
	+ c({\widetilde{u}}_{r_{us}}^{n+1}; {\widetilde{u}}_{r_{us}}^{n+1} , \phi_i)  - b(\phi_i, p_{r_p}^{n+1})
	= 0, \\[0.5cm]
b( {\widetilde{u}}_{r_{us}}^{n+1},  \psi_j ) = 0,
\end{cases}
            \label{eqn:ef-rom-1-r}\nonumber \\[0.3cm]
            &	\text{(II)}_r &\text{\emph{}} \quad
        	 \delta^2 \, a(\overline{u}^{n+1}_{r_{us}}, \phi_i)  +  \int_{\Omega} \overline{u}^{n+1}_{r_{us}} \cdot \phi_i \; dx = \int_{\Omega} \widetilde{u}^{n+1}_{r_{us}} \cdot \phi_i \; dx,\\
        	 %\text{\TI{Should it be a plus sign in front of $2 \delta^2$?}}\\
	\label{eqn:ef-rom-2-r} \nonumber \\[0.3cm]
            &	\text{(III)}_r &\text{\emph{}} \qquad 
        	   u^{n+1}_{r_{us}}
            = (1 - \chi) \widetilde{u}^{n+1}_{r_{us}}
            + \chi \overline{u}^{n+1}_{r_{us}}.
            \label{eqn:ef-rom-3}\nonumber
\end{eqnarray*}
%As did for velocity and pressure, the 
{The} reduced variables $\overline{u}_{r_{us}}$ and $\widetilde{u}_{r_{us}}$ can be expressed in $\mathbb U^{r}$ as \begin{equation}
	\overline{u}_{r_{us}} \doteq {\overline{u}}_{r_{us}}(x,t)
	= \sum_{j=1}^{r_{us}} a_j^{\overline{u}}(t) \phi_j(x) \quad \text{and} \quad \widetilde{u}_{r_{us}} \doteq {\widetilde{u}}_{r_{us}}(x,t)
	= \sum_{j=1}^{r_{us}} a_j^{\widetilde{u}}(t) \phi_j(x).
	\label{eqn:ws}
\end{equation}
%Namely, 
\begin{comment}
and, at each time instance $t^{n+1}$, the following algebraic system is solved:
\begin{equation}
\begin{cases}
    \displaystyle \frac{1}{\Delta t}\mathtt{M}({\mathtt {\widetilde{u}}}^{n+1} - \mathtt u^n) + 
    \nu \mathtt K {\mathtt {\widetilde{u}}}^{n+1}+ \mathtt C({\mathtt {\widetilde{u}}}^{n+1})\mathtt {\mathtt {\widetilde{u}}}^{n+1} - \mathtt {B^T} \mathtt p^{n+1} = 0,\\
    \mathtt B {\mathtt {\widetilde{u}}}^{n+1} = 0, \\
    \delta^2 \mathtt K {\mathtt {\overline{u}}}^{n+1} + 
    \mathtt M { \mathtt {\overline{u}}}^{n+1} = \mathtt M {\mathtt {\widetilde{u}}}^{n+1},\\
    % \text{\TI{Should it be a plus sign in front of $2 \delta^2$?}}\\
    \mathtt u^{n+1} = (1 -\chi){\mathtt {\widetilde{u}}}^{n+1} + \chi \overline{\mathtt u}^{n+1}.
    \end{cases}
    \label{eqn:EFR-EFR-system}
\end{equation}
The reduced unknown vectors $\mathtt {\widetilde{u}}^{n+1}$ and $\overline{\mathtt u}^{n+1}$ in $\mathbb R^{r_{us}}$ represent the evolved and filtered velocity, respectively. The matrices of system \eqref{eqn:EFR-EFR-system} are defined in Section \ref{sec:POD_t}.
\end{comment}
The use of the EFR-EFR strategy is beneficial %and necessary 
when dealing with large Reynolds numbers, since it alleviates numerical oscillations that %may 
often arise at the ROM level in the convection-dominated regime.
\section{%The 
A new feedback control strategy for %FOM and ROM
high Reynolds numbers}
%Setting {for the Full and Reduced Order %Model 
%{Models}
%}
\label{sec:problem}
%This section describes the NSE under a linear feedback control action at the continuous and FOM level.
{In this section, we propose a novel feedback control strategy for high Reynolds numbers.
In Section~\ref{sec:problem-continuous}, we present the new control strategy for the continuous case, and in Section~\ref{sec:problem-discrete}, we present it for the discrete case.  
}
 %\fb{Myself and TI are thinking about splitting this Section in two parts. Maybe split at 3.3?} 
 %\fb{Stress more that Theorems 2, 4, and related ones are novel}
% \ms{Done}
% \ti{Should we split this Section into two subSections?  I'm thinking 5.1 Continuous Case, and 5.2 Discrete Case.}
% \ms{Yes, It is a good idea}
\subsection{Continuous formulation}
\label{sec:problem-continuous}
Let $u \in \mathbb U$ and $p \in \mathbb Q$ be the state velocity and pressure variables. \reviewerA{Moreover, {let} $f \in L^2((0,T); L^2(\Omega))$ %represents 
{be} a distributed control law to be defined}. The action of $f$ steers $u$ toward a target divergence-free velocity $U$. The three variables verify the weak formulation of the %Navier-Stokes (NS) equations 
{NSE} almost everywhere for $t \in (0,T)$:
\begin{equation}
\label{eq:general_problem}
\begin{cases}
    \left \langle u_t, v\right \rangle + \nu a(u, v) + c(u; u, v) + b(v, p) = \left \langle f, v\right \rangle & \forall v \in H^1_{\Gamma_D}(\Omega),\\
    b(u,q) =0 & \forall q \in L^2(\Omega),
    \end{cases}
\end{equation}
where the forms are the ones introduced in Section \ref{sec:POD_t}.
% \ti{Define bilinear and trilinear forms.}
The problem features Dirichlet and Neumann 
%\ti{So we do use homogeneous BCs!} \ms{my bad, I rephrased}
%\ti{Do we really need homogeneous BCs?} \ms{I realized that the only thing you need is that $u$ and $U$ have the same boundary conditions. In the numerical experiment, this is evident and theoretically it comes from the fact that you have to choose $v = (u - U) \in H^1_{\Gamma_D}(\Omega)$ and this is true if and only if $u_{|_{\Gamma_D}} = U_{|_{\Gamma_D}}$. }
boundary conditions, and an initial velocity $u_0$. %\fb{a pagina 1 lo spazio $H^1_{\Gamma_D}$ era stato definito come quello che soddisfava le condizioni non omogenee (che in effetti è un po' strano e non standard)}. 
%In this context, we aim 
{We leverage the control} 
to reach the following goal:
\begin{equation}
\label{eq:convergence}
    \frac{d}{dt}\norm{u - U}^2_{L^2} \leq 0 \quad \text{a.e. in } (0,T).
\end{equation}
Namely, we want the solution, in time, to %be 
{become} more and more similar to the desired state $U$. This goal \eqref{eq:convergence} is achieved by {using} a linear feedback control in the first variable $u$, i.e., by %means of  
{appropriately choosing} $f = f(u,U)$. %\fb{justify the adjective linear: does it mean that $f$ is linear in the first argument?}.
If \eqref{eq:convergence} %is verified 
{holds} 
for a solution of \eqref{eq:general_problem}, then $(u,p,f)$ is an \emph{admissible solution}.\\ %\\
\reviewerA{From now on, we consider a desired state $U$ that satisfies the following conditions:
%We consider a divergence-free vector field $U$
\begin{itemize}
\item[$\circ$] It is divergence free and belongs to the following space~\cite{GUNZBURGER2000803}
$$
\mathbb U_d = \{ U \; : \; U \in C((0,T), H^1_{u_D}(\Omega) \cap H^2(\Omega) \text{ and } U_t \in C((0,T), H^1(\Omega))\},
$$
\item[$\circ$] features the same Dirichlet boundary condition %of 
{as} the velocity field $u$.
\end{itemize}
If these assumptions are satisfied, $U$ is said to be an \emph{admissible desired state}. We assume $U$ to be admissible for Theorem \ref{theo:gunz1}, Theorem \ref{theo:gunz3}, Theorem \ref{theo:our1}, Theorem \ref{theo:our2} and Theorem \ref{theo:our5}.}
%\ti{Do we need $U_t \in C((0,T), H^1(\Omega))$?  Do we need $H^2$?} \ms{These assumptions are made by Gunzburger too and they are related to the regularity theory for NSE - not sure :) - }
%\ti{I think we should specify the boundary conditions for U since this is important in (12)-(13).  Do we need to associate a pressure, $P$, for that purpose?} \ms{I don't think so. It is related only to the Dirichlet boundary condition}.
%In our setting, the desired state $U$ 
The force {associated to} $U$ is defined as
\begin{equation}
    \label{eq:F}
    F = U_t - \nu \Delta U + (U \cdot \nabla)U.
\end{equation}

%{We stress that we deal with a desired state that asymptotically reaches a steady state. 
%In other words, 
%{Thus,} for every $ \epsilon > 0$ there exists $\hat t > 0$ such that $|U_t| \leq \epsilon $}. 

%\\ 

We aim at building a linear feedback control law to deal with high $Re$ values. % of $Re$. 
First, %of all, 
%\ti{Why to this end?} \ms{My bad, I rephrased}
we recall two inequalities that %come in handy 
{will be used} in what follows \cite{quarteroni2008numerical}:
%\ti{Should we cite Temam~\cite{temam2001navier} for the continuity?} \ms{maybe \cite{quarteroni2008numerical}?}
\begin{itemize}
    \item[$\circ$]\emph{Poincar\'e inequality}: $\norm{u}_{L^2}^2 \leq C_P \norm{\nabla u}_{L^2}^2$, for all $u \in H^1_0(\Omega)$,
    \item[$\circ$]\emph{continuity of} $c(\cdot; \cdot, \cdot)$: $c(u;w,v) \leq K_0 \norm{\nabla u}_{L^2}\norm{\nabla w}_{L^2}\norm{\nabla  v}_{L^2}$ for all $u,w,v \in H^1(\Omega)$.
\end{itemize}
We %refer to 
{denote} the inverse of the Poincar\'e constant %as 
with $C_0 = {C_P}^{-1}$ and %to $K_0$ as 
the continuity constant of $c(\cdot; \cdot, \cdot)$ {with $K_0$}.
To build the linear feedback law, %we need,
we %take 
{draw} inspiration from \cite{GUNZBURGER2000803}, where the following theorem is stated for homogeneous %problems 
{boundary conditions} over the whole boundary $\partial \Omega$:

\begin{theorem}
\label{theo:gunz1}
\cite[Theorem 1]{GUNZBURGER2000803}
Let $f_A = f(u,U) = F - \gamma (u - U)$ be the chosen linear feedback control law with $\gamma > M = \max\{0, -C_0(\nu - K_0\norm{U}_{L^{\infty}((0,T); H^1(\Omega))})\}$. Assuming  $(\nu - K_0\norm{U}_{L^{\infty}((0,T); H^1(\Omega))}) > 0$ 
%\ti{Wouldn't this assumption imply $M=0$? }\ms{yes, but in their proof, they forgot the assumption $(\nu - K_0\norm{U}_{L^{\infty}((0,T); H^1(\Omega))}) > 0$ .}
and $u_{|_{\partial \Omega}} = U_{|_{\partial \Omega}} = 0$,
\begin{itemize}
    \item[(i)] if $(u,p,f_A)$ is a solution to \eqref{eq:general_problem}, then the solution is admissible, i.e., \eqref{eq:convergence} holds,
    \item[(ii)] and 
    $$
    \norm{u(t) - U(t)}^2_{L^2} \leq \norm{u_0 - U_0}^2_{L^2} e^{-2(\gamma - M)t} \quad \text{a.e. } t \in (0,T),
    $$
    i.e., the convergence is exponential.
\end{itemize}
\end{theorem}
\begin{remark}
We stress that, in Theorem \ref{theo:gunz1}, the assumption %\\ 
$(\nu - K_0\norm{U}_{L^{\infty}((0,T); H^1(\Omega))}) > 0$ is essential to reach both %theses 
(i) and (ii). %: 
% The interested reader may refer to the original paper for the proof. 
Indeed, the hypothesis $(\nu - K_0\norm{U}_{L^{\infty}((0,T); H^1(\Omega))}) > 0$ is exploited to make use of a Poincar\'e's inequality necessary to prove the exponential convergence. 
The interested reader may refer to the original paper for the proof. %\\ 

In our {high $Re$} setting, where small values of $\nu$ are investigated, Theorem \ref{theo:gunz1} does not always hold. To overcome this issue, in the novel Theorem~\ref{theo:our2}, we propose an important practical extension of the results in \cite{GUNZBURGER2000803}, %looking for 
{putting forth} a new definition for {the control} $f$ to deal with convection-dominated problems (i.e., higher $Re$) and the presence of mixed boundary conditions. %Indeed, the following holds.
\end{remark}
\begin{theorem}
\label{theo:our2}
Let $f_B = f(u,U) = F + (u - U)\nabla U - \gamma (u - U)$ be %the chosen 
{a} linear feedback control law with $\gamma \geq 0$. Assuming $u_{|_{\Gamma_D}} = U_{|_{\Gamma_D}}$ and $(u \cdot n) \geq 0$ on $\Gamma_N$,
\begin{itemize}
    \item[(i)] if $(u,p,f_B)$ is a solution to \eqref{eq:general_problem}, then the solution is admissible, i.e., \eqref{eq:convergence} holds,
    \item[(ii)] and 
    $$
    \norm{u(t) - U(t)}^2_{L^2} \leq \norm{u_0 - U_0}^2_{L^2} e^{-2(\gamma + C_0 \nu) t} \quad \text{a.e. } t \in (0,T),
    $$
    i.e., %\ 
    the convergence is exponential.
\end{itemize}
\end{theorem}
\begin{remark} {We emphasize that} Theorem \ref{theo:our2} %has no assumptions 
{does not have any restrictions} on the kinematic viscosity $\nu$, i.e., does not depend on the Reynolds number. %, and 
{Furthermore, it does not depend} on the Dirichlet boundary condition $\Gamma_D$ {either}. %On the contrary, 
{In contrast,} Theorem \ref{theo:gunz1} holds 
{only} for %some 
{large} $\nu$ {values (which satisfy the constraint in Theorem~\ref{theo:gunz1})} and for homogeneous Dirichlet conditions all over the boundary, i.e., $\Gamma_D = \partial \Omega$.
%\ti{How does $f_A$ depend on $\Gamma_D$?  We need to explicitly state that in Theorem 1 and/or Remark 1.} \ms{I added the information}
%The proposed extension 
{Thus, the new control in Theorem~\ref{theo:our2} fundamentally}
changes the nature of $f$ {in Theorem~\ref{theo:gunz1}} to deal with larger Reynolds {numbers,} %with respect to the one presented in \cite{GUNZBURGER2000803} 
%and the exponential convergence is recovered, 
%\ti{Wasn't it recovered for $f_A$, too?} \ms{not for large $Re$, due to the assumption over $\nu$}
as we will show in the numerical results of Section \ref{sec:examplecontrolFOM}.\end{remark}
\begin{proof}

First, %of all, 
%let us 
{we} define $w = u - U$. We want to estimate the term $$
 \frac{d}{dt}\norm{u - U}^2_{L^2} =  \frac{d}{dt}\norm{w}^2_{L^2}.
$$
%Let us 
{We} add $\pm \left \langle U_t, v \right \rangle$, $\pm \nu a(U, v)$, $\pm c(U; U, v)$ to the left-hand side of the momentum equation of \eqref{eq:general_problem} {and obtain}
%. We obtain, with $w = u - U$,
\begin{align*}
& \left \langle w_t, v\right \rangle + \nu a(w, v) + \underbrace{c(u; u, v) - c(U; U,v)}_{c(w; w, v) + c(w; U, v) + c(U; w, v)} +  b(v, p)  \\
    & \qquad \qquad \qquad
    + \underbrace{\left \langle U_t, v\right \rangle + \nu a(U, v) + c(U; U, v)}_{\left \langle F , v \right \rangle, \text{ see \eqref{eq:F}}}
     = \left \langle f, v\right \rangle \quad \forall v \in H^1_0(\Omega).
\end{align*}
%Let us 
{We} choose $v = w$ %, we have, 
{and obtain,}
thanks to the divergence-free property of $w$,
\begin{align}
\label{eq:step1}
 \left  \langle w_t, w\right \rangle + \nu a(w, w) + c(w; w, w) + c(w; U, w) + c(U; w, w)
     - \left \langle f - F, w\right \rangle = 0.
\end{align}
%Let us 
{We} note that 
\begin{equation}
    c(w;w,w) = \frac{1}{2} \int_{\Gamma_D \cup \Gamma_N } w^2(w \cdot n)\; dx  \quad \text{and} \quad c(U;w,w) = \frac{1}{2} \int_{\Gamma_D \cup \Gamma_N } w^2(U \cdot n)\; dx.
\end{equation}
%Exploiting 
{Using} $w \in H^1_{\Gamma_D}(\Omega)$ and the %assumption of 
outflow boundary {condition}, 
%\ti{How do we get this? Do we make extra assumptions?} \ms{I added the assumption in the theorem}
we have
%\ti{Is this the ``trick'' you were talking about?} \ms{I discovered it is not a real ``trick" but just an observation over the outflow condition}
\begin{equation}
\label{eq:cwww}
    c(w;w,w) + c(U;w,w) = \frac{1}{2} \int_{\Gamma_N } w^2(u \cdot n)\; dx \geq 0. 
\end{equation}
Moreover, %having chosen 
{choosing} $f = f_B$, we %realize that
{obtain}
\begin{equation}
\label{eq:fB_relation}
\langle f, w \rangle =  \langle F, w \rangle + c(w; U, w) - \gamma \norm{w}_{L^2}^2.
\end{equation}
Plugging \eqref{eq:fB_relation} into \eqref{eq:step1}, applying Poincar\'e's inequality %over 
{to} $\nu a(w,w)$, using %relations 
{inequality}~\eqref{eq:cwww}, and noticing that 
%$
$  \displaystyle \left  \langle w_t, w\right \rangle = \frac{1}{2}\frac{d}{dt}\norm{w}^2_{L^2},$
%$
we obtain 
%\ti{Should the Poincare constant be squared?}
%\begin{align}
%\label{eq:thesi1}
$ \displaystyle
\frac{1}{2}\frac{d}{dt}\norm{w}^2_{L^2} + \nu C_0 \norm{w}_{L^2}^2
     + \gamma \norm{w}_{L^2}^2 \leq 0,
 %\end{align}
 $
 %that is, 
 i.e., 
 % \begin{align}
 $ \displaystyle
        \frac{1}{2}\frac{d}{dt}\norm{w}^2_{L^2} \leq - (\gamma + \nu C_0 )\norm{w}_{L^2}^2
     < 0.$
% \end{align}
%In this way, 
{Thus,} we proved (i). %\\ 
Thesis (ii) is a consequence of %the 
Gronwall's inequality, %obtaining
{which yields}
 \begin{equation}
    \norm{u(t) - U(t)}^2_{L^2} \leq \norm{u_0 - U_0}^2_{L^2} e^{-2 (\gamma + \nu C_0) t} \quad \text{a.e. } t \in (0,T).
\end{equation}
\end{proof}
\begin{remark}
We stress that the assumption of $(u \cdot n) \geq 0$ on $\Gamma_D$ does not allow backflow on the Neumann boundary. This assumption might be restrictive in some settings and real-life scenarios. The hypothesis can be removed if: 
%\begin{itemize}
%\item[$\circ$] 
(i) we assume Dirichlet conditions over the whole boundary and $u_{|_{\partial \Omega}} = U_{|_{\partial \Omega}}$, or
%\item[$\circ$] %defining the 
(ii) we define a nonlinear control law $\tilde{f}_B$ that, in weak form, reads as follows: for all $v \in H^1_{\Gamma_D}(\Omega)$, %reads
$$\left < \tilde f_B, v \right >  = \left < f_B, v \right > + \int_{\Gamma_N}(u - U)(u \cdot n)v \; dx.$$
%\end{itemize}
%Furthermore, 
{We also emphasize that} the outflow assumption is not needed at the discrete level (see Theorem \ref{theo:our1}).
\end{remark}
\reviewerA{
\begin{remark}
\label{rem:construction_control}
The control laws $f_A$ and $f_B$ are built with a \emph{constructive} strategy, with the main goal of verifying relation \eqref{eq:convergence}. We also note that we are not solving an optimal control problem since the control law is not related to any functional to minimize. However, we are proving a stabilizing control that lets the solution converge toward a desired configuration exponentially in time.
\end{remark}
}

\subsection{Discrete formulation}
\label{sec:problem-discrete}

{In this section, we extend the novel feedback control strategy for high Reynolds numbers from the continuous setting (see Section~\ref{sec:problem-continuous}) to the discrete setting.  
To this end, we use the same framework as that employed 
}
%Let us now focus on the discrete version of the problem. The setting is the one already described 
for the uncontrolled problem in Section \ref{sec:efrFOMuncontrolled}. %Namely, once %performed 
%we perform a $\mathbb P^2 - \mathbb P^1$ Taylor-Hood approximation in space and {implicit Euler} in time, we want to solve the following system in the weak form: 
{Specifically, we use a $\mathbb P^2 - \mathbb P^1$ Taylor-Hood spatial discretization and the implicit Euler time discretization, which yields the following system:}
\begin{equation}
\begin{cases}
	\displaystyle \int_{\Omega}\frac{u^{n+1} - u^{n}} {\Delta t} \cdot v_h \; dx
	+ \nu a(u^{n+1},  v_h) + \widetilde{c}(u^{n+1}; u^{n+1}, v_h)  & 
	       \\ \qquad  \qquad \qquad \qquad \qquad- b(v_n, p^{n+1})
	= \displaystyle \int_{\Omega} f^{n+1}(u^{n+1}, U^n, U^{n+1}) \cdot v_h \; dx & \forall v_h \in \mathbb U^{N_h^u}, \\
b( u^{n+1},  q_h ) = 0, & \forall q_h \in \mathbb Q^{N_h^p},
\end{cases}
\label{eq:controlFE}
\end{equation}
where $ \widetilde{c}(\cdot; \cdot, \cdot)$ is the skew-symmetric approximation of the form ${c}(\cdot; \cdot, \cdot)$ (%we refer the reader to 
{see}  \cite{layton2008introduction,temam2001navier}) defined as 
\begin{equation}
    \label{eq:NSskew}
    \widetilde{c}(u; v, z)
    = \frac{1}{2} [c(u; v, z) - c(u;z,v)] \quad \forall \, u,v,z \in U^{N_h^u}.
\end{equation}
We recall that the superscript $^{n}$ denotes a variable evaluated at the time $t^n$. %\\
We note that the discrete control $f^{n+1}$ depends on the velocity variable at time $t^{n+1}$ and on the desired profile at times $t^{n}$ and $t^{n+1}$. Indeed, the new control law we propose is related to {the forcing term} $F$ %as 
introduced in \eqref{eq:F}. In %a 
the fully-discrete setting, we approximate $F$ as % follows: 
%\ti{To be consistent, should we replace the $h$ superscript with the $n$ superscript?  For example, should $F^{n+1}$ be replaced with $F^{n+1}$?}
%\begin{equation}
%\label{eq:Fh}
$F^{n+1} = \displaystyle \frac{U^{n+1} - U^{n}}{\Delta t} + %v 
{\nu} \Delta U^{n+1} + (U^{n+1} \cdot \nabla)U^{n+1}. 
$ %\end{equation}
%\\ 
We also recall that {$U$} is taken divergence-free at each time instance, %namely 
{which implies that} $b(U^{n+1}, q_h)=0$ for all $q_h \in \mathbb Q^{N_h^p}$ and $n = 0, \dots, N_T$. In the fully-discrete framework, the definition of \emph{admissible solution} translates into
\begin{equation}
    \label{eq:FEconvergence}
    \norm{u^{n + 1} - U^{n + 1}}_{L^2}^2 \leq \norm{u^{n} - U^{n}}_{L^2}^2 \quad \text{ for } \quad  n=0,\dots, N_T - 1.
\end{equation}
In this setting, %in \cite{GUNZBURGER2000803} they propose 
the following theorem {is proved in \cite{GUNZBURGER2000803}}:
\begin{theorem}
\label{theo:gunz3}
\cite[Theorem 3]{GUNZBURGER2000803}
Let $f_A^{n+1} = f^{n+1}(u^{n+1}, U^n, U^{n+1}) = F^{n+1} - \gamma (u^{n+1} - U^{n+1})$ 
be the chosen linear feedback control law with $\gamma > M = \max\{0, -C_0(\nu - K_0\norm{U}_{L^{\infty}((0,T); H^1(\Omega))})\}$. Assuming  $(\nu - K_0\norm{U}_{L^{\infty}((0,T); H^1(\Omega))}) > 0$ and $u_{|_{\partial \Omega}} = U_{|_{\partial \Omega}} = 0$,
\begin{itemize}
    \item[(i)] if $(u^{n+1},p^{n+1},f_A^{n+1})$ is a solution to \eqref{eq:controlFE}, then the solution is admissible, i.e.,  \eqref{eq:FEconvergence} holds,
    \item[(ii)] and 
    $$
    \norm{u^{n+1} - U^{n+1}}^2 \leq \left ( \frac{1}{1+2\Delta t (\gamma + C_0(\nu - K_0\norm{U}_{L^{\infty}((0,T); H^1(\Omega)})}\right )^{n+1}\norm{u_0^h - U_0^h}^2.
    $$
\end{itemize}
\end{theorem}
%As already stated at the continuous level, we want to generalize the theorem 
In the %next 
novel Theorem~\ref{theo:our1}, we generalize Theorem~\ref{theo:gunz3} 
to cases where $(\nu - K_0\norm{U}_{L^{\infty}((0,T); H^1(\Omega)})$ can be negative, i.e., {to the case $\nu \ll 1$, which is generally the case of interest in realistic settings}. 
Before stating the theorem, we define %weakly 
{the weak form of} the control $f_B$ at time $t^{n+1}$ as %:
\begin{equation}
  \label{eq:weakFB}  
\left < f_B^{n+1}, v \right > = \left < F^{n+1}_B, v \right > + \tilde c(u^{n+1} - U^{n+1}; U^{n+1}, v) - \gamma \left <u^{n+1} - U^{n+1}, v \right >,
\end{equation}
for all $v \in H^1_{\Gamma_D}(\Omega)$,
where 
\begin{equation}
\label{eq:FBpart1}
\left < F^{n+1}_B, v \right > = \int_{\Omega}\frac{U^{n+1} - U^{n}} {\Delta t} v \; dx + \nu a(U^{n+1}, v) + \widetilde{c}(U^{n+1}; U^{n+1}, v).
\end{equation}

\begin{theorem}
\label{theo:our1}

Let $f_B^{n+1}$ be the %chosen 
{novel} linear feedback control law defined %as 
in \eqref{eq:weakFB} with $\gamma \geq 0$ and $u_{|_{\Gamma_D}} = U_{|_{\Gamma_D}}$. %, then
\begin{itemize}
    \item[(i)] {If} $(u^{n+1},p^{n+1},f_B^{n+1})$ is a solution to \eqref{eq:controlFE}, then the solution is admissible, i.e., \eqref{eq:FEconvergence} holds,
    \item[(ii)] and 
$$
    \norm{u^{n+1} - U^{n+1}}^2_{L^2} \leq \left ( \frac{1}{1+2\Delta t (\gamma + C_0\nu)}\right )^{n+1}\norm{u_0^h - U_0^h}^2_{L^2}.
    $$
\end{itemize}
\end{theorem}
\begin{proof}
For the sake of notation, %let us 
{we} define $w^{n+1} = u^{n+1} - U^{n+1}$. Notice that, thanks to the assumption $u_{|_{\Gamma_D}} = U_{|_{\Gamma_D}}$, $w^{n+1} \in \mathbb U^{N_h^u} \subset H^1_{\Gamma_D}(\Omega)$. We want to bound the term $\norm{w^{n+1}}_{L^2}$ with $\norm{w^{n}}_{L^2}$.
To this end, we add the following terms to the left-hand side of the momentum equation of \eqref{eq:controlFE}: 
\begin{equation}
\label{eq:lhsadd}
\pm \int_{\Omega} \frac{U^{n+1}}{\Delta t}v_h \; dx, \quad \pm \int_{\Omega}\frac{U^{n}}{\Delta t}v_h\; dx, \quad \pm \nu a(U^{n+1}, v_h), \quad \text{and} \quad \pm \widetilde{c}(U^{n+1}; U^{n+1}, v_h). \end{equation}
{Recalling}
%obtaining, recalling 
that $w^{n+1} = u^{n+1} - U^{n+1}$, {we obtain}
\begin{equation}
\label{eq:step1fd}
	\displaystyle \int_{\Omega}\frac{w^{n+1} - w^{n}} {\Delta t} v_h \; dx
	+ \nu a(w^{n+1},  v_h) + \underbrace{\widetilde{c}(u^{n+1}; u^{n+1}, v_h) - \widetilde{c}(U^{n+1}; U^{n+1},v)}_{\tilde{c}(w^{n+1}; w^{n+1}, v_h) + \tilde{c}(w^{n+1}; U^{n+1}, v_h) + \tilde{c}(U^{n+1}; w^{n+1}, v_h)}
	\end{equation}
	$$
	\qquad \qquad \qquad+ b(v_h, p^n)
	+ \underbrace{\int_{\Omega}\frac{U^{n+1} - U^{n}} {\Delta t} v_h \; dx + \nu a(U^{n+1}, v_h) + \widetilde{c}(U^{n+1}; U^{n+1}, v_h)}_{\left <F^{n+1}, v_h \right > \text{ see \eqref{eq:FBpart1}}}
	= \displaystyle \int_{\Omega} f^{n+1} v_h \; dx.
$$

Let us choose $v_h = w^{n+1}$. %Thanks to the definition of 
Since $\widetilde{c}(\cdot;\cdot,\cdot)$ is skew-symmetric, we have $\widetilde{c}(u; v, v) = 0$ for every $u, v$ in $\mathbb U^{N^u_h}$. Moreover, notice that $w^{n+1}$ is %\ti{weakly?} 
divergence-free as {the difference} of two %\ti{weakly?} \ms{$u$ is divergence-free at the strong level and we assumed $U$ divergence-free at the strong level too} 
divergence-free functions. Then, equation \eqref{eq:step1fd} becomes
\begin{align}
\label{eq:step1FE}
 & \displaystyle \int_{\Omega}\frac{w^{n+1} - w^{n}} {\Delta t} w_{n+1} \; dx
	 + \nu a(w^{n+1}, w^{n+1}) \\ \nonumber 
	& \qquad \qquad \qquad \qquad+ \widetilde{c}(w^{n+1}; U^{n+1}, w^{n+1})
     - \int_{\Omega} (f^{n+1} - F^{n+1})w^{n+1} \; dx = 0.
\end{align}
%Having chosen 
{Choosing} $\displaystyle \int_{\Omega} f^{n+1} v_h \; dx 
= \left < f^{n+1}, v_h \right > = 
 \left < f_B^{n+1}, v_h \right > $ {in~\eqref{eq:controlFE}}, we %realize that
{obtain}
\begin{equation}
\label{eq:fB_relationFE}
\int_{\Omega}f^{n+1}w^{n+1}\; dx = \left < F^{n+1}, w^{n+1} \right > + \widetilde{c}(w^{n+1}; U^{n+1}, w^{n+1}) - \gamma \norm{w^{n+1}}_{L^2}^2.
\end{equation}
%where,
%\ti{Why do we have this approximation?} \ms{using skew-symmetric approximation proposed in \cite{temam2001navier,layton2008introduction}, we have}
%\ti{I think we need to be careful here.  The skew-symmetric form is an approximation of the trilinear form, but they're not equal.  How big is the difference between the two?}
%$$
%\widetilde{c}(w^{n+1}; %U^{n+1}, w^{n+1}) \sim 
%\int_{\Omega} (u^{n+1} - U^{n+1})\cdot \nabla U^{n+1} w^{n+1} \; dx = c(w^{n+1}; U^{n+1}, w^{n+1})
%$$
%\ti{
%My feeling is that the easiest solution would be to change the linear feedback control law in Theorem 4 from 

%$$
%f_B^{n+1} 
%=  f^{n+1}(u^{n+1}, U^n, U^{n+1})  
%= F^{n+1} + (u^{n+1} - U^{n+1})\nabla U^{n+1} - \gamma (u^{n+1} - U^{n+1}),
%$$ 
%which corresponds to $c$, %to

%$$
%f_C^{n+1} 
%= f^{n+1}(u^{n+1}, U^n, U^{n+1}) 
%= F^{n+1} 
%+ \frac{1}{2} (u^{n+1} - U^{n+1})\nabla U^{n+1} 
%- \frac{1}{2} U^{n+1} \nabla (u^{n+1} - U^{n+1}) 
%- \gamma (u^{n+1} - %U^{n+1}),
%$$
%which corresponds to $\widetilde{c}$.
%\ms{you are totally right. I'll change it}
%}
Plugging \eqref{eq:fB_relationFE} into \eqref{eq:step1FE}, applying Poincar\'e's inequality %over 
{to} $\nu a(w^{n+1},w^{n+1})$,  and 
noticing that 
%\ti{Please check the next inequality since I changed it slightly.}
$$
  \displaystyle \int_{\Omega}\frac{w^{n+1} - w^{n}} {\Delta t} \cdot w^{n+1} \; dx 
  = \frac{1}{2\Delta t} \norm{w^{n+1}}_{L^2}^2 
  + \frac{1}{2\Delta t} \norm{w^{n+1}-w^{n}}_{L^2}^2 
  %- \int_{\Omega}\frac{w^{n} \cdot w^{n+1}} {\Delta t} \; dx ,
 - \frac{1}{2\Delta t} \norm{w^{n}}_{L^2}^2,
$$
%and completing the square, 
we have
\begin{align*}
\frac{1}{2\Delta t}\norm{w^{n+1}}_{L^2}^2 +\frac{1}{2\Delta t}\norm{w^{n+1} - w^{n}}_{L^2}^2   + \nu C_0 \norm{w^{n+1}}_{L^2}^2
     + \gamma \norm{w^{n+1}}_{L^2}^2 & \leq \frac{1}{2\Delta t}\norm{w^{n}}_{L^2}^2.
\end{align*}
Since $\norm{w^{n+1} - w^{n}}_{L^2}^2 %> 
{\geq} 0$, we obtain
\begin{align}
\label{eq:thesi1}
     \norm{w^{n+1}}_{L^2}^2 & \leq 
     \left (\frac{1}{1 + 2\Delta t(\gamma + C_0 \nu)}\right )\norm{w^{n}}_{L^2}^2 ,
\end{align}
%Namely, we proved (i). 
{which proves (i).}
Thesis (ii) is obtained applying {inequality~\eqref{eq:thesi1} recursively.} %the proof in a recursive way,
%\ti{What do mean by recursive way?} \ms{repeating the estimate to the right-hand side on $\norm{w^{n}}_{L^2}^2$ for $n=0, \dots, N_T - 1$.}
\end{proof}

\begin{remark}[The role of $\gamma$]
\label{remark:gamma}
In the controlled setting, the parameter $\gamma$ represents an \reviewerA{a priori chosen} penalization parameter with respect to the control action. %Indeed, a 
{A} large value of $\gamma$ allows a faster convergence towards the desired state. {However, %other cases are worth mentioning, such as small values of $\gamma$. Indeed, 
in real applications, a large value of $\gamma$ translates into a large %r  
physical and economic effort in controlling the system. %Usually, 
{Since} one wants to spend %less 
as few resources as possible to reach the goal,} %}. Thus, 
we %want to 
investigate the cases where $\gamma \rightarrow 0$ and how %it 
{this} affects the convergence rate. Moreover, we stress that, since we are working with large $Re$, $\gamma \rightarrow 0$ translates into a system that may feature numerical instabilities, %that might need a 
{which may require} further stabilization besides the control action.
This feature has already been observed in the optimal control framework for convection-dominated advection-diffusion equations, see, e.g., \cite{zoccolan2,zoccolan1}.
\end{remark}

\section{%EFR Control for FOM and ROM
{EFR stabilization for the new feedback control}}
    \label{sec:efr-control}

{
As explained in Remark~\ref{remark:gamma}, when $\gamma \rightarrow 0$, further stabilization may be needed in addition to the control stabilizing effect.
In this section, we leverage the EFR strategy to stabilize the novel feedback control introduced in Section~\ref{sec:problem} in the convection-dominated setting (i.e., for large Reynolds numbers).
Specifically, in Section~\ref{sec:efr-control-efr}, we outline and analyze the new feedback control strategy with EFR stabilization at the FOM level.
In Section~\ref{sec:aEFRsec}, we propose a new adaptive EFR strategy to improve the accuracy of the feedback control at the FOM level.
Finally, we compare the feedback control without (Section~\ref{sec:POD_t_c}) and with (Section~\ref{sec:differential-filter}) EFR stabilization at the ROM level.
}
%This Section focuses on the EFR strategy in the controlled NS setting, both at the FOM and ROM level.
\subsection{EFR algorithm}
    \label{sec:efr-control-efr}
%For the control %problem 
%{case,} the EFR strategy is analogous to the uncontrolled case. 
{The EFR strategy for the control case is similar to the EFR strategy for the uncontrolled case.}
The main difference is %in 
the presence of $f^{n+1}(\widetilde{u}^{n+1}, U^n, U^{n+1})$ on the right-hand side of the evolve step. Indeed,
%relying on 
{using the} implicit Euler {method} for the time %evolution, 
discretization with {the relaxation parameter} $\chi \in [0,1]$, %relaxation parameter, 
the EFR approach for the controlled system at the time $t^{n+1}$ reads:
%\ti{I think we should use $\widetilde{u}$ for step I and $\overline{u}$ for filtering in step II, which is the traditional notation. Overline is reserved for filtering :)}
\begin{eqnarray}
         &	\text{(I)}_c& \text{\emph{ Evolve}:} \quad 
\begin{cases}
        	 \displaystyle \frac{{\widetilde{u}}^{n + 1} - u^n}{\Delta t} + (\widetilde{u}^{n+1} \cdot \nabla)    \widetilde{u}^{n+1} & \\  \qquad \qquad \qquad \qquad \quad - \nu \Delta \widetilde{u}^{n+1} + \nabla p^{n+1} = f^{n+1}(\widetilde{u}^{n+1}, U^n, U^{n+1}) & \text{in } \Omega , \vspace{1mm}\\
\nabla \cdot \widetilde{u}^{n+1} = 0 & \text{in } \Omega , \vspace{1mm}\\
\widetilde{u}^{n+1} = u_D^{n+1} & \text{on } \Gamma_D , \vspace{1mm}\\
\displaystyle -p^{n+1} n + \frac{\partial {\widetilde{u}^{n+1}}}{\partial n} = 0  & \text{on } \Gamma_N . \\
\end{cases}
            \label{eqn:ef-rom-1}\nonumber \\[0.3cm]
            &	\text{(II)}_c &\text{\emph{ Filter:}} \quad
\begin{cases} 
        	 -\delta^2 \, \Delta \overline{u}^{n+1} +  \overline{u}^{n+1} = \widetilde{u}^{n+1}& \text{in } \Omega , \vspace{1mm}\\
\overline{u}^{n+1} = u^{n+1}_D &\text{on } \Gamma_D , \vspace{1mm}\\
\displaystyle \frac{\partial \overline{u}^{n+1}}{\partial  n} = 0  & \text{on } \Gamma_N .
\end{cases}
	\label{eqn:ef-rom-2} \nonumber \\[0.3cm]
            &	\text{(III)}_c &\text{ \emph{ Relax:}} \qquad 
        	   u^{n+1}
            = (1 - \chi) \, \widetilde{u}^{n+1}
            + \chi \, \overline{u}^{n+1} \, .
            \label{eq:controlEFR}\nonumber
\end{eqnarray}
The only difference %with 
{from} the uncontrolled setting is in step (I), since only there the control action is present. %We stress that also for the EFR strategy, a convergence result can be proved.
Next, we prove a convergence result for the new feedback control strategy with EFR stabilization, which is outlined in steps {(I)$_{c}$--(III)$_{c}$}.

\begin{theorem}
\label{theo:our5}

%\ti{Should the control depend on the filtered velocity?} \ms{we can investigate that, paper $10^{th}$? :) I'm kidding, with this formulation we recover an exponential decay also for the EFR-FOM approach}
Let \eqref{eq:weakFB} %$f_B^{n+1} =  f^{n+1}(\widetilde{u}^{n+1}, U^n, U^{n+1})  = F^{n+1} + (\widetilde{u}^{n+1} - U^{n+1})\nabla U^{n+1} - \gamma (\widetilde{u}^{n+1} - U^{n+1})$ 
be the chosen linear feedback control law with $\gamma \geq 0$, {$0 < \chi < 1$}, $\overline{C}(\widetilde{u}^k) = C( \delta h^2_{max} +  h^3_{max} +  \delta^2 \norm{\Delta \widetilde{u}^{k}}_{L^2})$ for some positive constant $C$ and $k =0, \dots, n+1$. %Assuming 
{We assume that} $\widetilde{u}_{|_{\Gamma_D}} = U_{|_{\Gamma_D}}$. %, then
% \begin{itemize}
%     \item[(i)] if
{If} $(\widetilde{u}^{n+1},p^{n+1},f_B^{n+1})$ is a solution to step (I)$_c$ and $u^{n+1}$ is the %evolved-filter-relax 
    {EFR} velocity defined in (III)$_c$,  then, for any $0 < \varepsilon < 1$, the following holds:  
\begin{align}
    \norm{u^{n+1} - U^{n+1}}^2_{L^2} & \leq \left ( \frac{1}{(1 - \varepsilon)(1+2\Delta t (\gamma + C_0\nu))}\right )^{n+1}\norm{u_0^h - U_0^h}^2_{L^2} \nonumber\\
    & \label{eq:bound} \qquad \qquad \qquad + \frac{\chi^2}{\varepsilon (1 - \varepsilon)(1+2\Delta t (\gamma + C_0\nu))}\large{ \sum_{i = 1}^{n}} \overline{C}(\widetilde{u}^{i})^2 +  \frac{\chi^2\overline{C}(\widetilde{u}^{n+1})^2}{\varepsilon}.
    \end{align}
%\end{itemize}

\begin{comment}
\begin{itemize}
    \item[(i)] if $(\widetilde{u}^{n+1},p^{n+1},f_B^{n+1})$ is a solution to step (I)$_c$ and $u^{n+1}$ is the %evolved-filter-relax 
    {EFR} velocity defined in (III)$_c$, then the following holds:  
\begin{align}
    \norm{u^{n+1} - U^{n+1}}^2_{L^2} & \leq \left ( \frac{1}{(1 - \chi)(1+2\Delta t (\gamma + C_0\nu))}\right )^{n+1}\norm{u_0^h - U_0^h}^2_{L^2} \nonumber\\
    & \label{eq:bound} \qquad \qquad \qquad + \sum_{i = 0}^{n}\left ( \frac{4\chi }{(1 - \chi)(1+2\Delta t (\gamma + C_0\nu))}\right )^{i + 1} \norm{U^{n-i}}^2_{L^2}\\ 
    &  \qquad \qquad \qquad \qquad \qquad \qquad \qquad \qquad \qquad \qquad \qquad \qquad \qquad { + 4\chi \norm{U^{n+1}}^2_{L^2}} \nonumber .
    \end{align}
   
\end{itemize}
\end{comment}

%fb{why isn't $\frac{\chi - \chi^2}{1-\chi}$ simplified to $\chi$?}
\end{theorem}

\begin{proof}
% \ms{I think that the magenta term confirms that our problem can feature a plateau behaviour, as depicted in Figure \ref{fig:convergencegamma0001}}
% \ti{I think that even the sum looks problematic.  It seems that we cannot prove that the EFR solution is admissible, i.e., it satisfies (20).  Is that correct?} \ms{you are right, I tried to explain ourselves in Remark 5 (just above)}
We %here 
define the variables $\widetilde{w}^{n+1} = \widetilde{u}^{n+1} - U^{n+1}$ and
$w^{n}= u^{n} - U^n$. Both %the 
variables are in $H^1_{\Gamma_D}(\Omega)$. The proof is analogous to the %one 
{proof} of Theorem \ref{theo:our1}. Indeed, we add the quantities in \eqref{eq:lhsadd} to the left-hand side of the evolve step in (I)$_c$. %With the same computations, choosing 
{Choosing} $v_h = \widetilde{w}^{n+1}$ and %by 
{using the} definition of $f_{B}$, we obtain the relation
\begin{align}
\label{eq:thesi1_EFR}
     \norm{\widetilde{w}^{n+1}}_{L^2}^2 & \leq 
     \left (\frac{1}{1 + 2\Delta t(\gamma + C_0 \nu)}\right )\norm{w^{n}}_{L^2}^2.
\end{align}
It remains to prove that $\norm{\widetilde{w}^{n+1}}_{L^2}^2 \geq C_1 
\norm{{w}^{n+1}}_{L^2}^2 - C_2$, for some $C_1>0$ and $C_2 \geq 0$. 
We start with the relax step of the EFR algorithm
$    u^{n+1}
    = (1 - \chi) \, \widetilde{u}^{n+1}
    + \chi \, \overline{{u}}^{n+1}
    \label{eqn-remark:traian-1}
$
and subtract $U^{n+1}$ from both sides:
\begin{align*}
    w^{n+1} = u^{n+1} - U^{n+1}
    & = (1 - \chi) \widetilde{u}^{n+1}
    + \chi \overline{{u}}^{n+1}
    - (1 - \chi) U^{n+1}
    - \chi U^{n+1} 
    \nonumber \\
   & = (1 - \chi) \widetilde{w}^{n+1}
   + \chi (\overline{{u}}^{n+1} - U^{n+1})
    \nonumber \\
   & = (1 - \chi) \widetilde{w}^{n+1}
   + \chi (\overline{{u}}^{n+1} - \widetilde{u}^{n+1})
   + \chi (\widetilde{u}^{n+1} - U^{n+1})
    \nonumber \\
   & = (1 - \chi) \widetilde{w}^{n+1}
   + \chi (\overline{{u}}^{n+1} - \widetilde{u}^{n+1})
   + \chi \widetilde{w}^{n+1}
    \nonumber  = \widetilde{w}^{n+1}
   + \chi (\overline{{u}}^{n+1} - \widetilde{u}^{n+1}).
    \nonumber 
\end{align*}
By the triangle inequality and by Lemma 2.12 in \cite{layton2008numerical}, %\ms{Traian, I'll cite your paper \cite{xie2018numerical} in the ROM part, see Remark \ref{rem:T}} \ti{There is no need to cite our paper!}
we obtain, for some positive constant $C$,
\begin{align}
\label{eq:C_delta_h_Delta}
    \norm{w^{n+1}}_{L^2} 
   & \leq \norm{\widetilde{w}^{n+1}}_{L^2}
   + \chi \norm{\overline{u}^{n+1} - \widetilde{u}^{n+1}}_{L^2} \nonumber \\
   & \leq \norm{\widetilde{w}^{n+1}}_{L^2} + \chi C( \delta h^2_{max} +  h^3_{max} +  \delta^2 \norm{\Delta \widetilde{u}^{n+1}}_{L^2}).\\
\end{align}
   For the sake of clarity, let us define
   $\overline{C}(\widetilde{u}^{n+1}) = C( \delta h^2_{max} +  h^3_{max} +  \delta^2 \norm{\Delta \widetilde{u}^{n+1}}_{L^2})$.
%Relation 
{Inequality} \eqref{eq:C_delta_h_Delta} implies, exploiting Young's inequality with $0 < \varepsilon < 1$, 
\begin{align*}
\norm{\widetilde{w}^{n+1}}^2_{L^2} \geq (\norm{w^{n+1}}_{L^2} - \chi \overline{C}(\widetilde{u}^{n+1}))^2
    & = \norm{w^{n+1}}^2_{L^2} + \chi^2\overline{C}(\widetilde{u}^{n+1})^2 - 2\chi\overline{C}(\widetilde{u}^{n+1}) \norm{w^{n+1}}_{L^2} \\ 
    & \geq \norm{w^{n+1}}^2_{L^2} + \chi^2\overline{C}(\widetilde{u}^{n+1})^2 - {\norm{w^{n+1}}^2_{L^2}}{\varepsilon} - \frac{\chi^2 \overline{C}(\widetilde{u}^{n+1})^2}{\varepsilon}\\
    & = (1 - \varepsilon) \left ( \norm{w^{n+1}}^2_{L^2} - \frac{\chi^2 \overline{C}(\widetilde{u}^{n+1})^2}{\varepsilon} \right ). %\\
\end{align*}
Exploiting this relation in \eqref{eq:thesi1_EFR}, we obtain
\begin{equation}
\label{eq:ricursiveERF}
    \norm{{w}^{n+1}}_{L^2}^2 \leq 
     \left (\frac{1}{(1-\varepsilon)(1 + 2\Delta t(\gamma + C_0 \nu))}\right )\norm{w^{n}}_{L^2}^2 + \frac{\chi^2 \overline{C}(\widetilde{u}^{n+1})^2}{\varepsilon}. 
\end{equation}
Applying \eqref{eq:ricursiveERF} %also 
to $w^n$, we have  
\begin{align*}
\label{eq:ricursiveERF2}
    \norm{{w}^{n+1}}_{L^2}^2 \leq 
     & \left (\frac{1}{(1-\varepsilon)(1 + 2\Delta t(\gamma + C_0 \nu))}\right )^2\norm{w^{n - 1}}_{L^2}^2 \\
     & \qquad \qquad \qquad \qquad + \frac{\chi^2 \overline{C}(\widetilde{u}^{n})^2}{\varepsilon(1-\varepsilon)(1 + 2\Delta t(\gamma + C_0 \nu))}   + \frac{\chi^2\overline{C}(\widetilde{u}^{n+1})^2}{\varepsilon}.
\end{align*}
Finally, applying relation \eqref{eq:ricursiveERF} recursively, we obtain the thesis.
%\ms{Dear Traian, I have a question: is this estimate ``a priori"? Everything is computable on the rhs (I mean, I cannot compute some constants, but the Laplacian of the discrete solution is computable). Am I missing something? }
%\ti{That's a good question, Maria :)  If you're referring to \eqref{eq:bound}, I don't think that's an a priori bound since it depends on $\widetilde{u}^{n+1}$, which is not known before running the code.  In order for the bound to be a priori, all the quantities on the RHS must depend on quantities that can be calculated before the code is run, e.g., $u$.}
\end{proof}

\begin{remark}[on the EFR convergence] \label{rem:conv}
From Theorem \ref{theo:our5}, we observe that EFR strategy worsens the bound in Theorem \ref{theo:our1} since no exponential convergence and admissibility are guaranteed. 
Next, we analyze the three terms separately.
%\begin{enumerate}
%\item 
%$\displaystyle \left ( \frac{1}{(1 - \varepsilon)(1+2\Delta t (\gamma + C_0\nu))}\right )^{n+1}\norm{u_0^h - U_0^h}^2_{L^2}$: %it 
{It} is clear that, for %the 
{an} arbitrary small $\varepsilon$, {a} small $\gamma$, and %for 
{a} large Reynolds 
number, %this 
the term $\displaystyle \left ( \frac{1}{(1 - \varepsilon)(1+2\Delta t (\gamma + C_0\nu))}\right )^{n+1}\norm{u_0^h - U_0^h}^2_{L^2}$ scales as %follows:} %\vspace{2mm}
%$$\left ( \frac{1}{(1 - \varepsilon)(1+2\Delta t (\gamma + C_0\nu))}\right )^{n+1}\norm{u_0^h - U_0^h}^2_{L^2} \sim \norm{u_0^h - U_0^h}^2_{L^2}.$$
$\displaystyle \norm{u_0^h - U_0^h}^2_{L^2}.$
    %\item
    %$ \displaystyle \frac{\chi^2}{\varepsilon (1 - \varepsilon)(1+2\Delta t (\gamma + C_0\nu))} \sum_{i = 1}^{n}\overline{C}(\widetilde{u}^{i})^2 $: %no hints can be guessed by this term, 
    The scaling for the term $\displaystyle \frac{\chi^2}{\varepsilon (1 - \varepsilon)(1+2\Delta t (\gamma + C_0\nu))} \sum_{i = 1}^{n}\overline{C}(\widetilde{u}^{i})^2$ is not clear,
    since, for $k=0,\dots, n$, $\overline{C}(\widetilde{u}^k)$ %, for $k=0,\dots, n$ 
    can be large. Everything depends on %scale observations between 
    the {scalings of the} various constants and, thus, on the problem at hand. %\vspace{2mm}
    %\item 
    %$\displaystyle \frac{\chi^2\overline{C}(\widetilde{u}^{n+1})^2}{\varepsilon}$: %the 
    The considerations %of 
    for the second term %{(2)} 
    also apply for the term $\displaystyle \frac{\chi^2\overline{C}(\widetilde{u}^{n+1})^2}{\varepsilon}$.  
%\end{enumerate}
A possible choice for $\varepsilon$ is $\varepsilon = \chi$. Indeed, in our numerical frameworks $\chi \ll 1$. This is not a restrictive hypothesis and %it 
is a common choice in literature~\cite{bertagna2016deconvolution,Strazzullo20223148}, as we will %address 
{explain} in Section \ref{sec:aEFRsec}. %In this 
This way, the considerations %on 
for the first term %(1) 
are still valid, and the following scalings for the second and third terms %(2) and (3) %it holds:
hold: 
$\displaystyle 
%\text{(2)} \sim 
\frac{\chi}{(1 - \chi)(1+2\Delta t (\gamma + C_0\nu))} \sum_{i = 1}^{n}\overline{C}(\widetilde{u}^{i})^2 \sim \chi \sum_{i = 1}^{n}\overline{C}(\widetilde{u}^{i})^2$ %\quad \text{and \quad (3)} \sim 
and ${\chi\overline{C}(\widetilde{u}^{n+1})^2}$, respectively. 
%$$
This way, the small value of $\chi$ can balance the possibly large values of $\overline{C}(\widetilde{u}^k)$, for $k=0, \dots, n+1$. 
%Anyway, the observations might translate into 
{These observations suggest possible} slow convergence or plateau phenomena even for small values of $\chi$, as we will see in the numerical results presented in Section \ref{sec:examplecontrolFOM}. 
%However, we 
{We} propose a solution to this issue in Section \ref{sec:aEFRsec}. %\ms{Do you think this remark is useful? Maybe we can get rid of it.}
%\ti{I actually like it, although it does seem to open a new can of worms :)  But we'll have to deal with these issues in the numerical Section anyway, so why not use theory to foreshadow that? Plus, this remark motivates the introduction of aEFR in the next Section.}
\end{remark}

\subsection{{A new adaptive EFR (aEFR) algorithm}}
    \label{sec:aEFRsec}

{
The numerical investigation in Section \ref{sec:examplecontrolFOM} shows that employing EFR with small relaxation parameters alleviates the spurious numerical oscillations and allows us to reach the desired state faster than noEFR in the first part of the simulation.
However, as explained in Remark~\ref{rem:conv},  
}
Theorem \ref{theo:our5} does not guarantee an exponential convergence as stated in Theorem \ref{theo:our1} for the noEFR strategy. %However, numerical observations in Section \ref{sec:examplecontrolFOM}, show that employing EFR with small relaxation parameters alleviates spurious oscillations and this allows it to reach the desired state faster than noEFR in the first part of the simulation. 
% We remark that we always work with the assumption $\chi$ to be small. It is a classic choice in literature, see e.g.\ \cite{bertagna2016deconvolution,strazzullo2022consistency}, where they advocate $\chi \sim K \Delta t$, for some $K > 0$. Nevertheless, the thesis of Theorem \ref{theo:our5} shows some %criticalities 
% {potential issues} in the convergence rate, see remark \ref{rem:conv}. %Indeed, the term $4\chi\norm{U^{n+1}}^2_{L^2}$ may be large even for $n \rightarrow \infty$ and small $\chi$. % \fb{I guess this sentence is saying ``looking at the theorem, we can see what could go wrong'', but I would rephrase it}, while this is not true for the noEFR process. 

This is the reason why, {in this section,} we propose an adaptive-EFR (aEFR) strategy, {which is outlined in Algorithm \ref{alg:aEFR}}. %Namely, 
{In the new aEFR strategy,} given a tolerance $\tau$, we apply EFR if $\norm{w^n}_{L^2} = \norm{u^{n} - U^n}_{L^2}^2 \geq \tau $. %otherwise, 
{Otherwise, we apply the} standard controlled NSE simulation \eqref{eq:controlFE}. %is performed 
%\ti{Please check equation numbers above.}
%\fb{can I play devil's advocate here? Phrased in this way I am afraid that a nasty reviewer will come up with the following comment: $U$ is basically the NSE solution you expect at steady state, so if you can turn off EFR when $u^n$ is very close to $U$ and your simulation still ``works'', then it means that your Re was not so large (say, to necessitate EFR) in the first place. The way I would answer my nasty self is that there may well be situations (like ours) in which EFR is needed in the transient but maybe not at steady state, as hinted in the next sentence. Traian, do you have a suggestion on how to phrase these first few sentences of aEFR in a way that people cannot come up with this comment?}. 
%\ti{Good question, Francesco.  How about we illustrate this numerically in Section 4, and here we just comment on that?}
%\ti{I tried to do that below.} \ms{It sounds convincing to me :) }
%This technique will 
The main goal of the new aEFR strategy is to alleviate %allow alleviating 
the oscillations in the first part of the time evolution %while, 
and, when a good %convergence to 
approximation of $U^n$ is reached, %it recovers 
to recover the exponential bound stated in Theorem \ref{theo:our1}. 
The aEFR strategy aims at tackling those settings in which, e.g., the EFR is needed in the transient regime but not for the steady state.
An example of this type of setting is that used in the numerical investigation in Section \ref{sec:examplecontrolFOM}.
%\fb{this somehow is saying that if there exists $\hat{n} \in \mathbb{N}$ such that $\norm{w^{\hat{n}}} < \tau$, then we expect that $\norm{w^n} < \tau$ for all $n \geq \hat{n}$? Or am I reading this the wrong way?}\ms{What you said is true. Moreover, $\norm{w^n} \leq e^{\text{some coefficient}}\norm{w^{\hat{n}}}$}. 
%\ti{In that case, should we explicitly state this somewhere?} \ms{I put the reference to the theorem.}
{We note that, although the} %The 
choice of {the parameter $\tau$ in the new aEFR strategy} is problem dependent, %and arbitrary. However, 
{it} can be guided by a threshold value %of accuracy 
{for the difference} between the controlled solution and the desired state. 
Indeed, the numerical results of Section \ref{sec:examplecontrolFOM} will show that
employing {the new aEFR algorithm} allows:
%\begin{itemize}
%\item[$\circ$] 
(i) to alleviate the numerical oscillations while there is room for improvement in reaching the desired state, and
%\item[$\circ$] 
(ii) to recover the exponential convergence expected by Theorem \ref{theo:our1}.
%\end{itemize}

\begin{algorithm}
\caption{aEFR}\label{alg:aEFR}
\begin{algorithmic}[1]
\State{$u_0, u_{in}, \tau$}\Comment{Inputs needed}
\For{$n \in \{1, \dots, N_T\}$}\Comment{Time loop}
\If{$\norm{u^{n} - U^n}_{L^2}^2 \geq \tau$}
\State{(I)$_c$ + (II)$_c$ + (III)$_c$} \Comment{EFR simulation}
\Else
\State{Solve \eqref{eq:controlFE} }
\Comment{Standard controlled %NS simulation
{NSE simulation}}
\EndIf
\EndFor
\end{algorithmic}
\end{algorithm}

\subsection{The %POD-Galerkin 
{EFR-noEFR} (and aEFR-noEFR) approach}
\label{sec:POD_t_c}
{In this section, we use the new feedback control with EFR stabilization at the FOM level, but not at the ROM level.}
To build the reduced bases, we apply the same POD-based ROM strategy %proposed 
{outlined} in Section \ref{sec:POD_t}. Namely, we %propose 
{use} a standard POD for the pressure variable, and a POD with supremizer enrichment technique for the velocity. The main difference is that {in the EFR-noEFR approach,} the snapshots are provided by the solution of the %aforementioned 
regularized control problem {described in Section~\ref{sec:efr-control-efr}}.
Also in this case, we call the reduced spaces as $ \mathbb U^{r_{us}}$ and $\mathbb Q^{r_p}$, {which are} spanned by $\{\phi_i \}_{i=1}^{r_{us}}$ and $\{\psi_i \}_{i=1}^{r_{p}}$, respectively. After the building phase, for each time instance $t^{n+1}$, a standard Galerkin projection is performed in the controlled NSE framework, i.e., we solve
\begin{equation}
\begin{cases}
	\displaystyle \int_{\Omega}\frac{u_{r_{us}}^{n+1} - u_{r_{us}}^{n}} {\Delta t} \phi_i \; dx
	+ \nu a(u_{r_{us}}^{n+1},  \phi_i)
	+ c(u_{r_{us}}^{n+1}; u_{r_{us}}^{n+1} , \phi_i) \\ \qquad \qquad \qquad \qquad \qquad \qquad \qquad - b(\phi_i, p_{r_{p}}^{n+1})
	= \displaystyle \int_{\Omega} f^{n+1}(u^{n+1}_{r_{us}}, U^n, U^{n+1})\phi_i \; dx, \\
b( u_{r_{us}}^{n+1},  \psi_j ) =0,
\end{cases}
\label{eq:controlEFRnoEFR}
\end{equation}
for all $ i = 1, \ldots, r_{us},$ and $j = 1, \ldots, r_p$. Here, $u_{r_{us}}$ and $p_{r_p}$ are the reduced variables as defined in Section \ref{sec:POD_t}. The strategy is summarized in Algorithm \ref{alg:EFRnoEFR}. We stress that %the 
an adaptive version, aEFR-noEFR, can be easily devised, %see 
{as illustrated in} {Algorithm} \ref{alg:aEFRnoEFR}.
\begin{algorithm}[H]
\caption{EFR-noEFR}\label{alg:EFRnoEFR}
\begin{algorithmic}[1]
\State{$u_0, u_{in}, N_{u}, N_{p}$}\Comment{Inputs needed}
\For{$n \in \{1, \dots, N_T\}$}\Comment{Time loop}
\State{(I)$_{c}$ + (II)$_{c}$ + (III)$_{c}$} \Comment{EFR simulation}
\EndFor
\State{$\{u_i\}_{i=1}^{N_{u}} \subseteq \{u^k\}_{k=1}^{N_{T}}$ 
\quad $\{p_i\}_{i=1}^{N_{p}} \subseteq \{p^k\}_{k=1}^{N_{T}}$} \Comment{Snapshots}
\State{$\mathbb U^{r_{us}} \doteq \text{POD}(\{u_i\}_{i=1}^{N_{u}})\oplus \text{POD}(\{S(p_i)_{i=1}^{N_{u}})\}$ } \Comment{Supremizer enrichment for velocity space}
\State{$\mathbb Q^{r_p} \doteq \text{POD}(\{p_i\}_{i=1}^{N_p})$ } \Comment{Standard POD for pressure}
\For{$n \in \{1, \dots, N_T\}$}\Comment{Time loop}
\State{Solve \eqref{eq:controlEFRnoEFR} } \Comment{%Perform a 
Standard controlled %NS equation
{NSE simulation}}
\EndFor
\end{algorithmic}
\end{algorithm}

\begin{algorithm}[H]
\caption{aEFR-noEFR}\label{alg:aEFRnoEFR}
\begin{algorithmic}[1]
\State{$u_0, u_{in}, N_{u}, N_{p}, \tau$}\Comment{Inputs needed}
\State{Apply aEFR (Algorithm \ref{alg:aEFR})}
\State{$\{u_i\}_{i=1}^{N_{u}} \subseteq \{u^k\}_{k=1}^{N_{T}}$ 
\quad $\{p_i\}_{i=1}^{N_{p}} \subseteq \{p^k\}_{k=1}^{N_{T}}$} \Comment{Snapshots}
\State{$\mathbb U^{r_{us}} \doteq \text{POD}(\{u_i\}_{i=1}^{N_{u}})\oplus \text{POD}(\{S(p_i)_{i=1}^{N_{u}})\}$ } \Comment{Supremizer enrichment for velocity space}
\State{$\mathbb Q^{r_p} \doteq \text{POD}(\{p_i\}_{i=1}^{N_p})$ } \Comment{Standard POD for pressure}
\For{$n \in \{1, \dots, N_T\}$}\Comment{Time loop}
\State{Solve \eqref{eq:controlEFRnoEFR} } \Comment{%Perform a 
Standard controlled %NS equation
{NSE simulation}}
\EndFor
\end{algorithmic}
\end{algorithm}

\subsection{EFR-EFR (and aEFR-aEFR) approach}
	\label{sec:differential-filter}
%Also in the controlled case, the EFR-EFR strategy employs the EFR algorithm 
% {at both the FOM and the ROM levels.}
{In this section, we use the new feedback control with EFR stabilization at both the FOM and the ROM level.}
The reduced spaces are built as in Section \ref{sec:POD_t_c}. However, the Galerkin projection is performed for the three steps of the EFR as follows:
%\ti{Should the control depend on the ROM solution?}
\begin{eqnarray*}
         &	\text{(I)}_{cr}& \text{\emph{}} \quad 
\begin{cases}
	\displaystyle \int_{\Omega}\frac{\widetilde{u}_{r_{us}}^{n+1} - u_{r_{us}}^{n}} {\Delta t} \phi_i \; dx
	+ \nu a({\widetilde{u}}_{r_{us}}^{n+1},  \phi_i)
	+ c({\widetilde{u}}_{r_{us}}^{n+1}; {\widetilde{u}}_{r_{us}}^{n+1} , \phi_i)  \\ \qquad \qquad \qquad \qquad \qquad \qquad \qquad - b(\phi_i, p_{r_p}^{n+1})
	= \displaystyle \int_{\Omega} f^{n+1}(\widetilde{u}_{r_{us}}^{n+1}, U^n, U^{n+1})\phi_i \; dx, \\
b( {\widetilde{u}}_{r}^{n+1},  \psi_j ) = 0,
\end{cases}
            \label{eqn:ef-rom-1-r}\nonumber \\[0.3cm]
            &	\text{(II)}_{cr} &\text{\emph{}} \quad
        	 \delta^2 \, a(\overline{u}^{n+1}_{r_{us}}, \phi_i)  +  \int_{\Omega} \overline{u}^{n+1}_{r_{us}}  \phi_i \; dx = \int_{\Omega} \widetilde{u}^{n+1}_{r_{us}}  \phi_i \; dx,\\
        	 %\text{\TI{Should it be a plus sign in front of $2 \delta^2$?}}\\
	\label{eqn:ef-rom-2-r} \nonumber \\[0.3cm]
            &	\text{(III)}_{cr} &\text{\emph{}} \qquad 
        	   u^{n+1}_{r_{us}}
            = (1 - \chi) \widetilde{u}^{n+1}_{r_{us}}
            + \chi \overline{u}^{n+1}_{r_{us}},
            \label{eq:controEFEEFR}\nonumber
\end{eqnarray*}
for all $ i = 1, \ldots, r_{us},$ and $j = 1, \ldots, r_p$. As %did for 
{in the case of} the reduced velocity $u_{r_{us}}$ and %the 
reduced pressure $p_{r_p}$, the reduced variables $\overline{u}_{r_{us}}$ and $\widetilde{u}_{r_{us}}$ are %the ones 
{those} described in Section \ref{sec:efrefr}. The {EFR-EFR} approach is described in {Algorithm} \ref{alg:EFREFR}. Also in this case, we propose the adaptive version, aEFR-aEFR, in Algorithm \ref{alg:aEFRaEFR}. The criterion chosen to apply or not the EFR strategy at the reduced level is
$\norm{\mathsf Q_{\mathbb U^{r_{us}}}^Tu^{n}_r - U^n} \geq \tau$, where $\mathsf Q_{\mathbb U^{r_{us}}}$ is the basis matrix related to the velocity space, %that 
{which} projects back the reduced solution to the FOM space. %\\

For the sake of clarity, {in Table \ref{tab:Acro},} we summarize all the {FOM and ROM} acronyms and the corresponding features. %in Table \ref{tab:Acro}.
\begin{algorithm}[H]
\caption{EFR-EFR}\label{alg:EFREFR}
\begin{algorithmic}[1]
\State{$u_0, u_{in}, N_{u}, N_{p}$}\Comment{Inputs needed}
\For{$n \in \{1, \dots, N_T\}$}\Comment{Time loop}
\State{(I)$_c$ + (II)$_c$ + (III)$_c$} \Comment{EFR simulation}
\EndFor
\State{$\{u_i\}_{i=1}^{N_{u}} \subseteq \{u^k\}_{k=1}^{N_{T}}$ 
\quad $\{p_i\}_{i=1}^{N_{p}} \subseteq \{p^k\}_{k=1}^{N_{T}}$} \Comment{Snapshots}
\State{$\mathbb U^{r_{us}}\doteq \text{POD}(\{u_i\}_{i=1}^{N_{u}})\oplus \text{POD}(\{S(p_i)_{i=1}^{N_{u}})\}$ } \Comment{Supremizer enrichment for velocity}
\State{$\mathbb Q^{r_p} \doteq \text{POD}(\{p_i\}_{i=1}^{N_p})$ } \Comment{Standard POD for pressure}
\For{$n \in \{1, \dots, N_T\}$}\Comment{Time loop}
% \If{$\norm{\mathsf Q_{\mathbb U^{r_us}}^Tu^{n}_r - U^n} \geq \tau$}
\State{(I)$_{cr}$ + (II)$_{cr}$ + (III)$_{cr}$} \Comment{EFR simulation}
%\Else
% \State{Solve \eqref{eq:controlFE} }
%\EndIf
\EndFor
\end{algorithmic}
\end{algorithm}

\begin{algorithm}[H]
\caption{aEFR-aEFR}\label{alg:aEFRaEFR}
\begin{algorithmic}[1]
\State{$u_0, u_{in}, N_{u}, N_{p}, \tau$}\Comment{Inputs needed}
\State{Apply aEFR (Algorithm \ref{alg:aEFR})}
\State{$\{u_i\}_{i=1}^{N_{u}} \subseteq \{u^k\}_{k=1}^{N_{T}}$ 
\quad $\{p_i\}_{i=1}^{N_{p}} \subseteq \{p^k\}_{k=1}^{N_{T}}$} \Comment{Snapshots}
\State{$\mathbb U^{r_{us}} \doteq \text{POD}(\{u_i\}_{i=1}^{N_{u}})\oplus \text{POD}(\{S(p_i)_{i=1}^{N_{u}})\}$ } \Comment{Supremizer enrichment for velocity}
\State{$\mathbb Q^{r_p} \doteq \text{POD}(\{p_i\}_{i=1}^{N_p})$ } \Comment{Standard POD for pressure}
\For{$n \in \{1, \dots, N_T\}$}\Comment{Time loop}
\If{$\norm{\mathsf Q_{\mathbb U^{r_{us}}}^Tu^{n}_r - U^n}^2_{L^2} \geq \tau$}
\State{(I)$_{cr}$ + (II)$_{cr}$ + (III)$_{cr}$} \Comment{EFR simulation}
\Else
\State{Solve {\eqref{eq:controlEFRnoEFR}}} \Comment{Standard controlled NSE simulation}
\EndIf
\EndFor
\end{algorithmic}
\end{algorithm}

\begin{table}[H]
\caption{{FOM and ROM} acronyms. Gray cells %on 
in ROM-columns indicate that no reduction in performed in the algorithm.}
\label{tab:Acro}
\resizebox{\textwidth}{!}{
\begin{tabular}{|c|c|c|c|c|c|c|}
\hline
Acronym & FOM regularization& ROM regularization& FOM adaptivity & ROM adaptivity & Eq. or algorithm\\ \hline
noEFR &   &\cellcolor{gray!50} &  & \cellcolor{gray!50} & eq. \eqref{eq:controlFE}\\ \hline
EFR  &  $\checkmark$ & \cellcolor{gray!50} & & \cellcolor{gray!50} & (I)$_c$ + (II)$_c$ + (III)$_c$\\ \hline
aEFR  &  $\checkmark$ & \cellcolor{gray!50} &  $\checkmark$ & \cellcolor{gray!50} & Algorithm \ref{alg:aEFR} \\ \hline
EFR-noEFR &  $\checkmark$ & & & &
Algorithm \ref{alg:EFRnoEFR}\\ \hline

EFR-EFR &  $\checkmark$ & $\checkmark$ & & & Algorithm \ref{alg:EFREFR}\\ \hline

aEFR-noEFR &  $\checkmark$ &  & $\checkmark$ &  & Algorithm \ref{alg:aEFRnoEFR}\\ \hline

aEFR-aEFR &  $\checkmark$ & $\checkmark$ & $\checkmark$ & $\checkmark$ & Algorithm \ref{alg:aEFRaEFR}\\ \hline

\end{tabular}}
\end{table}

\begin{remark}
\label{rem:T}
No {theoretical} convergence study is performed at the ROM level. However, we think that it is possible to exploit the FOM results to prove a ROM version of Theorem \ref{theo:our1} and Theorem \ref{theo:our5}, exploiting the techniques presented in \cite{xie2018numerical}.
\end{remark}

\section{%The need of EFR and EFR-EFR in the controlled setting
{Numerical Results}}
\label{exp1:fafb}
% This section focuses on the \ti{numerical investigation of the} %need \ti{``need" is too strong, I think} \ms{effect of?} for 
% \ti{important role played by} the EFR algorithm %for 
% \ti{in the} controlled Navier-Stokes \ti{at large Reynolds numbers} both at the FOM and at the ROM level (sections \ref{sec:examplecontrolFOM} and \ref{sec:examplecontrolROM}, respectively). Before this investigation, we numerically analyze the role of the $f_A$ and $f_B$ control in the case of convection-dominated Navier-Stokes equations (section \ref{sec:experiment1}). In the following, we perform three \ti{numerical} experiments, summarized in Table \ref{tab:Exp}:
%\ti{I think we should add corresponding section numbers in the above paragraph.}\ms{done}
{
In this section, we perform a numerical investigation of the novel feedback control and the EFR strategy used in the convection-dominated regime.
To this end, in Section~\ref{sec:experiment1}, we compare the new feedback control strategy, $f_B$, with the standard control approach, $f_A$, in the convection-dominated regime of the NSE at the FOM level.
In Section~\ref{sec:examplecontrolFOM}, we investigate the role played by the EFR stabilization in the FOM of the convection-dominated NSE with a small control parameter $\gamma$ (which is common in realistic applications).
Finally, in Section~\ref{sec:examplecontrolROM}, we investigate the EFR stabilization in the ROM of the convection-dominated NSE with a small control parameter $\gamma$.%\\

Specifically, we perform three numerical experiments, which are summarized in Table \ref{tab:Exp}:
}
\begin{itemize}
    \item [$\circ$] \emph{Experiment 1}. A numerical comparison between %noEFR with 
    {the} $f_A$ and $f_B$ control {strategies} %for several values of $\gamma$. 
    {in the convection-dominated regime at the FOM level, with no EFR stabilization.}
    %We 
    {The goal of this experiment is to}  investigate whether $f_B$ yields more accurate results than $f_A$.
    \item [$\circ$] \emph{Experiment 2}. A numerical comparison between noEFR, EFR,  and aEFR at the FOM level with $f_B$ control and small $\gamma$ {values}. In this test, we investigate whether EFR-based strategies are useful in the controlled setting {with realistic control parameter values}.
    \item [$\circ$] \emph{Experiment 3}. A numerical comparison between noEFR and aEFR at the ROM level, %once performed aEFR to obtain the POD basis functions. 
    {with POD basis functions obtained from an EFR stabilized FOM.}
    This numerical test investigates whether EFR-based strategies are beneficial at the {controlled} ROM level.
\end{itemize}

\begin{table}[H]
\caption{%Scheme of the proposed experiments.
{Summary of the three numerical experiments.}}
\label{tab:Exp}
\centering 
\begin{tabular}{|c|c|c|c|c|c|c|c|}
\hline
\cellcolor{gray!40}& FOM& ROM& $f_A$ & $f_B$ & noEFR& EFR & aEFR      \\ \hline
Experiment 1  & $\checkmark$  & &  $\checkmark$ & $\checkmark$ & $\checkmark$ & &   \\ \hline
Experiment 2  & $\checkmark$  & &  & $\checkmark$ & $\checkmark$ & $\checkmark$& $\checkmark$  \\ \hline
Experiment 3  &   & $\checkmark$ &  & $\checkmark$ & $\checkmark$  & & $\checkmark$   \\ \hline
\end{tabular}
\end{table}
Moreover, the interested reader may find further investigations on the role of the controlled snapshots, on the predictive regime and on reduction both in time and in $\gamma$ in Appendices \ref{rem:uncontrolled}, \ref{rem:pred} and \ref{rem:gammapar}, respectively. 
\subsection{%noEFR and 
$f_A$ vs $f_B$ numerical comparison {for unstabilized FOM (noEFR)}(Experiment 1)}
\label{sec:experiment1}
%We propose a first example to understand the role of $\gamma$, the \emph{control parameter}, in the setting of convection-dominated flows and whether $f_B$ yields more accurate results than $f_A$.\\
{In this numerical experiment, we investigate whether the novel feedback control, $f_B$, yields more accurate results than the classical control approach, $f_A$, at the FOM level.
To this end, we consider an unstabilized FOM (noEFR), and realistic, small values for the control parameter, $\gamma$.}

%\ti{Cite references where this test problem was used.}
The spatial domain is $\Omega \doteq \{(0, 2.2) \times (0, 0.41)\} \setminus \{(x, y)\in \mathbb R^2 \text{ such that } (x - 0.2)^2 + (y - 0.2)^2 - 0.05^2 = 0 \}$, represented in Figure \ref{fig:domain}. 
The Dirichlet boundary condition is 
\begin{equation}
\label{eq:inlet}
u_D = 
\begin{cases}
0 & \text{ on } \Gamma_W, \\
u_{\text{in}} = \displaystyle \left  ( \frac{6}{0.41^2}y(0.41 - y), 0 \right) & \text{ on } \Gamma_{\text{in}}, 
\end{cases}
\end{equation}
where $\Gamma_W$ (solid teal boundary in Figure \ref{fig:domain}) is the union of the bottom ($\Gamma_B$) and %the 
top ($\Gamma_T$) walls of the channel together with the walls of the cylinder ($\Gamma_C$). The inlet condition $u_{\text{in}}$ is applied to $\Gamma_{\text{in}} = \{0\} \times [0, 0.41]$, (magenta dashed line in Figure \ref{fig:domain}). The prescribed inlet condition is given by equation \eqref{eq:inlet}.

% \begin{figure}[H]
% \centering
% \hspace{-8mm}\includegraphics[scale=0.27]{img/mesh/mesh.png}
% \caption{{The mesh.}
% }
% \label{fig:mesh}
% \end{figure}
\begin{figure}[H]
\begin{center}
\begin{tikzpicture}[scale=5.0]

\filldraw[color=teal!80, fill=gray!10, very thick](0.2,0.2) circle (0.05);
\filldraw[color=magenta!90, very thick, dashed](0,0) -- (0.,0.41);
\filldraw[color=teal!80, fill=gray!10, very thick](0,0.41) -- (2.2,0.41);
\filldraw[color=teal!80, fill=gray!10, very thick](0,0.) -- (2.2,0.);
\filldraw[color=black!80, fill=gray!10, very thick, dotted](2.2,0.) -- (2.2,0.41);

\node at (-.1,0.25){\color{black}{$\Gamma_{\text{in}}$}};
\node at (2.3,0.25){\color{black}{$\Gamma_{N}$}};
\node at (0,-.05){\color{black}{$(0,0)$}};
%\node at (3,-.3){\color{black}{$(1,0)$}};
%\node at (7.5,-.3){\color{black}{$(1 + \mu_2,0)$}};
%
\node at (0,0.45){\color{black}{$(0,0.41)$}};
\node at (2.35,0.45){\color{black}{$(2.2,0.41)$}};
\node at (2.3,-0.05){\color{black}{$(2.2,0)$}};
\node at (0.35,.2){\color{black}{$ \Gamma_C$}};
\node at (0.99,.45){\color{black}{$ \Gamma_T$}};
\node at (0.99,-.05){\color{black}{$ \Gamma_B$}};

\end{tikzpicture}
\end{center}
\caption{Schematic representation of the domain $\Omega$. $\Gamma_D = \Gamma_{\text{in}} \cup \Gamma_C \cup \Gamma_T \cup \Gamma_B$. No-slip conditions are applied over the teal solid boundary. The inlet condition is applied to $\Gamma_{\text{in}}$, i.e., the magenta dashed line. The %``do-nothing" 
{``free flow"} boundary condition is applied on $\Gamma_N$, i.e., the black dotted line. %features Neumann boundary conditions.
}
\label{fig:domain}
\end{figure}

For example, this is the case when we set $\nu = 10^{-4}$, %and 
$T=8$, %for 
{and} $\Delta t = 4\cdot 10^{-4}$. %For this numerical test, we refer to \cite{Hinze2000273}. However, we are working with a %significantly larger 
%Reynolds number that is { one order of magnitude higher}: $Re = 1000$.
{Our computational setting is similar to that in \cite{Hinze2000273}, although  the Reynolds number used in our numerical investigation (i.e., $Re=1000$) is one order of magnitude higher.}
We refer the reader to \cite{ST96} for a complete description of the benchmark in an uncontrolled setting, i.e., for $f_A=f_B=0$.%\\ 
Our goal is to reach a desired profile $U$ by means of the linear feedback control actions described in Theorem \ref{theo:gunz1} and Theorem \ref{theo:our1}, namely $f_A$ and $f_B$, respectively.
We recall that, defining $F$ as in \eqref{eq:F},
$f_A = F - \gamma(u - U)$ and $f_B = F + (u - U)\nabla U - \gamma(u - U).$
For the sake of clarity, we say that the system (or the velocity) is A-controlled and B-controlled when $f_A$ and $f_B$ are employed, respectively. 
In this test case, the desired velocity profile solves the following steady Stokes problem in $\Omega$:
%\ti{What's the motivation for choosing the Stokes solution?} \ms{I added the motivation just below the Stokes euqations}
\begin{equation}
\label{eq:Stokes}
\begin{cases}
 - \Delta U + \nabla P = 0 & \text{in }  \Omega, \\
\nabla \cdot U = 0  & \text{in }  \Omega, \\
U = u_D & \text{on } \Gamma_D, \\
\displaystyle -P n + \frac{\partial U}{\partial n} = 0  & \text{on }  \Gamma_N.
\end{cases}
\end{equation}
\begin{figure}
\centering
\hspace{-8mm}\includegraphics[scale=0.175]{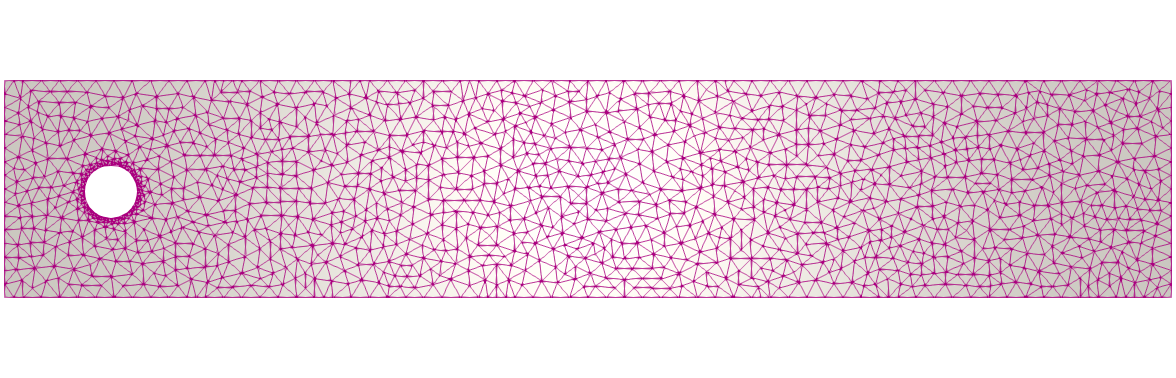}
\caption{{The mesh.}
}
\label{fig:mesh}
\end{figure}

%The solution of \eqref{eq:Stokes} is depicted in Figure \ref{fig:U1}. \ms{The Stokes solution represents a laminar and beneficial desired state for the system we are investigating.} \\
{The Stokes solution, which is displayed in Figure \ref{fig:U1}, is the laminar state that we want the controlled system to achieve.}
First, %of all, 
we stress that the problem we are dealing with needs stabilization if no control is applied,%\footnote{The reader interested in the need for EFR in the uncontrolled context may refer to Appendix \ref{app:EFRapp}.}, 
since numerical oscillations arise even for small $t$. This is shown in Figure \ref{fig:nocontrolvscontrol} (top), where we %show 
{display} the results at $t=0.4$ for a simulation without stabilization. The bottom panels in %the same figure 
{Figure \ref{fig:nocontrolvscontrol}} show %some 
{two} cases in which control is applied. %, instead. 
Moreover, $\Gamma_N$ (black dotted boundary in Figure \ref{fig:domain}) features %``do-nothing'' 
{``free flow"} boundary conditions. 
The initial condition is $u_0 = (0,0)$. The mesh parameters of the triangular mesh (depicted in Figure \ref{fig:mesh}) are %:
 $h_{min} = 4.46 \cdot 10^{-3}$ and $h_{max} = 4.02 \cdot 10^{-2}$. After a Taylor-Hood $\mathbb P^2 - \mathbb P^1$ FE discretization for velocity and pressure, respectively, we obtain a FE space with $N_h \doteq N_h^{u} + N_h^p = 14053$ degrees of freedom. With %this 
 {these} mesh parameters, as stated in \cite{Strazzullo20223148}, we are working in a marginally-resolved regime and, for high $Re$ values, the simulation does not {accurately} capture the {flow} features,  %of the flow, %presenting 
 {displaying} numerical oscillations.
\begin{figure}[H]
  \centering
  \includegraphics[width=0.49\textwidth]{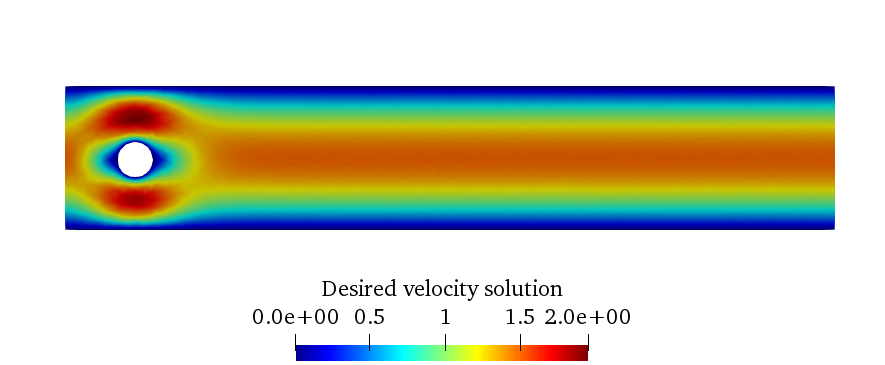}
  \caption{Experiment 1. Desired velocity profile $U$: steady Stokes flow ($\nu$ = 1).}
  \label{fig:U1}
\end{figure}
From the bottom right panel, it is clear that the action of $f_B$ stabilizes the flow in order to reach the desired configuration $U$. This does not hold for the A-controlled velocity (%see 
{displayed in} the bottom left panel in Figure \ref{fig:nocontrolvscontrol}), where spurious numerical oscillations are %inherited from the uncontrolled system. 
{still displayed}. 

%\\ 

\begin{figure}[H]
  \centering
  \includegraphics[width=0.48\textwidth]{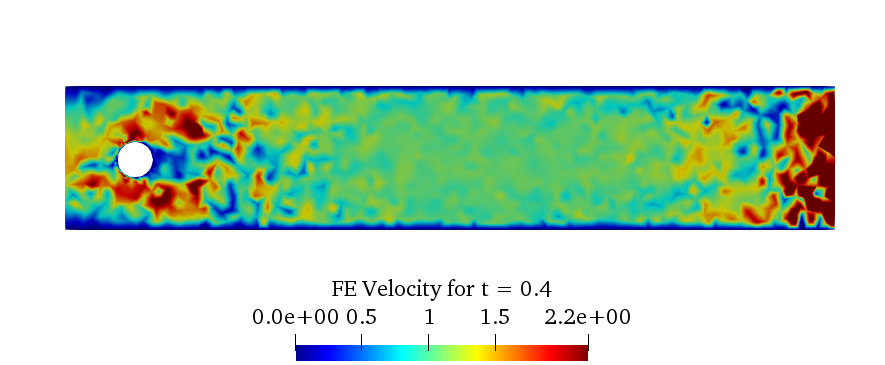}  \\
  \includegraphics[width=0.48\textwidth]{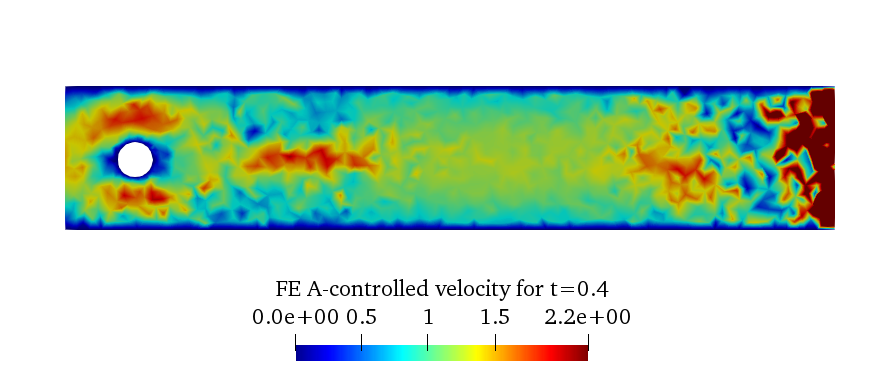}\includegraphics[width=0.48\textwidth]{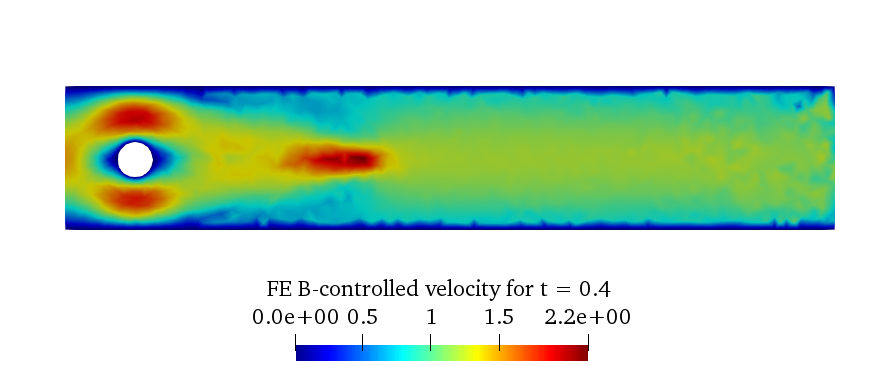}
  \caption{Experiment 1. Top plot: FE velocity profile for $t=0.4$. No regularization is used. Bottom plot: A-controlled {(left)} and B-controlled {(right)} versions of the velocity profile for $t=0.4$ and $\gamma=1$. %, left and right plots respectively.
  }
  \label{fig:nocontrolvscontrol}
\end{figure}
This result agrees with our motivation %of 
{for} introducing $f_B$ as a \emph{stronger} control, especially tailored to treat %more complex phenomena 
{higher Reynolds number flows} with $\nu \ll 1$, while $f_A$ was defined to deal with {laminar} settings where $(\nu - K_0\norm{U}_{L^{\infty}((0,T);H^1(\Omega))})>0$, namely $\nu \sim 1$ {(see Theorem~\ref{theo:gunz1})}. 
In Figure \ref{fig:exp1}, we represent the exponential convergence for {the} A-controlled and B-controlled systems, {which} we %theoretically 
proved in Section \ref{sec:problem}.
We show the value of {the \emph{tracking error},} $E_U(t) = \norm{u(t) -U(t)}_{L^2(\Omega)}^2$, for $\gamma \in \{5, 1\}$ and $t \in (0, 8)$ in the right plot of Figure \ref{fig:exp1}. 
Similarly, the left panel shows the temporal evolution of $E_U(t)$ when $\gamma \in \{50, 25\}$ and $t \in (0, 0.3)$, where the considerably shorter time interval is justified by the fact that convergence is reached way before the final time $T=8$. Indeed, for $\gamma = 50$, $E_U(t)$ is below machine precision already for $t=0.324$. For {the} large values of $\gamma \in \{50, 25\}$, the two controls are comparable in terms of convergence. Differences between the A-controlled and B-controlled cases appear %instead 
for smaller values of $\gamma$, i.e.,  $\gamma=5$ and $\gamma=1$. While for $\gamma = 5$ the exponential convergence is preserved for both %the 
controls, we can observe that {\it $f_A$ does not converge to $U$ for $\gamma=1$, whereas $f_B$ quickly converges to $U$}. For these reasons, from now on, we will investigate the performances of the $f_B$ control, {which} %that %seems to be 
{is} more suited to problems where $\gamma \ll 1$, i.e., when the control action might not {be sufficient to stabilize the simulation in the convective-dominated regime.} %and the balance the convective-dominated nature of the Navier-Stokes equations for $Re=1000$. 
%\ti{Maybe clarify even more? :)} \ms{rephrased}

\begin{figure}[H]
  \centering
  \includegraphics[width=0.5\textwidth]{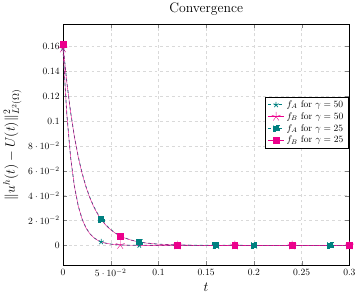}
  \includegraphics[width=0.47\textwidth]{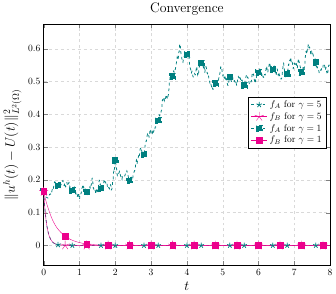}
  \caption{Experiment 1. Convergence results for $f_A$ (teal dashed line) and $f_B$ (magenta solid line) for %s
  $\gamma \in \{50,25, 5, 1\}.$ For $\gamma=50, 25, 5$, the two approaches coincide. }
  \label{fig:exp1}
\end{figure}

\subsection{noEFR vs.\ EFR and aEFR numerical comparison (Experiment 2)}
\label{sec:examplecontrolFOM}
Following the discussion above, {in what follows}, we only consider the control $f_B$, %since its definition takes into account the possibility of dealing with high $Re$.
{which is more appropriate for high $Re$ flows.}
%\ti{We should say this explicitly in the previous Section.} \ms{I added that}
A natural question arises: What happens to $f_B$ as $\gamma \rightarrow 0$, namely for a weak control action, {which is common in realistic applications}. %\ti{i.e., in the low control limit?)} \ms{Yes, but I am not sure about the terminology ``low control" limit. I used, ``weak control action"}. 

\begin{figure}
  \centering
  \includegraphics[width=0.43\textwidth]{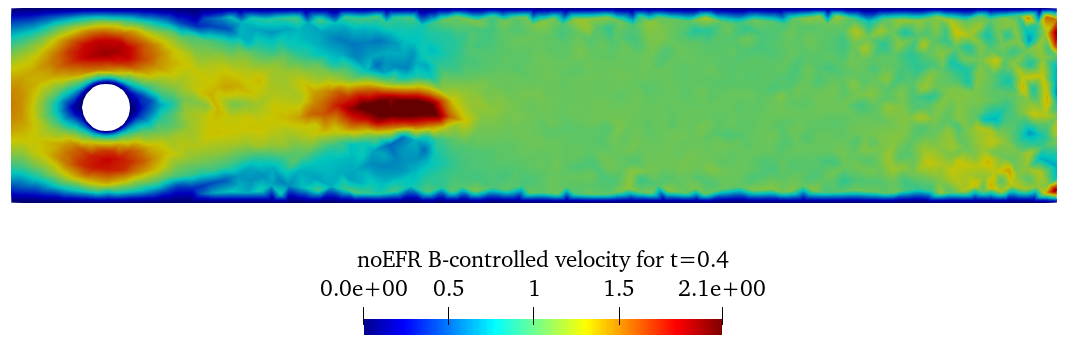}\includegraphics[width=0.43\textwidth]{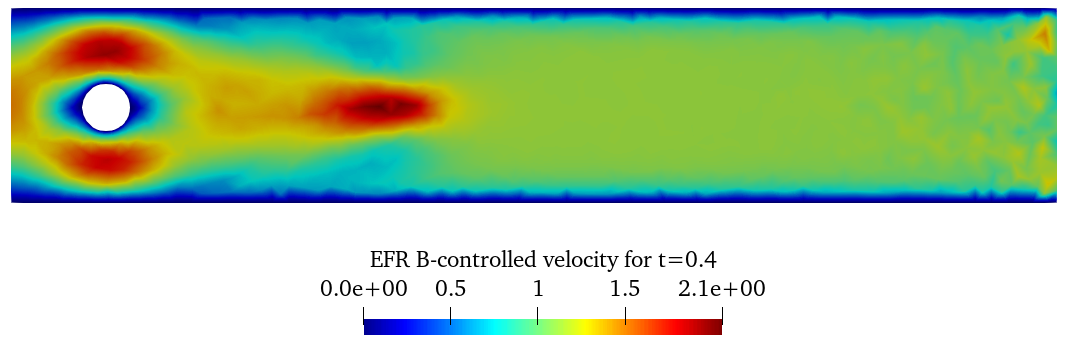}\\
  \includegraphics[width=0.43\textwidth]{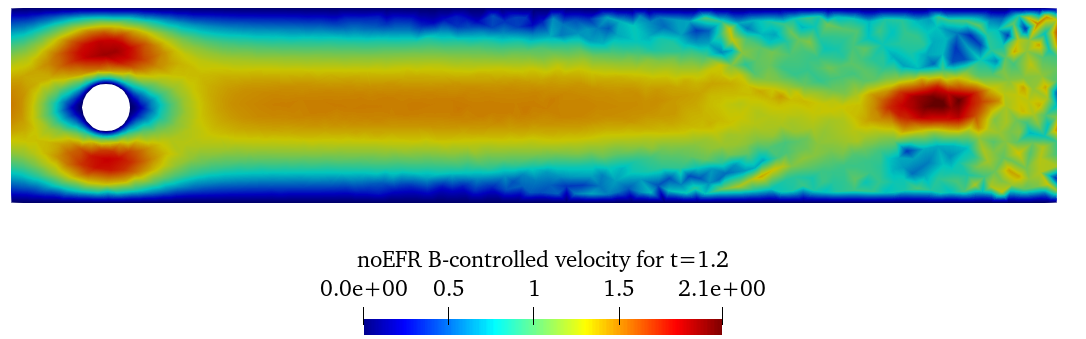}  
\includegraphics[width=0.43\textwidth]{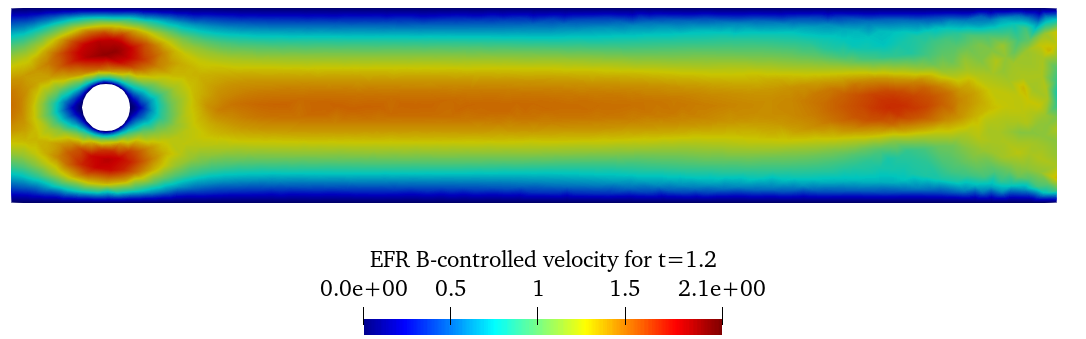}\\ \includegraphics[width=0.43\textwidth]{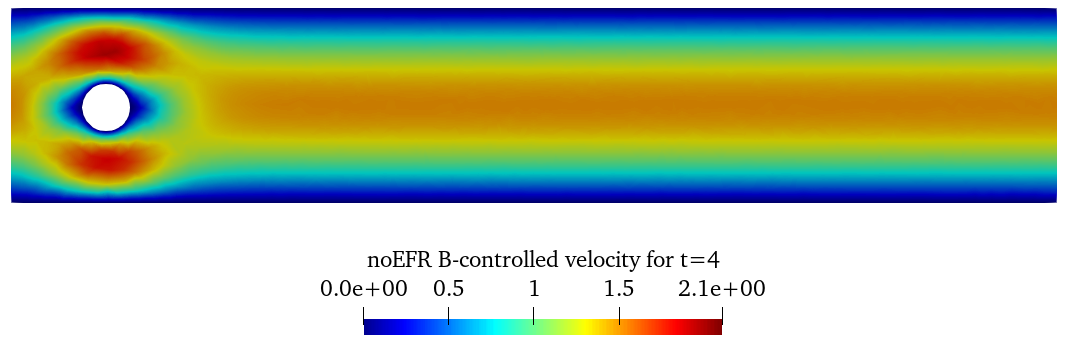}
\includegraphics[width=0.43\textwidth]{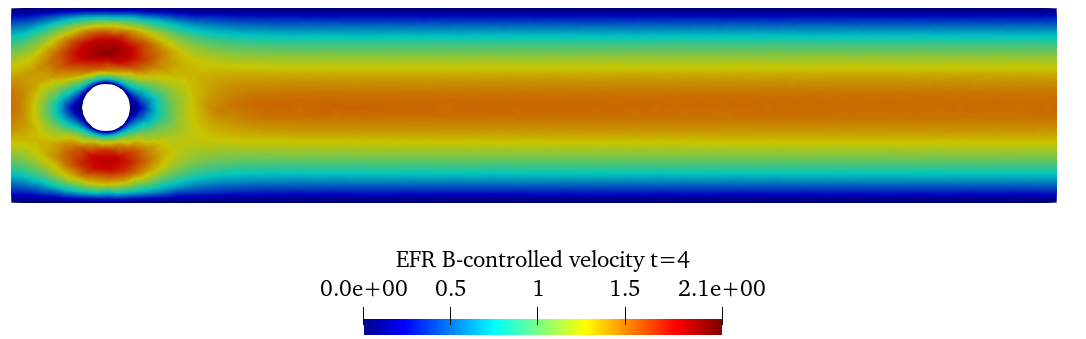}
  \caption{Experiment 2. Representative noEFR (left plots) and EFR (right plots) B-controlled velocity for $t\in \{ 0.4, 1.2, 4\}$ (top to bottom) and $\gamma=10^{-4}$. No qualitative difference can be observed between EFR and aEFR strategy; %: 
  {for clarity,} the latter %approach is 
  {results are} not shown.}
  \label{fig:compEFRaEFR}
\end{figure}
We investigate this question in the same numerical setting proposed in Section \ref{sec:experiment1}, with the only exception that $T=4$. %Indeed 
{The reason for this choice is that} 
the simulation is trivial for $t > 4$, when the behaviour of the controlled system is laminar and similar to the desired state $U$, which numerical experiments in the previous Section show to be reached already for $t \sim 2$. We show the results in Figures  \ref{fig:compEFRaEFR} and \ref{fig:convergencegamma0001} for $\gamma=10^{-4}$, which is a %value of 
$\gamma$ {value} considerably smaller than the %ones 
{values} employed in the previous experiment. It is clear that %for 
such a small {$\gamma$} value %$\gamma=10^{-4}$ 
%we are slowing 
{slows} down the convergence of the controlled solution $u$ to the desired state $U$ with respect to the cases in %the previous experiment.
{Experiment 1.}
From a qualitative point of view, such a conclusion is supported by a comparison of the top left plot of Figure \ref{fig:compEFRaEFR} with the bottom right plot of Figure \ref{fig:nocontrolvscontrol}, representing B-controlled velocities at time instance $t=0.4$ for $\gamma=10^{-4}$ and $\gamma=1$, respectively. The larger value $\gamma=1$ stabilizes numerical oscillations, which are instead more visible when $\gamma=10^{-4}$, most of all at the end of the channel. These oscillations propagate in time within the channel (see the left center plot of Figure \ref{fig:compEFRaEFR}) until the convergence is reached at $t=4$ (see the left bottom plot of Figure \ref{fig:compEFRaEFR}). The spurious oscillations related to the small values of $\nu$ and $\gamma$ affect the exponential convergence of the problem, as depicted in Figure \ref{fig:convergencegamma0001}. Let us focus on the left plot \reviewerA{depicting the derivative of $E_U(t)$ with respect to time}: \reviewerA{observing the oscillating behaviour of the derivative, it is clear that the exponential {convergence} of the noEFR approach is affected for $t \in (0,2)$ by the spurious numerical noise, slowing down the achievement of the goal}.
This suggests that, there is a need for some other kind of stabilization, {different from the control action}, to get rid of the spurious numerical instabilities and reach the goal within a certain threshold in a %smaller amount of 
{shorter} time. We employ the EFR strategy and analyze its effects on the control setting. For the EFR {approach}, we use the same parameters %used in Section \ref{sec:exampleFOM}, i.e.,
$\delta = C_\delta \delta^*$ with $C_{\delta}=\sqrt{11}$, %and 
$\delta^* = 4.46\cdot 10^{-3}$, and $\chi=5 \cdot \Delta t$. %The reason behind this choice \ms{is related to a comparison between noEFR $L^2$ velocity norm for a fine mesh and EFR applied to the mesh of Figure \ref{fig:mesh} in the uncontrolled setting. However, for the sake of brevity, we omit those numerical results}.
{In the uncontrolled setting, these parameters yielded the most accurate results on the mesh in Figure \ref{fig:mesh} (numerical results not included). 
In the controlled framework, we employ the same parameters.}
%is %explained in Section \ref{sec:exampleFOM} for the uncontrolled setting. In the controlled framework, we %decide to 
%employ the same parameters. \\

In the right plots of Figure \ref{fig:compEFRaEFR}, we see how the EFR strategy reduces the numerical oscillations, as expected. We stress that this is also helpful in the convergence behaviour. Indeed, lower values of the objective $E_U$ are reached in a %smaller amount of 
{shorter} time, as can be observed from the dashed orange line of Figure \ref{fig:convergencegamma0001}.
\begin{figure}[H]
  \centering
  \includegraphics[width=0.482\textwidth]{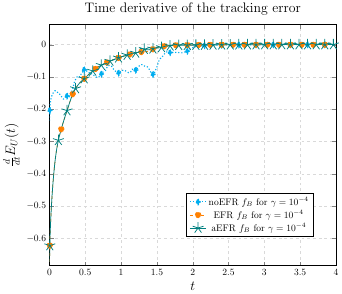}
  \includegraphics[width=0.475\textwidth]{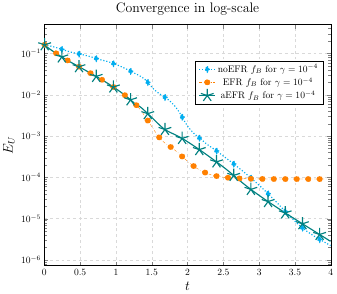}
  \caption{\reviewerA{Experiment 2. Time derivative of $E_U(t)$ (left plot) and log-scale convergence results (right plot) for noEFR (dotted cyan line), EFR (orange dashed line), and aEFR (solid teal line) $f_B$ controlled  velocity with $\gamma=10^{-4}$.}}
  \label{fig:convergencegamma0001}
\end{figure}
 We remark that, {although} %also in the EFR case, 
 we do not recover the exponential decay {in the EFR case, we do obtain significantly better results than the noEFR results:} %, however:
\begin{itemize}
    \item[$\circ$]The EFR approach alleviates the spurious numerical oscillations. Thus, we can expect a better convergence rate with respect to noEFR in the first part of the time interval. %However, 
    {We note, however, that} the exponential convergence of Theorem \ref{theo:our2} is not guaranteed, as shown in Theorem \ref{theo:our5}.
    \item[$\circ$] The EFR strategy reaches the steady state profile with an accuracy of $O(10^{-4})$ more than %a half-second 
    {$0.5$ time units} 
    before the noEFR approach, i.e., {at} $t=2.5$ instead of $t=2.9$.
\end{itemize}
Despite the increased accuracy in the first part of the simulation, the EFR control reaches a plateau for $t>2.5$, and $E_U \sim 10^{-4}$ after that time instance. Indeed, once the convergence is reached, i.e., when the profile is laminar, the EFR is unnecessarily over-diffusing. This leads to a more inaccurate reconstruction with respect to the noEFR approach.
However, for small $\gamma$ and $\nu$ values, noEFR does not show an exponential convergence at the beginning of the simulation, while EFR does. The %issue of the 
EFR stagnation {issue} is overcome through the employment of {the new} aEFR approach with $\tau \sim 0.006$. The choice of the $\tau$ {value} is guided by the EFR simulation: We %wanted to stop 
{stopped} the EFR simulation before {it reached} the plateau. % phenomenon). \\

By definition, aEFR performs as EFR until $E_U \geq \tau$, i.e., for $t< 1.46$. After that time instance, the standard controlled Navier-Stokes simulation is performed. For the sake of brevity, we do not show the representative solutions for the aEFR strategy, since {there are} no {visible} qualitative differences %can be seen 
from the right plots of Figure \ref{fig:compEFRaEFR}. The advantages of using aEFR are clear in the right plot of Figure \ref{fig:convergencegamma0001} (solid teal line), where the exponential convergence is recovered as in the noEFR system avoiding the plateau effect of the EFR strategy. The aEFR approach guarantees $E_U \sim 10^{-4}$  for $t = 2.65$ (before the noEFR approach, which reaches %such 
a {similar} threshold around $t=2.9$). For $t>2.65$, the {aEFR strategy achieves} exponential convergence, %is recovered 
and by $T=4$ noEFR and aEFR are comparable. %\\
\\
\reviewerA{In Table \ref{tab:cpu}, we show the CPU time related to the noEFR, EFR, and aEFR simulations. The regularized strategies do not represent a huge computational burden to deal with when compared with the noEFR strategy. This is due to the modularity and the simplicity of the EFR strategy which (i) solves a linear problem to filter the solution and (ii) performs a convex combination between the evolved and filtered velocity. From Table \ref{tab:cpu}, we can observe that EFR and aEFR are even faster than noEFR: the regularization helps the convergence of the Newton solver, decreasing the computational time of the simulation.}\\
Based on these results, %from now on 
{in what follows,} 
we will exclusively employ the aEFR strategy. 
\reviewerA{

\begin{table}[]
\caption{\reviewerA{CPU time (in seconds). Comparison between noEFR, EFR, and aEFR.}}
\label{tab:cpu}
\centering
\begin{tabular}{l|l|l|l|}
\cline{2-4}
                               & noEFR & EFR & aEFR \\ \hline
\multicolumn{1}{|l|}{CPU time} &   11558     &  11236   &  9060    \\ \hline
\end{tabular}
\end{table}
}
\reviewerA{
\begin{remark}
    We stress that the results of our numerical investigation do not depend on the linear nature of the desired state $U$. Indeed, comparable results are obtained if the solution is steered toward the solution $U$ of the steady state incompressible Navier-Stockes equations: nonlinear system
    \begin{equation}
\label{eq:NS_desired}
\begin{cases}
 - \Delta U + (U \cdot \nabla) U + \nabla P = 0 & \text{in }  \Omega, \\
\nabla \cdot U = 0  & \text{in }  \Omega, \\
U = u_D & \text{on } \Gamma_D, \\
\displaystyle -P n + \frac{\partial U}{\partial n} = 0  & \text{on }  \Gamma_N.
\end{cases}
\end{equation}
For the sake of brevity, we do not include the results, since they do not add any additional information with respect to the results presented in Experiment 2.
\end{remark}}

\subsection{aEFR-noEFR vs. aEFR-aEFR %and  
Numerical Comparison (Experiment 3)}
\label{sec:examplecontrolROM}
In this section, we compare two approaches:
(i) aEFR-noEFR, i.e., {the} standard G-ROM with no regularization at the ROM level; and 
(ii) aEFR-aEFR, which %is 
{employs} the %standard G-ROM with an 
aEFR regularization at the ROM level.
In the FOM setting, we perform the aEFR strategy, as presented in Section \ref{sec:examplecontrolFOM}. All the parameters of the FOM discretization are defined in Section \ref{sec:experiment1}. Furthermore, for the aEFR-aEFR {strategy}, the parameters used at the FOM level are %exploited also 
{also used} at the ROM level, in agreement with the FOM-ROM consistency strategy proposed in \cite{Strazzullo20223148}. \reviewerA{Namely, we employ $\nu=10^{-4}$ (i.e., $Re=1000$), $\gamma=10^{-4}, T=4, C_{\delta}=\sqrt{11}, \delta^*=4.46\cdot 10^{-3}$, and $\chi=5\cdot \Delta t$, with $\Delta t = 4 \cdot 10^{-4}$.}
For the ROM, we %applied 
{apply} a POD approach over $1000$ equispaced snapshots in {the time interval}  $(0,4)$. The ROM solution is represented by $r_u = 20$ basis functions for the velocity and $r_p=r_s=1$ basis functions for the pressure and the supremizer. These values follow a retained %energy 
{information} criterion. 
In Table \ref{tab:E3energy}, we observe that $r_u=20$ preserves $98\%$ of the velocity %energy 
{information} of the system, % $98\%$, %i.e.\ 
{which is} the value we choose as a retained %energy 
{information} threshold. For the pressure, the value $r_p = 1$ is already enough to %describe the 
{represent} $99\%$ of the pressure retained %energy. \ti{What is the pressure energy?} 
{information}. 
For the sake of completeness, {in Table \ref{tab:E3energy},} we report %there 
the retained {information} values for $r_u$ in the set $\{1, 10, 15, 20\}$. %\\

We test the accuracy of the method by means of the relative errors for velocity and pressure %for 
{at} each time instance, defined, respectively, as

\begin{equation}
\label{eq:err_ROM}
E_u(t^n) = \frac{\norm{u^{n} - u_r^n}_{L^2(\Omega)}^2}{\norm{u^n}_{L^2(\Omega)}^2} \quad \text{ and } \quad E_p(t^n) = \frac{\norm{p^{n} - p_r^n}_{L^2(\Omega)}^2}{\norm{p^n}_{L^2(\Omega)}^2},
\end{equation}
 where $u^n$ is the FOM aEFR velocity at time $t_n$, and $u_r^n$ is the reduced velocity (either aEFR-noEFR or aEFR-aEFR). %, and similar 
{Similar} notations are employed for the pressure. %\\

In Figure \ref{fig:exp3comp}, we plot the log-relative errors $E_u$ (left plot) and $E_p$ (right plot) for the aEFR-noEFR and aEFR-aEFR approaches. The employment of aEFR strategy at the ROM level is beneficial in terms of accuracy for both variables. Focusing on the velocity, the gain in accuracy reaches two orders of magnitude for $t=3$, where $E_u$ reaches $10^{-3}$ and $10^{-5}$ for aEFR-noEFR and aEFR-aEFR, respectively. The same happens with the pressure %for 
{at} $t=1.5$, where $E_p$ reaches $10^{2}$ and $1$ for aEFR-noEFR and aEFR-aEFR, respectively. %\\ 
For the velocity, %for 
{at} all the time instances, the aEFR-aEFR is more accurate. We can state the same for the pressure, except for $t \in (0.9, 1.3)$. \\
\reviewerA{The benefits of using aEFR-aEFR are suggested in Table \ref{tab:rom_table_with_gamma}, where we plot $E_u$ and $E_p$ averaged in time for $r_u \in \{5,10,15,20\}$ and $\gamma \in \{10^{-1},10^{-2},10^{-3},10^{-4},10^{-5}\}$. We remark that we leave $r_p=r_s=1$ fixed since one basis function was sufficient to recover the 99\% of the pressure energy. We remark that the snapshots are collected for $\gamma=10^{-4}$. From Table \ref{tab:rom_table_with_gamma}, we conclude that aEFR-aEFR gives good results in terms of extrapolation of controls for all the $\gamma$ and $r_u$ values with errors ranging between $10^{-2}$ and $10^{-3}$. This also holds true for the aEFR-noEFR setting.\\
We observe a different behavior for the pressure: both strategies fail to recover the pressure field with blow-up when using a small value of $r_u$. However, aEFR-aEFR still performs better than aEFR-noEFR in the pressure context, reducing the errors by one order of magnitude for $r_u=15,20$ and all $\gamma$ value. It is also clear that for $\gamma \neq 10^{4}$ aEFR-aEFR and aEFR-noEFR are comparable in terms of velocity field reconstruction.}
}
\begin{figure}
  \centering
  \includegraphics[width=0.49\textwidth]{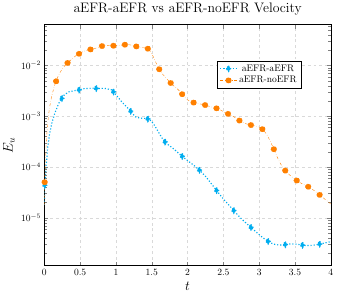}
  \includegraphics[width=0.49\textwidth]{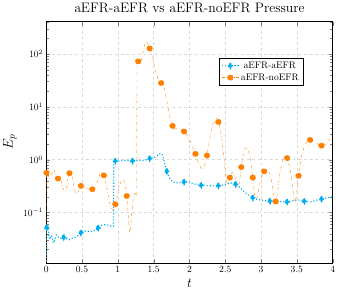}
  \caption{Experiment 3. Relative errors between aEFR-noEFR (dotted cyan line) and aEFR-aEFR (dashed orange line) $f_B$ controlled velocity {(left)} and pressure {(right)} with $\gamma=10^{-4}$. %, right and left plot, respectively
  }
  \label{fig:exp3comp}
\end{figure}
 \begin{table}[H]
\caption{Experiment 3. %Energy 
{Information} retained %from 
{by} the POD basis functions for velocity and pressure.}
\label{tab:E3energy} \centering
\begin{tabular}{|c|c|c|}
\hline
$r_u = r_p$ & Velocity information            & Pressure information               \\ \hline
1  & $53\%$  & $99\%$    \\ \hline
10   & $95\%$  & $99\%$    \\ \hline
15   & $97\%$ & $99\%$     \\ \hline
20  & $98\%$ & $99\%$   \\ \hline
\end{tabular}
\end{table}
\begin{figure}
  \centering
  \includegraphics[width=0.44\textwidth]{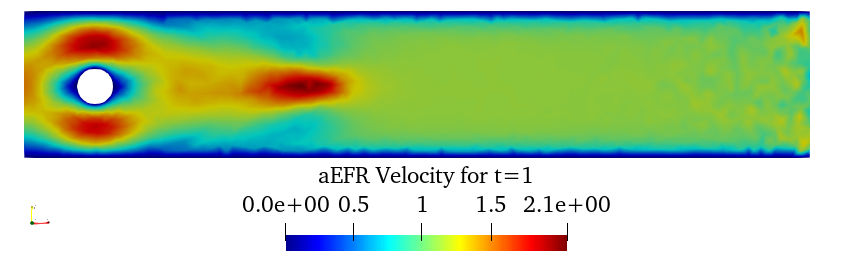}  \\
  \includegraphics[width=0.44\textwidth]{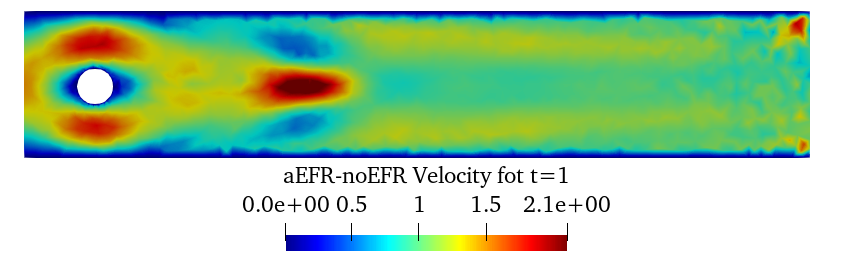}\includegraphics[width=0.44\textwidth]{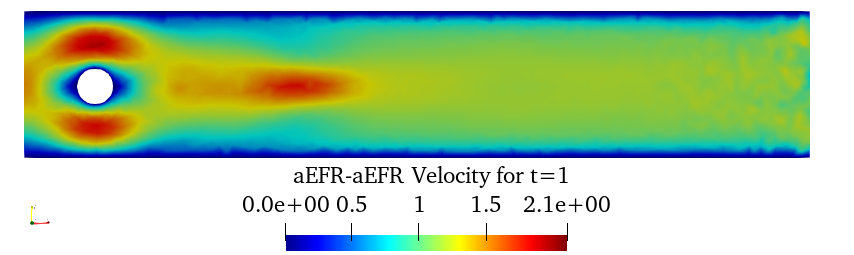}
  \caption{Experiment 3. Top plot: aEFR velocity profile for $t=1$. Bottom plots: aEFR-noEFR {(left)} and aEFR-aEFR {(right)} velocity profiles for $t=1$ and $\gamma=10^{-4}$. %, left and right plots respectively.
  }
  \label{fig:exp3vcomp}
\end{figure}
\begin{figure}
  \centering
  \includegraphics[width=0.44\textwidth]{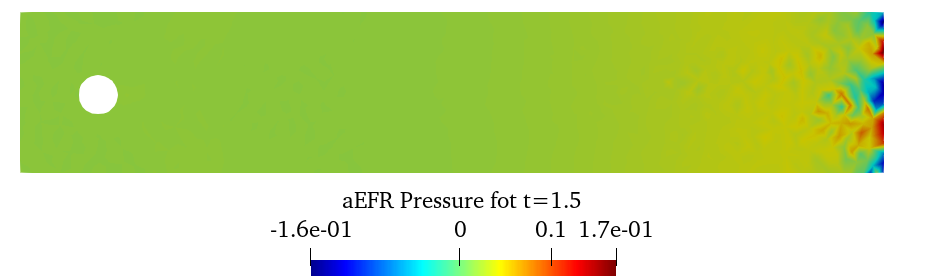}  \\
  \includegraphics[width=0.44\textwidth]{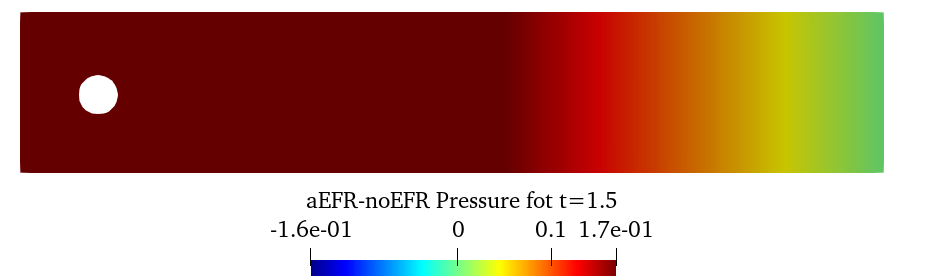}\includegraphics[width=0.44\textwidth]{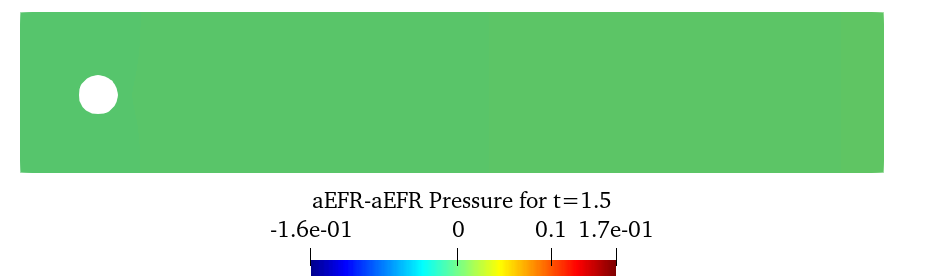}
  \caption{Experiment 3. Top plot: aEFR pressure profile for $t=1$. Bottom plot: aEFR-noEFR {(left)} and aEFR-aEFR {(right)} pressure profiles for $t=1$ and $\gamma=10^{-4}$. %, left and right plots respectively. 
  We stress that the values of the aEFR-aEFR pressure profile %is 
  {are} not exactly zero, but %the values 
  range between $10^{-3}$ and $10^{-2}$. %\fb{move this sentence to the caption?}. \ms{Done}
  }
  \label{fig:exp3pcomp}
\end{figure}

% \usepackage{multirow}
\begin{table}[]
\caption{\reviewerA{Average relative errors for aEFR-noEFR and aEFR-aEFR. The velocity and pressure errors are shown with respect to the basis number $r_u$ and the parameter $\gamma$.}}
\label{tab:rom_table_with_gamma}
\resizebox{\textwidth}{!}{
\begin{tabular}{cc|cccc||cccc|}
\cline{3-10}
                                                  &           & \multicolumn{4}{c||}{Average $E_u$}                                                                           & \multicolumn{4}{c|}{Average $E_p$}                                                                           \\ \cline{2-10} 
\multicolumn{1}{c|}{}                             & $\gamma$  & \multicolumn{1}{c|}{$r_u = 5$} & \multicolumn{1}{c|}{$r_u =10$} & \multicolumn{1}{c|}{$r_u =15$} & $r_u =20$ & \multicolumn{1}{c|}{$r_u = 5$} & \multicolumn{1}{c|}{$r_u =10$} & \multicolumn{1}{c|}{$r_u =15$} & $r_u =20$ \\ \hline
\multicolumn{1}{|c|}{\multirow{5}{*}{aEFR-noEFR}} & $10^{-1}$ & \multicolumn{1}{c|}{1.28e-1}          & \multicolumn{1}{c|}{9.48e-3}          & \multicolumn{1}{c|}{7.40e-3}          &     7.29e-3      & \multicolumn{1}{c|}{2.25e+4}          & \multicolumn{1}{c|}{1.14e+1}          & \multicolumn{1}{c|}{2.16e+1}          &       6.95e+0    \\ \cline{2-10} 
\multicolumn{1}{|c|}{}                            & $10^{-2}$ & \multicolumn{1}{c|}{1.61e-2}          & \multicolumn{1}{c|}{1.21e-2}          & \multicolumn{1}{c|}{9.96-3}          &    9.46e-3       & \multicolumn{1}{c|}{4.23e+4}          & \multicolumn{1}{c|}{3.03e+1}          & \multicolumn{1}{c|}{5.56e+1}          &      1.70e+1     \\ \cline{2-10} 
\multicolumn{1}{|c|}{}                            & $10^{-3}$ & \multicolumn{1}{c|}{1.65e-2}          & \multicolumn{1}{c|}{1.24e-2}          & \multicolumn{1}{c|}{1.02e-2}          &    9.69e-3       & \multicolumn{1}{c|}{4.47e+4}          & \multicolumn{1}{c|}{3.22e+0}          & \multicolumn{1}{c|}{5.78e+1}          &      1.75e+1     \\ \cline{2-10} 
\multicolumn{1}{|c|}{}                            & $10^{-4}$ & \multicolumn{1}{c|}{1.30e-2}          & \multicolumn{1}{c|}{9.97e-3}          & \multicolumn{1}{c|}{8.76e-3}          &    7.85e-3       & \multicolumn{1}{c|}{6.42e+3}          & \multicolumn{1}{c|}{9.75e+0}          & \multicolumn{1}{c|}{2.12e+1}         &      1.58e+1     \\ \cline{2-10} 
\multicolumn{1}{|c|}{}                            & $10^{-5}$ & \multicolumn{1}{c|}{1.65e-2}          & \multicolumn{1}{c|}{1.25e-2}          & \multicolumn{1}{c|}{1.03e-2}          &      9.73e-3    & \multicolumn{1}{c|}{4.30e+4}          & \multicolumn{1}{c|}{3.25e+1}          & \multicolumn{1}{c|}{5.78e+1}         &        1.82e+1   \\ \hline \hline
\multicolumn{1}{|c|}{\multirow{5}{*}{aEFR-aEFR}}  & $10^{-1}$ & \multicolumn{1}{c|}{1.02e-2}          & \multicolumn{1}{c|}{8.87e-3}          & \multicolumn{1}{c|}{8.63e-3}          &     8.43e-3      & \multicolumn{1}{c|}{6.40e+4}          & \multicolumn{1}{c|}{8.50e+0}          & \multicolumn{1}{c|}{4.43e+0}          &   2.29e+0        \\ \cline{2-10} 
\multicolumn{1}{|c|}{}                            & $10^{-2}$ & \multicolumn{1}{c|}{1.16e-2}          & \multicolumn{1}{c|}{1.02e-2}          & \multicolumn{1}{c|}{9.94e-3}          &     9.73e-3      & \multicolumn{1}{c|}{1.16e+5}          & \multicolumn{1}{c|}{1.44e+1}          & \multicolumn{1}{c|}{6.24e+0}          &      4.94e+0     \\ \cline{2-10} 
\multicolumn{1}{|c|}{}                            & $10^{-3}$ & \multicolumn{1}{c|}{1.12e-2}          & \multicolumn{1}{c|}{1.04e-2}          & \multicolumn{1}{c|}{1.00e-2}          & 9.89e-3           & \multicolumn{1}{c|}{1.16e+5}          & \multicolumn{1}{c|}{1.43e+1}          & \multicolumn{1}{c|}{6.24e+0}          & 5.09e+0          \\ \cline{2-10} 
\multicolumn{1}{|c|}{}                            & $10^{-4}$ & \multicolumn{1}{c|}{2.43e-3}          & \multicolumn{1}{c|}{1.47e-3}          & \multicolumn{1}{c|}{1.11e-3}          &      8.86e-4     & \multicolumn{1}{c|}{1.58e+4}          & \multicolumn{1}{c|}{1.18e+0}          & \multicolumn{1}{c|}{1.61e+0} & 1.99e-1 \\ 
\cline{2-10}
\multicolumn{1}{|c|}{}                            & $10^{-5}$ & \multicolumn{1}{c|}{1.18e-2}          & \multicolumn{1}{c|}{1.02e-2}          & \multicolumn{1}{c|}{1.01e-2}          &    9.91e-3       & \multicolumn{1}{c|}{1.17e+5}          & \multicolumn{1}{c|}{1.23e+1}          & \multicolumn{1}{c|}{6.26e+0} & 5.12e+0  \\ 
\hline

\end{tabular}
}
\end{table}

\begin{figure}[H]
  \centering
  \includegraphics[width=0.49\textwidth]{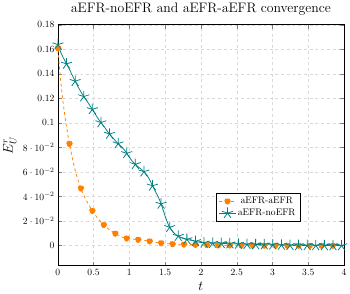}
  \includegraphics[width=0.49\textwidth]{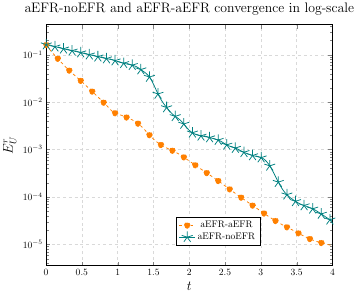}
  \caption{Experiment 3. Convergence results for aEFR-noEFR (solid teal line) and aEFR-aEFR (dashed orange line).  The left plot is %proposed 
  {replotted} %in 
  {on a} log-scale on the right.}
  \label{fig:conv_rom_aefr}
\end{figure}

In order to %state 
{illustrate} the better performance in terms of accuracy of aEFR-aEFR with respect to aEFR FOM solution, we show representative solutions for velocity and pressure in Figure \ref{fig:exp3vcomp} and Figure \ref{fig:exp3pcomp}, respectively. We compare %them 
{aEFR-aEFR and aEFR} with the aEFR-noEFR technique. For the velocity, we analyze {the results at} $t=1$, i.e., {at} a time instance where we can observe an increasing trend of $E_u$ for the aEFR-noEFR velocity. From the plot, we observe that the aEFR-noEFR velocity presents some spurious oscillations in the second half of the channel, i.e., for $x>1$.
{These oscillations are} alleviated by the aEFR-aEFR approach, leading to a solution that is more similar to the aEFR FOM {solution}. %one. \\
For the pressure, we see that the aEFR-noEFR strategy is off with respect to the aEFR solution. The aEFR-aEFR pressure is more similar and, moreover, alleviates the spurious oscillations that are still present at the end of the channel in the FOM simulation.
We conclude that the aEFR-aEFR %approach 
{strategy} is more accurate %with respect to 
{than the} aEFR-noEFR approach. %, underlining 
{We also emphasize} that no extra computational effort is needed {in the aEFR-aEFR strategy} since the filtering and the relaxation steps are not computationally expensive and %they 
are used only for the first part of the simulation. %\\
The benefits of using the aEFR-aEFR {strategy} are evident from the convergence plots in Figure \ref{fig:conv_rom_aefr}, where the {tracking error} $E_U^r = \norm{\mathsf Q_{\mathbb U^{r_{us}}}^Tu^{n}_r - U^n}^2_{L^2(\Omega)}$ is plotted against time. %\\
The aEFR-aEFR strategy stabilizes the system and reaches the desired state $U$ in a %smaller amount of 
{shorter} time, recovering the exponential convergence proved at the FOM level. In particular, we %can observe 
{note} that the aEFR-aEFR velocity reaches $E_U^r \sim 10^{-4}$ at $t=2.75$, almost a %second 
{time unit} before the aEFR-noEFR strategy. 
%In other words, 
{Thus}, {the} aEFR-aEFR approach is not only more accurate but {also} reaches the goal faster. These benefits are not related to an increased algorithmic complexity during the ROM evaluation. \reviewerA{Indeed, for both %the 
techniques the speedup\footnote{\reviewerA{The ROM simulation depends on the FOM dimension
affecting the speedup value. The online computational costs decrease if hyper-reduction techniques, such
as the empirical interpolation method (EIM), are employed. The interested reader may refer to \cite{barrault2004eim} or \cite[Chapter 5]{hesthaven2015certified}.
However, the application of this algorithm was beyond the scope of {this} work.}}, i.e., the number of ROM simulations one can perform in the time of a FOM simulation, is around $1.4$ when $r_u=20$. The speedup slightly increases up to $2$ when we reduce the basis number $r_u$, as expected. Thus, the aEFR-aEFR strategy is competitive with respect to {the} aEFR-noEFR technique in terms of computational costs.} 

\section{Conclusions}
    \label{sec:conclusions}

In this paper, we proposed a novel feedback control strategy for the incompressible Navier-Stokes equations (NSE) at high Reynolds numbers.
For the continuous case, in Theorem~\ref{theo:our2} we proved that the new feedback control strategy yields accurate results for high Reynolds numbers that were not covered by Theorem~\ref{theo:gunz1} %for the standard feedback control approach \ms{Maybe we should say ``
for the feedback law proposed in \cite{GUNZBURGER2000803}. %}.
For the discrete (finite element) case, we %again 
proved in Theorem~\ref{theo:our1} that the new feedback control strategy is accurate for high Reynolds numbers for which Theorem~\ref{theo:gunz3} {in \cite{GUNZBURGER2000803}} does not apply. 
In Section~\ref{exp1:fafb}, we compared the new feedback control strategy with the standard approach in the marginally-resolved numerical simulations of a two-dimensional flow past a circular cylinder at Reynolds numbers $Re=1000$.
The numerical results showed that the new feedback control yields more accurate results than the standard control.

Our second contribution is the development of an adaptive evolve-filter-relax (aEFR) regularization for both the FOM and the ROM feedback control settings.
Our numerical investigation shows that the novel aEFR
%that 
stabilizes marginally-resolved simulations at high Reynolds numbers and increases the accuracy of the new feedback control for realistic control parameters (i.e., when the magnitude of the feedback control is relatively low).
Specifically, the new aEFR strategy yields  more accurate results than the classical EFR approach, especially toward the end of the numerical simulation.  
Furthermore, both aEFR and EFR are significantly more accurate than the noEFR approach (i.e., when no regularization is used), especially at the beginning of the numerical simulation.

The first steps in the theoretical and numerical investigation of the novel feedback control and aEFR algorithm are encouraging. {Furthermore, although the new control law has been built \emph{a priori}, it represents a first step towards more complex applications: for example, they can be interpreted as a first approximation of nonlinear \emph{optimal} feedback controls. This will be a topic of future research}.
\reviewerAA{Furthermore, the law we propose is distributed over all the spatial domain and is impractical for NSE-based applications: a significant (and challenging) improvement would be to extend the theoretical and experimental analysis to boundary control, with \cite{raymond2005local,Raymond2006,RAYMOND2007627} as starting points.}
There are several other research directions for the further study of these strategies.
For example, other regularized ROMs (such as the Leray ROM~\cite{wells2017evolve}) could be analyzed and investigated numerically within the new feedback control framework.
The new feedback control and aEFR algorithm could also be studied in the numerical simulation and control of higher Reynolds number turbulent flows.
\reviewerAA{Furthermore, it would be interesting to compare the proposed feedback law with other feedback control strategies, such as model predictive control, and to test the robustness of the approach in noisy settings and under initial data and state perturbations.}
\reviewerB{Another possible future direction would be the employment of hyper-reduction techniques to achieve an efficient ROM capable of dealing with more complicated problems, such as parametrized NSE with varying Reynolds numbers. This would yield a complete model with the objective of a faster tuning of parameter optimization, capable of predicting unseen configurations, paving the way to adaptive and time-dependent control actions. Another valuable extension concerns partial or scarce knowledge about the target and the state, together with the addition of external source terms, and how these features would impact the control action.
Moreover, a valuable improvement concerns the understanding of the role of pressure and increasing its accuracy in the context of the regularized control strategies.}
{Finally, rigorous mathematical support should be provided for the EFR and aEFR strategies within the new feedback control framework in the ROM setting.}

\section*{Acknowledgments}
We acknowledge the European Union's Horizon 2020 research and innovation program under the Marie Skłodowska-Curie Actions, grant agreement 872442 (ARIA).
MS and CC thank the MIUR project ``Dipartimenti di Eccellenza 2018-2022'' (CUP E11G18000350001). MS and FB acknowledge the INdAM-GNCS project ``Metodi numerici per lo studio di strutture geometriche parametriche complesse'' (CUP E53C22001930001) and MS thanks the INdAM-GNCS Project ``Metodi di riduzione di modello ed approssimazioni di rango basso per problemi alto-dimensionali" (CUP\_E53C23001670001). 
MS acknowledges the ``20227K44ME - Full and Reduced order modelling of coupled systems: focus on non-matching methods and automatic learning (FaReX)" project – funded by European Union – Next Generation EU  within the PRIN 2022 program (D.D. 104 - 02/02/2022 Ministero dell’Università e della Ricerca). This manuscript reflects only the authors’ views and opinions and the Ministry cannot be considered responsible for them.
FB acknowledges the PRIN 2022 PNRR project ``ROMEU: Reduced Order Models for Environmental and Urban flows'' funded by the European Union -- NextGenerationEU under the National Recovery and Resilience Plan (NRRP), Mission 4 Component 2, CUP J53D23015960001. FB also thanks the project ``Reduced order modelling for numerical simulation of partial differential equations'' funded by Università Cattolica del Sacro Cuore. TI acknowledges {support through National Science Foundation grants DMS-2012253 and CDS\&E-MSS-1953113.} 

The computations in this work have been performed with RBniCS \cite{rbnics} library, which is an implementation in FEniCS \cite{fenics} of several reduced order modeling techniques; we acknowledge developers and contributors to both libraries. Computational resources were partially provided by HPC@POLITO, a project of Academic Computing within the Department of Control and Computer Engineering at the Politecnico di Torino (\href{http://hpc.polito.it}{http://hpc.polito.it}).

% \ti{to-do:
% Cite papers on: 
% (i) NS-$\alpha$; 
% (ii) optimal control + Regularized models, e.g., \cite{mallea2021optimal};
% %(iii) Cite Elnaz (ML).
% }
\bibliographystyle{abbrv}
\bibliography{bib_maria.bib,traian.bib,maria.bib}

\begin{thebibliography}{10}

\bibitem{rbnics}
{RBniCS} -- reduced order modelling in {FEniCS},
  \href{https://www.rbnicsproject.org/}{https://www.rbnicsproject.org/}.

\bibitem{GUNZBURGER2000803}
Analysis and approximation for linear feedback control for tracking the
  velocity in {N}avier–{S}tokes flows.
\newblock {\em Computer Methods in Applied Mechanics and Engineering},
  189(3):803--823, 2000.

\bibitem{afs19}
A.~Alla, M.~Falcone, and L.~Saluzzi.
\newblock An efficient {DP} algorithm on a tree-structure for finite horizon
  optimal control problems.
\newblock {\em SIAM Journal of Scientific Computing}, 41(4):2384--2406, 2019.

\bibitem{afs20}
A.~Alla, M.~Falcone, and L.~Saluzzi.
\newblock A tree structure algorithm for optimal control problems with state
  constraints.
\newblock {\em Rendiconti di Matematica e delle Sue Applicazioni}, 41:193--221,
  2020.

\bibitem{Alla20173091}
A.~Alla, M.~Falcone, and S.~Volkwein.
\newblock {Error analysis for POD approximations of infinite horizon problems
  via the dynamic programming approach}.
\newblock 55(5):3091--3115, 2017.

\bibitem{10.1007/978-3-319-23413-7_120}
A.~Alla and M.~Hinze.
\newblock {HJB-POD} feedback control for {N}avier-{S}tokes equations.
\newblock In G.~Russo, V.~Capasso, G.~Nicosia, and V.~Romano, editors, {\em
  Progress in Industrial Mathematics at ECMI 2014}, pages 861--868, Cham, 2016.
  Springer International Publishing.

\bibitem{alla2020}
A.~Alla and L.~Saluzzi.
\newblock A {HJB-POD} approach for the control of nonlinear {PDE}s on a tree
  structure.
\newblock {\em Appl Numer Math.}, 155:192--207, 2020.

\bibitem{audouze2009reduced}
C.~Audouze, F.~De~Vuyst, and P.~Nair.
\newblock Reduced-order modeling of parameterized {PDE}s using
  time--space-parameter principal component analysis.
\newblock {\em International Journal for Numerical Methods in Engineering},
  80(8):1025--1057, 2009.

\bibitem{ballarin2016fast}
F.~Ballarin, E.~Faggiano, S.~Ippolito, A.~Manzoni, A.~Quarteroni, G.~Rozza, and
  R.~Scrofani.
\newblock Fast simulations of patient-specific haemodynamics of coronary artery
  bypass grafts based on a {POD--Galerkin} method and a vascular shape
  parametrization.
\newblock {\em J. Comput. Phys.}, 315:609--628, 2016.

\bibitem{ballarin2015supremizer}
F.~Ballarin, A.~Manzoni, A.~Quarteroni, and G.~Rozza.
\newblock Supremizer stabilization of {POD--G}alerkin approximation of
  parametrized steady incompressible {N}avier--{S}tokes equations.
\newblock {\em Int. J. Numer. Meth. Engng.}, 102:1136--1161, 2015.

\bibitem{StrazzulloZuazua}
F.~Ballarin, G.~Rozza, and M.~Strazzullo.
\newblock Chapter 9 - {S}pace-time {POD}-{Galerkin} approach for parametric
  flow control.
\newblock In E.~Tr\`elat and E.~Zuazua, editors, {\em Numerical Control: Part
  A}, volume~23 of {\em Handbook of Numerical Analysis}, pages 307--338.
  Elsevier, 2022.

\bibitem{BanschBenner2015}
E.~B\"{a}nsch, P.~Benner, J.~Saak, and H.~K. Weichelt.
\newblock Riccati-based boundary feedback stabilization of incompressible
  {N}avier--{S}tokes flows.
\newblock {\em SIAM Journal on Scientific Computing}, 37(2):A832--A858, 2015.

\bibitem{COCV_2003__9__197_0}
V.~Barbu.
\newblock Feedback stabilization of {N}avier-{S}tokes equations.
\newblock {\em ESAIM: Control, Optimisation and Calculus of Variations},
  9:197--205, 2003.

\bibitem{Barbu2004}
V.~Barbu and R.~Triggiani.
\newblock Internal stabilization of {N}avier-{S}tokes equations with
  finite-dimensional controllers.
\newblock {\em Indiana University Mathematics Journal}, 53(5):1443--1494, 2004.

\bibitem{barrault2004eim}
M.~Barrault, Y.~Maday, N.~C. Nguyen, and A.~T. Patera.
\newblock An `empirical interpolation' method: Application to efficient
  reduced-basis discretization of partial differential equations.
\newblock {\em C. R. Acad. Sci. Paris, Ser. I}, 339:667--672, 2004.

\bibitem{Benner20201653}
P.~Benner, S.~Dolgov, A.~Onwunta, and M.~Stoll.
\newblock Low-rank solution of an optimal control problem constrained by random
  {N}avier-{S}tokes equations.
\newblock {\em International Journal for Numerical Methods in Fluids},
  92(11):1653--1678, 2020.

\bibitem{benner2019optimal}
P.~Benner and C.~Trautwein.
\newblock Optimal control problems constrained by the stochastic
  {N}avier--{S}tokes equations with multiplicative {L{\'e}vy} noise.
\newblock {\em Mathematische Nachrichten}, 292(7):1444--1461, 2019.

\bibitem{bertagna2016deconvolution}
L.~Bertagna, A.~Quaini, and A.~Veneziani.
\newblock {Deconvolution-based nonlinear filtering for incompressible flows at
  moderately large Reynolds numbers}.
\newblock {\em Int. J. Num. Meth. Fluids}, 81(8):463--488, 2016.

\bibitem{brands2016reduced}
B.~Brands, J.~Mergheim, and P.~Steinmann.
\newblock Reduced-order modelling for linear heat conduction with parametrised
  moving heat sources.
\newblock {\em GAMM-Mitteilungen}, 39(2):170--188, 2016.

\bibitem{burkardt2006pod}
J.~Burkardt, M.~Gunzburger, and H.~Lee.
\newblock {POD and CVT-based reduced-order modeling of Navier--Stokes flows}.
\newblock {\em Comput. Methods Appl. Mech. Engrg.}, 196(1-3):337--355, 2006.

\bibitem{Collis2004237}
S.~Collis, R.~Joslin, A.~Seifert, and V.~Theofilis.
\newblock Issues in active flow control: theory, control, simulation,' and
  experiment.
\newblock {\em Progress in Aerospace Sciences}, 40(4-5):237--289, 2004.

\bibitem{collis2002analysis}
S.~S. Collis and M.~Heinkenschloss.
\newblock Analysis of the streamline upwind/{P}etrov {G}alerkin method applied
  to the solution of optimal control problems.
\newblock {\em CAAM TR02-01}, 108, 2002.

\bibitem{Deckelnick2004297}
K.~Deckelnick and M.~Hinze.
\newblock Semidiscretization and error estimates for distributed control of the
  instationary {Navier-Stokes} equations.
\newblock {\em Numerische Mathematik}, 97(2):297 – 320, 2004.

\bibitem{Dede2012}
L.~Ded{\`e}.
\newblock Reduced basis method and error estimation for parametrized optimal
  control problems with control constraints.
\newblock {\em Journal of Scientific Computing}, 50(2):287--305, Feb 2012.

\bibitem{dolgov2022data}
S.~Dolgov, D.~Kalise, and L.~Saluzzi.
\newblock Data-driven tensor train gradient cross approximation for
  {H}amilton-{J}acobi-{B}ellman equations.
\newblock {\em arXiv preprint arXiv:2205.05109}, 2022.

\bibitem{ervin2012numerical}
V.~J. Ervin, W.~J. Layton, and M.~Neda.
\newblock Numerical analysis of filter-based stabilization for evolution
  equations.
\newblock {\em SIAM J. Numer. Anal.}, 50(5):2307--2335, 2012.

\bibitem{falcone2023approximation}
M.~Falcone, G.~Kirsten, and L.~Saluzzi.
\newblock Approximation of optimal control problems for the {N}avier-{S}tokes
  equation via multilinear {HJB-POD}.
\newblock {\em Applied Mathematics and Computation}, 442:127722, 2023.

\bibitem{fischer2001filter}
P.~F. Fischer and J.~Mullen.
\newblock Filter-based stabilization of spectral element methods.
\newblock {\em C. R. Acad. Sci. Paris S\'er. I Math.}, 332(3):265--270, 2001.

\bibitem{Fri95}
U.~Frisch.
\newblock {\em Turbulence, The {L}egacy of {A}.{N}.~{K}olmogorov}.
\newblock Cambridge University Press, Cambridge, 1995.

\bibitem{10.1007/978-3-642-59709-1_11}
A.~F{\"u}rsikov, M.~Gunzburger, L.~S. Hou, and S.~Manservisi.
\newblock Optimal control problems for the {N}avier-{S}tokes equations.
\newblock In H.-J. Bungartz, R.~H.~W. Hoppe, and C.~Zenger, editors, {\em
  Lectures on Applied Mathematics}, pages 143--155, Berlin, Heidelberg, 2000.
  Springer Berlin Heidelberg.

\bibitem{novo2}
B.~Garc{\'\i}a-Archilla and J.~Novo.
\newblock Error analysis of fully discrete mixed finite element data
  assimilation schemes for the navier-stokes equations.
\newblock {\em Advances in Computational Mathematics}, 46(4):61, 2020.

\bibitem{novo1}
B.~Garc\'{\i}a-Archilla, J.~Novo, and E.~S. Titi.
\newblock Uniform in time error estimates for a finite element method applied
  to a downscaling data assimilation algorithm for the navier--stokes
  equations.
\newblock {\em SIAM Journal on Numerical Analysis}, 58(1):410--429, 2020.

\bibitem{girfoglio2021pod}
M.~Girfoglio, A.~Quaini, and G.~Rozza.
\newblock {A POD-Galerkin reduced order model for a LES filtering approach}.
\newblock {\em J. Comput. Phys.}, 436:110260, 2021.

\bibitem{Girfoglio20210}
M.~Girfoglio, A.~Quaini, and G.~Rozza.
\newblock Pressure stabilization strategies for a {LES} filtering reduced order
  model.
\newblock {\em Fluids}, 6(9), 2021.

\bibitem{GIRFOGLIO2022105536}
M.~Girfoglio, A.~Quaini, and G.~Rozza.
\newblock A {POD}-galerkin reduced order model for the {N}avier–{S}tokes
  equations in stream function-vorticity formulation.
\newblock {\em Computers \& Fluids}, 244:105536, 2022.

\bibitem{GIRFOGLIO2023112127}
M.~Girfoglio, A.~Quaini, and G.~Rozza.
\newblock A hybrid projection/data-driven reduced order model for the
  {N}avier-{S}tokes equations with nonlinear filtering stabilization.
\newblock {\em Journal of Computational Physics}, 486:112127, 2023.

\bibitem{grimberg2020stability}
S.~Grimberg, C.~Farhat, and N.~Youkilis.
\newblock On the stability of projection-based model order reduction for
  convection-dominated laminar and turbulent flows.
\newblock {\em J. Comput. Phys.}, 419:109681, 2020.

\bibitem{gunzburger2020leray}
M.~Gunzburger, T.~Iliescu, and M.~Schneier.
\newblock A {L}eray regularized ensemble-proper orthogonal decomposition method
  for parameterized convection-dominated flows.
\newblock {\em IMA J. Numer. Anal.}, 40(2):886--913, 2020.

\bibitem{gunzburger2002perspectives}
M.~D. Gunzburger.
\newblock {\em Perspectives in flow control and optimization}.
\newblock SIAM, 2002.

\bibitem{gunzburger1991analysis}
M.~D. Gunzburger, L.~Hou, and T.~P. Svobodny.
\newblock Analysis and finite element approximation of optimal control problems
  for the stationary {N}avier-{S}tokes equations with distributed and {N}eumann
  controls.
\newblock {\em Mathematics of Computation}, 57(195):123--151, 1991.

\bibitem{GunzburgerOptimalControl2000}
M.~D. Gunzburger and S.~Manservisi.
\newblock Analysis and approximation of the velocity tracking problem for
  {N}avier-{S}tokes flows with distributed control.
\newblock {\em SIAM Journal on Numerical Analysis}, 37(5):1481--1512, 2000.

\bibitem{Heinkenschloss2008}
M.~Heinkenschloss, D.~C. Sorensen, and K.~Sun.
\newblock Balanced truncation model reduction for a class of descriptor systems
  with application to the {O}seen equations.
\newblock {\em SIAM Journal on Scientific Computing}, 30(2):1038--1063, 2008.

\bibitem{hesthaven2015certified}
J.~S. Hesthaven, G.~Rozza, and B.~Stamm.
\newblock {\em Certified Reduced Basis Methods for Parametrized Partial
  Differential Equations}.
\newblock Springer, 2015.

\bibitem{himpe2018hierarchical}
C.~Himpe, T.~Leibner, and S.~Rave.
\newblock Hierarchical approximate proper orthogonal decomposition.
\newblock {\em SIAM Journal on Scientific Computing}, 40(5):A3267--A3292, 2018.

\bibitem{doi:10.1137/S036301290241246X}
M.~Hinze.
\newblock Instantaneous closed loop control of the {N}avier--{S}tokes system.
\newblock {\em SIAM Journal on Control and Optimization}, 44(2):564--583, 2005.

\bibitem{Hinze2000273}
M.~Hinze and K.~Kunisch.
\newblock Three control methods for time-dependent fluid flow.
\newblock {\em Flow, Turbulence and Combustion}, 65(3-4):273--298, 2000.

\bibitem{hou}
L.~S. Hou and Y.~Yan.
\newblock Dynamics for controlled {N}avier--{S}tokes systems with distributed
  controls.
\newblock {\em SIAM Journal on Control and Optimization}, 35(2):654--677, 1997.

\bibitem{kaneko2020towards}
K.~Kaneko, P.-H. Tsai, and P.~Fischer.
\newblock Towards model order reduction for fluid-thermal analysis.
\newblock {\em Nucl. Eng. Des.}, 370:110866, 2020.

\bibitem{karcher_grepl_2014}
M.~K{\"a}rcher and M.~A. Grepl.
\newblock A posteriori error estimation for reduced order solutions of
  parametrized parabolic optimal control problems.
\newblock {\em ESAIM: Mathematical Modelling and Numerical Analysis},
  48(6):1615–1638, 2014.

\bibitem{karcher2018certified}
M.~K{\"a}rcher, Z.~Tokoutsi, M.~A. Grepl, and K.~Veroy.
\newblock Certified reduced basis methods for parametrized elliptic optimal
  control problems with distributed controls.
\newblock {\em Journal of Scientific Computing}, 75(1):276--307, 2018.

\bibitem{kunisch2008proper}
K.~Kunisch and S.~Volkwein.
\newblock Proper orthogonal decomposition for optimality systems.
\newblock {\em ESAIM: Math. Model. Numer. Anal.}, 42(1):1--23, 2008.

\bibitem{layton2008numerical}
W.~Layton, C.~C. Manica, M.~Neda, and L.~G. Rebholz.
\newblock Numerical analysis and computational testing of a high accuracy
  {L}eray-deconvolution model of turbulence.
\newblock {\em Num. Meth. P.D.E.s}, 24(2):555--582, 2008.

\bibitem{layton2008introduction}
W.~J. Layton.
\newblock {\em Introduction to the numerical analysis of incompressible viscous
  flows}, volume~6.
\newblock Society for Industrial and Applied Mathematics (SIAM), 2008.

\bibitem{Lee20212533}
H.-C. Lee.
\newblock Efficient computations for linear feedback control problems for
  target velocity matching of navier-stokes flows via pod and lstm-rom.
\newblock {\em Electronic Research Archive}, 29(3):2533 – 2552, 2021.

\bibitem{lee2021efficient}
H.-C. Lee.
\newblock Efficient computations for linear feedback control problems for
  target velocity matching of {N}avier-{S}tokes flows via {POD} and
  {LSTM}-{ROM}.
\newblock {\em Electronic Research Archive}, 29(3):2533--2552, 2021.

\bibitem{Leykekhman20122012}
D.~Leykekhman and M.~Heinkenschloss.
\newblock Local error analysis of discontinuous {G}alerkin methods for
  advection-dominated elliptic linear-quadratic optimal control problems.
\newblock {\em SIAM Journal on Numerical Analysis}, 50(4):2012--2038, 2012.

\bibitem{fenics}
A.~Logg, K.~Mardal, and G.~Wells.
\newblock {\em Automated Solution of Differential Equations by the Finite
  Element Method}.
\newblock Springer-Verlag, Berlin, 2012.

\bibitem{mallea2021optimal}
E.~Mallea-Zepeda, E.~Ortega-Torres, and {\'E}.~J. Villamizar-Roa.
\newblock {An optimal control problem for the Navier-Stokes-$\alpha$ system}.
\newblock {\em J. Dyn. Control Syst.}, pages 1--28, 2021.

\bibitem{mullen1999filtering}
J.~S. Mullen and P.~F. Fischer.
\newblock Filtering techniques for complex geometry fluid flows.
\newblock {\em Commun. Numer. Meth. Engng.}, 15(1):9--18, 1999.

\bibitem{negri2015reduced}
F.~Negri, A.~Manzoni, and G.~Rozza.
\newblock Reduced basis approximation of parametrized optimal flow control
  problems for the {S}tokes equations.
\newblock {\em Computers \& Mathematics with Applications}, 69(4):319--336,
  2015.

\bibitem{negri2013reduced}
F.~Negri, G.~Rozza, A.~Manzoni, and A.~Quarteroni.
\newblock Reduced basis method for parametrized elliptic optimal control
  problems.
\newblock {\em SIAM Journal on Scientific Computing}, 35(5):A2316--A2340, 2013.

\bibitem{noack2011reduced}
B.~R. Noack, M.~Morzynski, and G.~Tadmor.
\newblock {\em Reduced-Order Modelling for Flow Control}, volume 528.
\newblock Springer Verlag, 2011.

\bibitem{parish2023residual}
E.~J. Parish, M.~Yano, I.~Tezaur, and T.~Iliescu.
\newblock Residual-based stabilized reduced-order models of the transient
  convection-diffusion-reaction equation obtained through discrete and
  continuous projection.
\newblock {\em arXiv preprint, \url{http://arxiv.org/abs/2302.09355}}, 2023.

\bibitem{Pichi20221361}
F.~Pichi, M.~Strazzullo, F.~Ballarin, and G.~Rozza.
\newblock Driving bifurcating parametrized nonlinear {PDE}s by optimal control
  strategies: Application to {N}avier-{S}tokes equations with model order
  reduction.
\newblock {\em ESAIM: Mathematical Modelling and Numerical Analysis},
  56(4):1361 – 1400, 2022.

\bibitem{Pop00}
S.~Pope.
\newblock {\em Turbulent flows}.
\newblock Cambridge University Press, Cambridge, 2000.

\bibitem{quarteroni2009numerical}
A.~Quarteroni and S.~Quarteroni.
\newblock {\em Numerical models for differential problems}, volume~2.
\newblock Springer, 2009.

\bibitem{quarteroni2007reduced}
A.~Quarteroni, G.~Rozza, and A.~Quaini.
\newblock Reduced basis methods for optimal control of advection-diffusion
  problems.
\newblock In {\em Advances in Numerical Mathematics}, number
  CMCS-CONF-2006-003, pages 193--216. RAS and University of Houston, 2007.

\bibitem{quarteroni2008numerical}
A.~Quarteroni and A.~Valli.
\newblock {\em Numerical approximation of partial differential equations},
  volume~23.
\newblock Springer Science \& Business Media, 2008.

\bibitem{raymond2005local}
J.-P. Raymond.
\newblock Local boundary feedback stabilization of the {N}avier-{S}tokes
  equations.
\newblock {\em Control Systems: Theory, Numerics and Applications, Rome}, 30,
  2005.

\bibitem{Raymond2006}
J.-P. Raymond.
\newblock Feedback boundary stabilization of the two-dimensional
  {N}avier--{S}tokes equations.
\newblock {\em SIAM Journal on Control and Optimization}, 45(3):790--828, 2006.

\bibitem{RAYMOND2007627}
J.-P. Raymond.
\newblock Feedback boundary stabilization of the three-dimensional
  incompressible {N}avier–{S}tokes equations.
\newblock {\em Journal de Mathématiques Pures et Appliquées}, 87(6):627--669,
  2007.

\bibitem{rozza2007stability}
G.~Rozza and K.~Veroy.
\newblock On the stability of the reduced basis method for {S}tokes equations
  in parametrized domains.
\newblock {\em Comput. Methods Appl. Mech. Engrg.}, 196(7):1244--1260, 2007.

\bibitem{sabetghadam2012alpha}
F.~Sabetghadam and A.~Jafarpour.
\newblock $\alpha$ regularization of the {POD-G}alerkin dynamical systems of
  the {K}uramoto--{S}ivashinsky equation.
\newblock {\em Appl. Math. Comput.}, 218(10):6012--6026, 2012.

\bibitem{ST96}
M.~Sch$\ddot{\mbox{a}}$fer and S.~Turek.
\newblock The benchmark problem `flow around a cylinder' flow simulation with
  high performance computers {II}.
\newblock {\em in E.H. Hirschel (Ed.), Notes on Numerical Fluid Mechanics}, 52,
  Braunschweig, Vieweg:547--566, 1996.

\bibitem{stabile2019reduced}
G.~Stabile, F.~Ballarin, G.~Zuccarino, and G.~Rozza.
\newblock A reduced order variational multiscale approach for turbulent flows.
\newblock {\em Adv. Comput. Math.}, pages 1--20, 2019.

\bibitem{Strazzullo2}
M.~Strazzullo, F.~Ballarin, and G.~Rozza.
\newblock {POD}-{G}alerkin model order reduction for parametrized time
  dependent linear quadratic optimal control problems in saddle point
  formulation.
\newblock {\em Journal of Scientific Computing}, 83(55), 2020.

\bibitem{StrazzulloRB}
M.~Strazzullo, F.~Ballarin, and G.~Rozza.
\newblock A {C}ertified {R}educed {B}asis method for linear parametrized
  parabolic optimal control problems in space-time formulation.
\newblock Submitted, 2021,
  \href{https://arxiv.org/abs/2103.00460}{https://arxiv.org/abs/2103.00460}.

\bibitem{Strazzullo3}
M.~Strazzullo, F.~Ballarin, and G.~Rozza.
\newblock {POD}-{G}alerkin model order reduction for parametrized nonlinear
  time dependent optimal flow control: an application to {S}hallow {W}ater
  {E}quations.
\newblock {\em Journal of Numerical Mathematics}, 30(1):63--84, 2022.

\bibitem{Strazzullo20223148}
M.~Strazzullo, M.~Girfoglio, F.~Ballarin, T.~Iliescu, and G.~Rozza.
\newblock Consistency of the full and reduced order models for
  evolve-filter-relax regularization of convection-dominated,
  marginally-resolved flows.
\newblock {\em International Journal for Numerical Methods in Engineering},
  123(14):3148--3178, 2022.

\bibitem{ZakiaMaria}
M.~Strazzullo, Z.~Zainib, F.~Ballarin, and G.~Rozza.
\newblock Reduced order methods for parametrized nonlinear and time dependent
  optimal flow control problems: towards applications in biomedical and
  environmental sciences.
\newblock {\em Numerical Mathematics and Advanced Applications ENUMATH 2019},
  2021.

\bibitem{temam2001navier}
R.~Temam.
\newblock {\em Navier-Stokes equations: {T}heory and numerical analysis},
  volume~2.
\newblock American Mathematical Society, 2001.

\bibitem{wells2017evolve}
D.~Wells, Z.~Wang, X.~Xie, and T.~Iliescu.
\newblock An evolve-then-filter regularized reduced order model for
  convection-dominated flows.
\newblock {\em Int. J. Num. Meth. Fluids}, 84:598--–615, 2017.

\bibitem{xie2018numerical}
X.~Xie, D.~Wells, Z.~Wang, and T.~Iliescu.
\newblock Numerical analysis of the {L}eray reduced order model.
\newblock {\em J. Comput. Appl. Math.}, 328:12--29, 2018.

\bibitem{zoccolan2}
F.~Zoccolan, M.~Strazzullo, and G.~Rozza.
\newblock Stabilized weighted reduced order methods for parametrized
  advection-dominated optimal control problems governed by partial differential
  equations with random inputs.
\newblock {\em Journal of Numerical Mathematics}, 2024.

\bibitem{zoccolan1}
F.~Zoccolan, M.~Strazzullo, and G.~Rozza.
\newblock A streamline upwind petrov-galerkin reduced order method for
  advection-dominated partial differential equations under optimal control.
\newblock {\em Computational Methods in Applied Mathematics}, 2024.

\end{thebibliography}

\appendix

\section{Uncontrolled snapshots collection}
\label{rem:uncontrolled} 
In Experiment 3, we considered \emph{controlled} snapshots to build the reduced spaces. However, another way to proceed is to collect EFR snapshots from the uncontrolled problem \eqref{eq:NSE} and apply the control strategy only online. To test this approach, we collect $1000$ equispaced \emph{uncontrolled} snapshots, employing $\nu=10^{-4}$ (i.e., $Re=1000$), $T=4, C_{\delta}=\sqrt{11}, \delta^*=4.46\cdot 10^{-3}$, and $\chi=5\cdot \Delta t$, with $\Delta t = 4 \cdot 10^{-4}$.

\begin{figure}[H]
        \centering \includegraphics[width=0.49\textwidth]{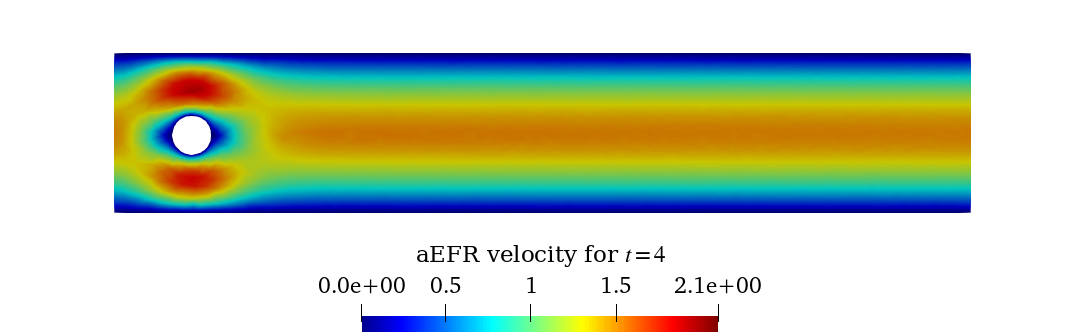}\\
        \centering \includegraphics[width=0.49\textwidth]{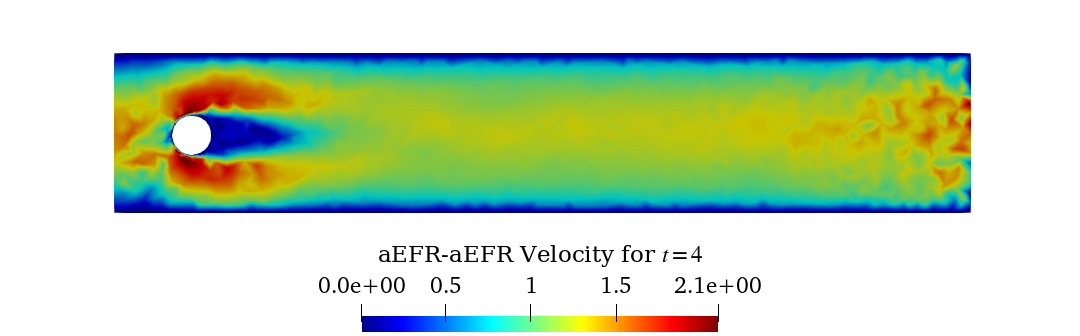}
        \centering \includegraphics[width=0.49\textwidth]{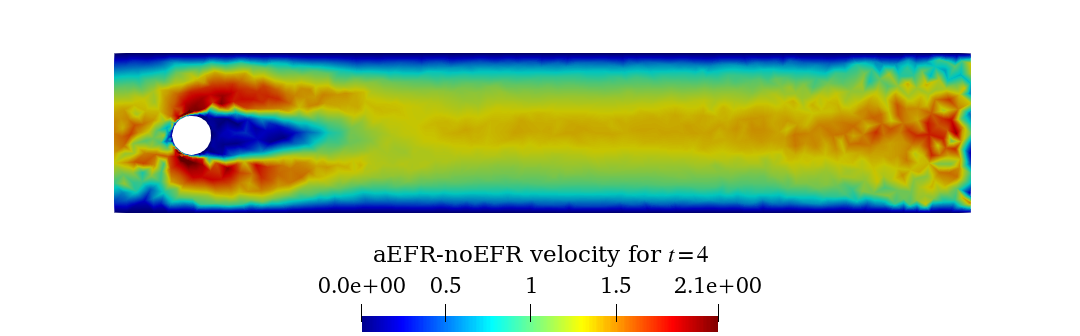}
        
\caption{{Uncontrolled snapshot collection. Top plot: aEFR velocity profile for $t=4$. Bottom plots: aEFR-noEFR {(left)} and aEFR-aEFR {(right)} velocity profiles for $t=4$ and $\gamma=10^{-4}$.}}

\label{fig:rem_comp}
\end{figure}

\begin{figure}[H]
        \centering \includegraphics[scale=0.195]{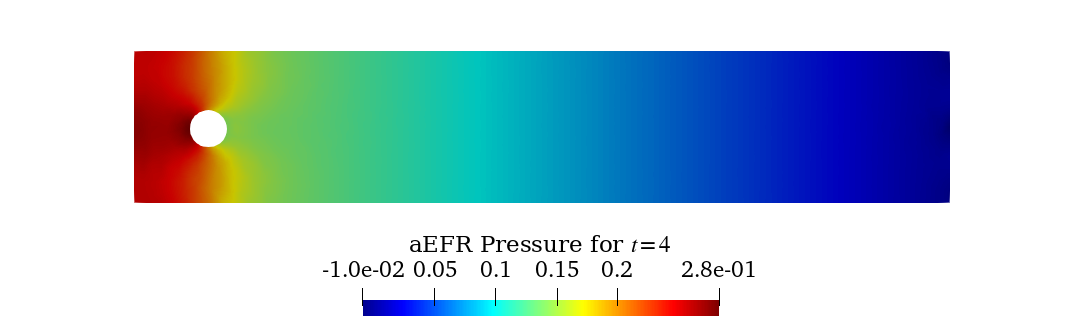}\\
        \centering \includegraphics[scale=0.195]{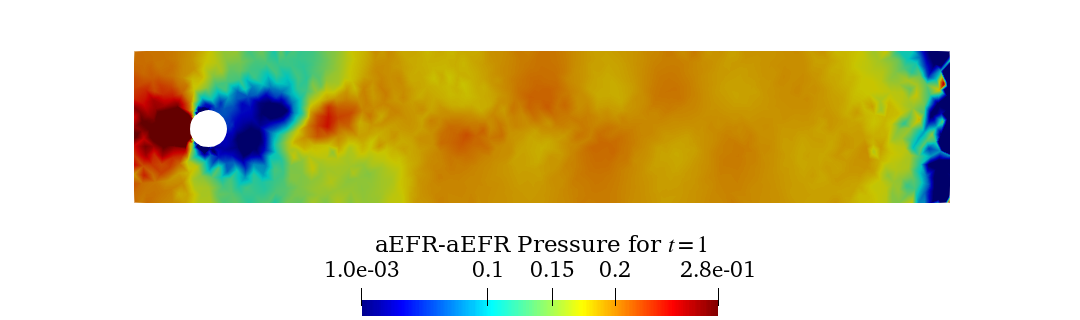}
        \centering \includegraphics[scale=0.195]{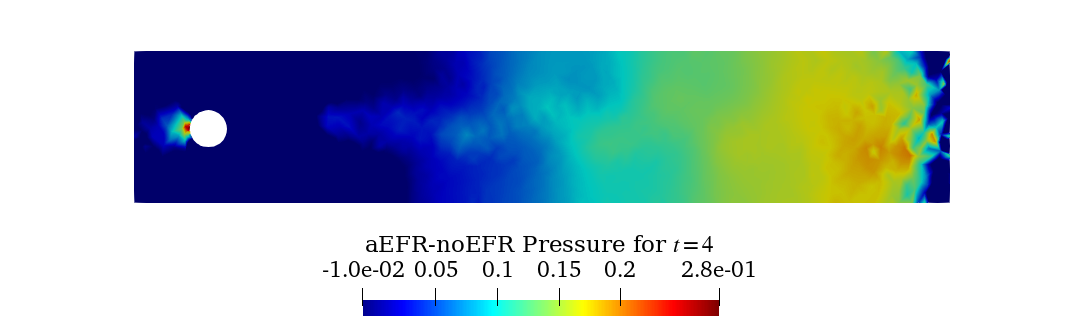}\\
        
\caption{{Uncontrolled snapshot collection. Top plot: aEFR pressure profile for $t=4$. Bottom plots: aEFR-noEFR {(left)} and aEFR-aEFR {(right)} pressure profiles for $t=4$ and $\gamma=10^{-4}$.}}
\label{fig:rem_comp_p}

\end{figure}

\begin{figure}[H]
        \centering \includegraphics[width=0.43\textwidth]{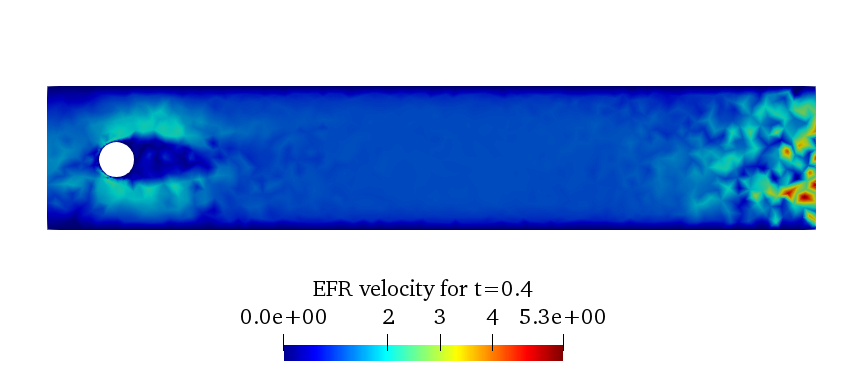} \\
        \centering \includegraphics[width=0.43\textwidth]{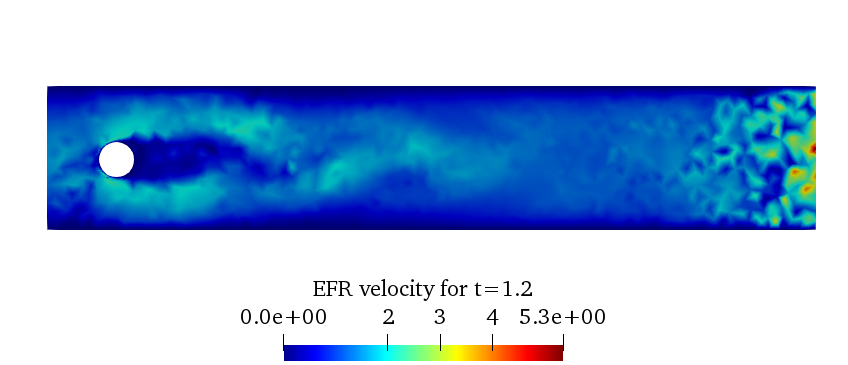}
        \centering \includegraphics[width=0.43\textwidth]{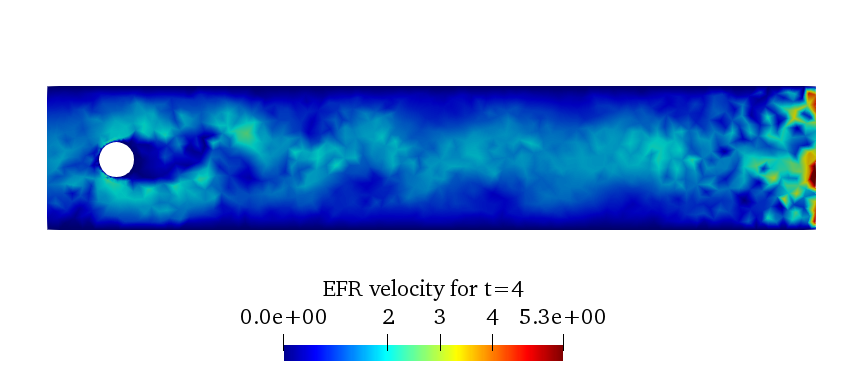}\\
        
\caption{Uncontrolled snapshot collection. Top: EFR snapshot for $t=0.4$. Bottom: EFR snapshot for $t=1.2$ (left) and $t=4$ (right)}

\label{fig:snap_no_control}
\end{figure}

The goal is to recover a \emph{controlled} solution with $r_u = 20$ basis functions for the velocity, and $r_p=r_s=1$ basis functions for the pressure and the supremizer, for $\gamma=10^{-4}$. The data are chosen to be consistent with Experiment 3.
First of all, we compare qualitatively the aEFR-aEFR and aEFR-noEFR solutions to the aEFR FOM solutions in Figure \ref{fig:rem_comp} and \ref{fig:rem_comp_p}, where we show representative solutions for velocity and pressure at $t=4$, respectively. Both fields present unphysical numerical oscillations. The use of uncontrolled snapshots gives more inaccurate results than Experiment 3 (the reader might compare the results to Figures \ref{fig:exp3vcomp} and \ref{fig:exp3pcomp}, respectively). This is not unexpected: the uncontrolled EFR snapshots present spurious oscillations, as shown in Figure \ref{fig:snap_no_control}, where we plot some representative velocity fields for various time instances. For the sake of brevity, we do not show the pressure snapshots, but they feature the same oscillatory behavior of the velocity profiles.

\begin{figure}[H]

  \centering
  \includegraphics[width=0.46\textwidth]{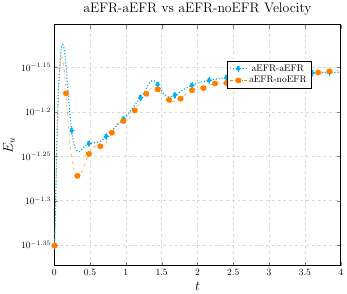} \hfill
  \includegraphics[width=0.45\textwidth]{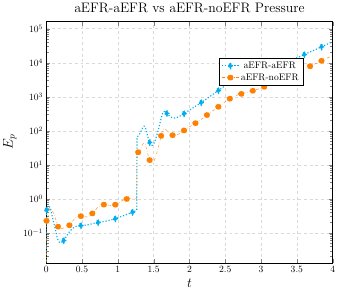}\\
  \caption{Experiment with uncontrolled snapshot collection. Relative errors between aEFR-aEFR (dotted cyan line) and aEFR-noEFR (dashed orange line) controlled velocity {(left)} and pressure {(right)} with $\gamma=10^{-4}$.
  }
  \label{fig:rem_err}
    \end{figure}
 
 \begin{figure}[H]

  \centering
  \includegraphics[width=0.47\textwidth]{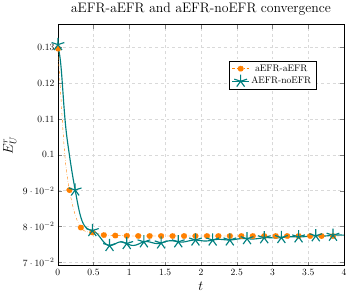}
  \includegraphics[width=0.49\textwidth]{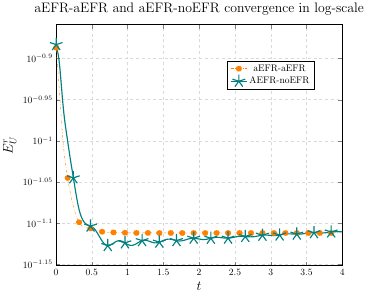}
  \caption{Experiment with uncontrolled snapshot collection. Convergence results for aEFR-noEFR (solid teal line) and aEFR-aEFR (dashed orange line).  The left plot is %proposed 
  {replotted} %in 
  {on a} log-scale on the right.}
  \label{fig:rem_conv}
\end{figure}

These oscillations spoil the ROM reconstruction, as testified by the plot of the relative errors between the controlled FOM reference and the controlled aEFR-aEFR and aEFR-noEFR solutions in Figure \ref{fig:rem_err}: the errors drastically increase when compared to Experiment 3 (most of all the pressure error). Furthermore, the convergence is affected by the choice of the snapshots: both aEFR-aEFR and aEFR-noEFR fail to reach the desired state, and the numerical convergence towards the target is lost.
In conclusion, aEFR-aEFR is competitive and effective when controlled snapshots are used to build the ROM space. If uncontrolled snapshots are used, numerical oscillations arise. We recall that using controlled snapshots is not affecting the computational efficiency of the proposed approach, due to the \emph{a priori} nature of the feedback law.
We believe that tuning $\delta^*$, $\chi$, and $\gamma$ may lead to better results. However, finding the optimal parameters for the uncontrolled and controlled EFR is beyond the scope of this contribution.

\section{Predictive regime}
\label{rem:pred}
    In Experiment 3, we considered the \emph{reconstructive regime}, i.e., we ran the reduced model for the same time window as the training phase. Next, we focus on the predictive capabilities of our approach. For the sake of clarity, we list again the parameters of the simulation: we collect $1000$ equispaced \emph{controlled} snapshots in $[0,T]$, with $T=4$, employing $\nu=10^{-4}$ (i.e., $Re=1000$), $\gamma=10^{-4}$, $C_{\delta}=\sqrt{11}, \delta^*=4.46\cdot 10^{-3}$, and $\chi=5\cdot \Delta t$, with $\Delta t = 4 \cdot 10^{-4}$. However, this time, we simulate the reduced controlled flow until $t=6$, i.e., for $t>T$.
  
    \begin{figure}[H]

  \centering
  \includegraphics[width=0.46\textwidth]{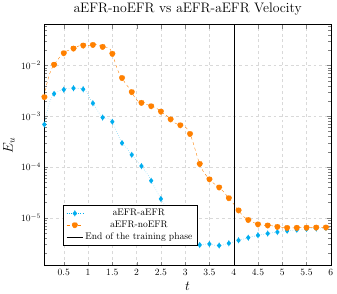} \hfill
  \includegraphics[width=0.45\textwidth]{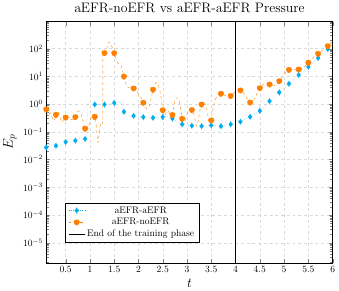}\\
  \caption{Experiment in the predictive regime. Relative errors between aEFR-aEFR (dotted cyan line) and aEFR-noEFR (dashed orange line) controlled velocity {(left)} and pressure {(right)} with $\gamma=10^{-4}$.
  }
  \label{fig:rem_err_pred}
    \end{figure}
\begin{figure}[H]

  \centering
  \includegraphics[width=0.47\textwidth]{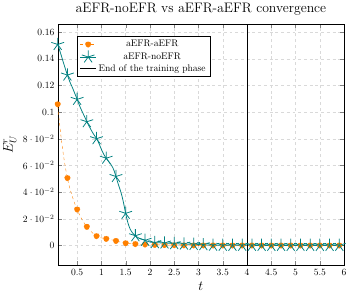}
  \includegraphics[width=0.49\textwidth]{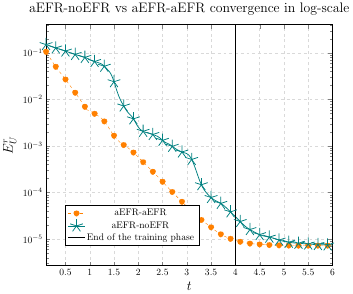}
  \caption{Experiment in the predictive regime. Convergence results for aEFR-noEFR (solid teal line) and aEFR-aEFR (dashed orange line).  The left plot is %proposed 
  {replotted} %in 
  {on a} log-scale on the right.}
  \label{fig:rem_conv_pred}
\end{figure}

 We perform a quantitative analysis in terms of relative errors with respect to the aEFR FOM solution in Figure \ref{fig:rem_err_pred}. This plot shows the reconstructive regime error until $t=4$ (the vertical line represents the time window considered for the POD basis construction) and the predictive regime error until $t=6$. We observe that the velocity field attains a degree of accuracy of the order $10^{-5}$. Moreover, by the end of the simulation, the two approaches coincide. This is not unexpected: the flow has been completely stabilized by $t=6$, and aEFR-aEFR and aEFR-noEFR are comparable. We note, however, that aEFR-aEFR is more accurate until $t=5$. 
In the right plot, we observe that the pressure error significantly increases. As discussed in Experiment 3, both algorithms do not recover the pressure accurately, but aEFR-aEFR is more accurate than aEFR-noEFR until $t=5.5$.
    For completeness, we show the tracking error in the predictive regime in Figure \ref{fig:rem_conv_pred}. We see that both approaches yield $E_U^r$ of the order of $10^{-5}$, despite no information of the controlled dynamic being collected in the POD phase for $t > 4$. For brevity, we do not show the qualitative results for $t>4$: the aEFR-aEFR and aEFR-noEFR velocity fields are indistinguishable from each other and from the desired state represented in Figure \ref{fig:U1}.  
 
\section{aEFR-aEFR reduction in time and in the parameter $\gamma$}
\label{rem:gammapar}
In Experiment 3, we consider the time as the \emph{only parameter}. To propose the aEFR-aEFR strategy as a way to predict control actions for several values of $\gamma$, we report some results for $\gamma$ as an additional parameter. In this numerical investigation, we use a nested-POD (n-POD) strategy, to avoid the computational burden of a standard POD approach. 
The n-POD performs (i) a first compression on each parametric trajectory considered in the building phase, retaining a first set of modes, and (ii) a final compression on these modes scaled by their singular values. More details on this approach, which goes under several names in the ROM community, can be found in \cite{audouze2009reduced,ballarin2016fast,brands2016reduced, himpe2018hierarchical}.
For the sake of clarity, we list again the parameters of the simulation: we collect $1000$ equispaced \emph{controlled} snapshots in time in the interval $[0,T]$, with $T=4$, employing $\nu=10^{-4}$ (i.e., $Re=1000$), $C_{\delta}=\sqrt{11}, \delta^*=4.46\cdot 10^{-3}$, and $\chi=5\cdot \Delta t$, with $\Delta t = 4 \cdot 10^{-4}$ for several values of $\gamma$. 
We consider different time trajectories given by 15 equispaced values of $\gamma$ in $[10^{-1}, 10^{-5}]$. In step (i), we perform a first compression retaining 60 modes for each $\gamma$ we considered in the offline phase. Then, in step (ii), the modes are further compressed for a reduced space with dimensions $r_u=20$ and $r_p=r_s=1$ for velocity, pressure, and supremizer, respectively. The information retained by the n-POD is analogous to the one of the POD approach for time reduction in Experiment 3, reported in Table \ref{tab:E3energy}. To test the accuracy of the n-POD, we analyze the relative errors $E_u$ and $E_p$ averaged in time for $r_u \in \{5,10,15,20\}$ and $\gamma \in \{5\cdot 10^{-2},10^{-2},10^{-3},10^{-4},5 \cdot 10^{-5}\}$. 
We recall that we fixed $r_p=r_s=1$ since the mode retains 99\% of the pressure information.

\begin{table}[H]
\caption{\reviewerB{Average relative errors for aEFR-aEFR for parametric $\gamma$. The velocity and pressure errors are shown with respect to the basis dimension $r_u$ and the parameter $\gamma$.}}
\label{tab:rom_table_with_gamma}
\resizebox{\textwidth}{!}{
\begin{tabular}{cc|cccc||cccc|}
\cline{3-10}
                                                  &           & \multicolumn{4}{c||}{Average $E_u$}                                                                           & \multicolumn{4}{c|}{Average $E_p$}                                                                           \\ \cline{2-10} 
\multicolumn{1}{c|}{}                             & $\gamma$  & \multicolumn{1}{c|}{$r_u = 5$} & \multicolumn{1}{c|}{$r_u =10$} & \multicolumn{1}{c|}{$r_u =15$} & $r_u =20$ & \multicolumn{1}{c|}{$r_u = 5$} & \multicolumn{1}{c|}{$r_u =10$} & \multicolumn{1}{c|}{$r_u =15$} & $r_u =20$ \\ \hline
\multicolumn{1}{|c|}{\multirow{5}{*}{aEFR-aEFR}} & $5 \cdot 10^{-2}$ & \multicolumn{1}{c|}{2.20e-3}          & \multicolumn{1}{c|}{1.39e-3}          & \multicolumn{1}{c|}{1.00e-3}          &     7.29e-4      & \multicolumn{1}{c|}{8.09e-4}          & \multicolumn{1}{c|}{7.18e-1}          & \multicolumn{1}{c|}{4.25e-1}          &       3.37e-1    \\ \cline{2-10} 
\multicolumn{1}{|c|}{}                            & $10^{-2}$ & \multicolumn{1}{c|}{2.17e-3}          & \multicolumn{1}{c|}{1.33e-3}          & \multicolumn{1}{c|}{9.48e-4}          &    7.53e-4       & \multicolumn{1}{c|}{4.54e-1}          & \multicolumn{1}{c|}{7.09e-1}          & \multicolumn{1}{c|}{4.18e-1}          &      3.35e-1     \\ \cline{2-10} 
\multicolumn{1}{|c|}{}                            & $10^{-3}$ & \multicolumn{1}{c|}{2.16e-3}          & \multicolumn{1}{c|}{1.32e-3}          & \multicolumn{1}{c|}{9.46e-4}          &    9.36e-4       & \multicolumn{1}{c|}{4.50e-1}          & \multicolumn{1}{c|}{7.07e-1}          & \multicolumn{1}{c|}{4.17e-1}          &      3.34e-1     \\ \cline{2-10} 
\multicolumn{1}{|c|}{}                            & $10^{-4}$ & \multicolumn{1}{c|}{2.14e-3}          & \multicolumn{1}{c|}{1.31e-3}          & \multicolumn{1}{c|}{9.35e-4}          &    7.40e-4       & \multicolumn{1}{c|}{4.54e-1}          & \multicolumn{1}{c|}{7.07e-1}          & \multicolumn{1}{c|}{4.17e-1}         &      3.33e-1     \\ \cline{2-10} 
\multicolumn{1}{|c|}{}                            & $5 \cdot 10^{-5}$ & \multicolumn{1}{c|}{2.17e-3}          & \multicolumn{1}{c|}{1.34e-3}          & \multicolumn{1}{c|}{9.34e-4}          &      7.41e-4    & \multicolumn{1}{c|}{4.55e-1}          & \multicolumn{1}{c|}{7.08e-1}          & \multicolumn{1}{c|}{4.17e-1}         &        3.33e-1   \\ \hline

\end{tabular}
}
\end{table}

 \begin{figure}[H]

  \centering
  \includegraphics[width=0.49\textwidth]{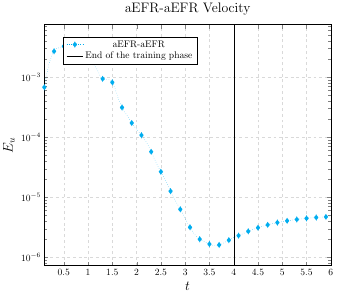} \hfill
  \includegraphics[width=0.49\textwidth]{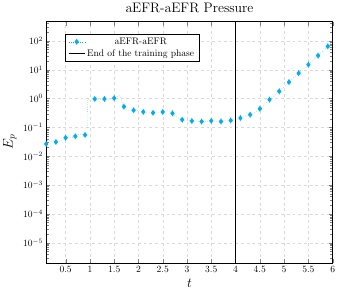}\\
  \caption{Experiment with reduction both in time and for $\gamma$. Relative errors between aEFR-aEFR controlled velocity {(left)} and pressure {(right)} with $\gamma=10^{-4}$.
  }
  \label{fig:rem_err_pred_par}
    \end{figure}
 
 \begin{figure}[H]

  \centering
  \includegraphics[width=0.49\textwidth]{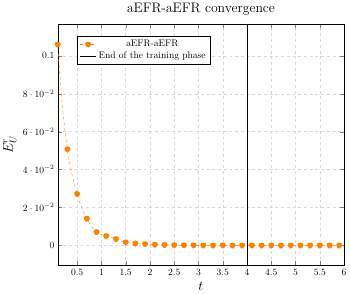}
  \includegraphics[width=0.49\textwidth]{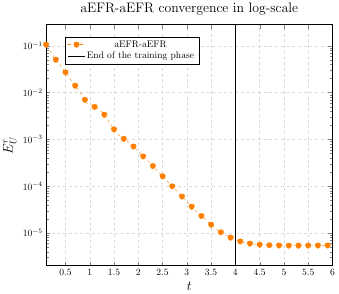}
  \caption{Experiment with reduction both in time and for $\gamma$. Convergence results for aEFR-aEFR for $\gamma=10^{-4}$.  The left plot is %proposed 
  {replotted} %in 
  {on a} log-scale on the right.}
  \label{fig:rem_conv_pred_par}
\end{figure}

The values of $\gamma$ investigated online were not considered in the building phase. We obtain good results in terms of extrapolation of the controlled velocity for all the $\gamma$ and $r_u$ values, with errors ranging between $2 \cdot 10^{-3}$ and $7 \cdot 10^{-4}$. As usual, the pressure error is larger and of the order of $10^{-1}$, in analogy with what was already observed in Experiment 3.
We conclude that the n-POD approach can predict the controlled velocity fields in the reconstructive regime with high accuracy.

We now fix $\gamma=10^{-4}$ as a representative parameter. For completeness, in Figure \ref{fig:rem_err_pred_par}, we plot the relative errors for aEFR-aEFR velocity and pressure fields compared to aEFR FOM velocity both in the reconstructive and the predictive regime. The results are comparable to the results obtained with reduction only in time. The same holds for the tracking error $E_U^r$ plotted in Figure \ref{fig:rem_conv_pred_par}: the convergence is similar to the results obtained in Experiment 3 and %Remark 
\ref{rem:pred}.

\end{document}